\documentclass[a4paper, 11pt, leqno]{amsart}

\usepackage[dvipsnames]{xcolor}

\usepackage{iftex}
\ifPDFTeX 
  \usepackage[utf8]{inputenc}
  \usepackage[T1]{fontenc}
\else 
  \usepackage{fontspec}
\fi
\usepackage{lmodern}

\usepackage{geometry}
\usepackage[textwidth=28mm]{todonotes}
\usepackage{mathtools}
\usepackage{yfonts}
\usepackage[foot]{amsaddr}

\usepackage{subcaption} 

\usepackage{amssymb} 

\usepackage[pdfencoding=auto]{hyperref}
\hypersetup{
    colorlinks,
    citecolor=green,
    filecolor=black,
    linkcolor=blue,
    urlcolor=black
}

\usepackage[english]{babel}

\usepackage{dsfont}
\usepackage{bm} 
\usepackage{comment}
\usepackage{enumerate}
\usepackage[shortlabels]{enumitem}
\usepackage{tikz}
\usepackage{pgfplots}
\pgfplotsset{compat=1.15}
\usetikzlibrary{arrows}
\usepackage{algorithm}
\usepackage{algorithmic}

\numberwithin{equation}{section}

\usepackage{tabularx}

\newcommand{\nocontentsline}[3]{}
\let\origcontentsline\addcontentsline
\newcommand\stoptoc{\let\addcontentsline\nocontentsline}
\newcommand\resumetoc{\let\addcontentsline\origcontentsline}

\theoremstyle{plain}
\newtheorem{thm}{Theorem}[section]

\newtheorem{rem}[thm]{Remark}

\theoremstyle{definition} 
 
\theoremstyle{remark}

\newcommand{\R}{\mathbb{R}}

\newcommand{\eps}{\varepsilon}
\newcommand{\Q}{\mathcal{Q}}
\newcommand{\tA}{\tilde{A}}
\newcommand{\ta}{\tilde{a}}
\newcommand{\tb}{\tilde{b}}
\newcommand{\tc}{\tilde{c}}
\newcommand{\bA}{\bar{A}}
\newcommand{\ba}{\bar{a}}
\newcommand{\bb}{\bar{b}}
\newcommand{\bc}{\bar{c}}

\DeclareMathOperator{\diverg}{div}

\definecolor{codegreen}{rgb}{0,0.6,0}
\definecolor{codegray}{rgb}{0.5,0.5,0.5}
\definecolor{codepurple}{rgb}{0.58,0,0.82}
\definecolor{backcolour}{rgb}{0.95,0.95,0.92}

\newcounter{review}

\makeatletter

\newcommand \listreviewname{List of Reviews}
\newcommand \listofreviews{\section*{\listreviewname} \@starttoc{tor}}
\makeatother

\definecolor{revisionColourOne}{RGB}{180, 0, 0}
\definecolor{revisionColourTwo}{RGB}{0, 0, 180}

\author{Domenico Caparello}
\author{Tommaso Tenna}
\address[Domenico Caparello]{Université Côte d’Azur, CNRS, LJAD, Parc Valrose, F-06108 Nice, France -- Department of Mathematics and Computer Science, University of Ferrara, Italy}
\address[Tommaso Tenna]{Université Côte d’Azur, CNRS, LJAD, Parc Valrose, F-06108 Nice, France -- Dipartimento di Matematica, Sapienza Università di Roma, P.le Aldo Moro 5, 00185 Rome, Italy}

\title[A High-Order Coupling Method for the ES-BGK model]{A Coupled IMEX Domain Decomposition Method for High-Order Time Integration of the ES-BGK Model of the Boltzmann Equation}
\date{}

\begin{document}


\maketitle

\begin{abstract}
In this paper, we propose a high-order domain decomposition method for the ES-BGK model of the Boltzmann equation, which dynamically detects regions of equilibrium and non-equilibrium. Our implementation automatically switches between Euler equations in regions where the fluid is at equilibrium, and the ES-BGK model elsewhere. The main challenge addressed in this work is the development of a coupled strategy between the macroscopic and the kinetic solvers, which preserves the overall temporal order of accuracy of the scheme. A coupled IMEX method is introduced across decomposed subdomains and solvers. This approach is based on a coupled IMEX method and allows high accuracy and computational efficiency. Several numerical simulations in two space dimensions are performed, in order to validate the robustness of our approach and the expected temporal high-order convergence.

\end{abstract}

\section{Introduction}
Numerical strategies for kinetic equations have been widely investigated over the last decades. These models describe the behavior of particles of rarefied gas undergoing instantaneous collisions and appear in a variety of systems, ranging from astrophysics and aerospace applications \cite{boscheridimarco2021, dimarco2018} to plasma physics phenomena \cite{dimarcomieussens2014, medaglia2023, albi2025}. The seminal Boltzmann Equation represents the cornerstone to describe the dynamics of a gas which is not in thermodynamical equilibrium, through a distribution function which depends on space, velocity and time \cite{cercignani1988}. This means that for its nature, this statistical distribution function lives in a seven dimensional space (three dimensions both for the physical space and the velocity space, plus the time variable). Simplifications of such systems in a bi-dimensional framework still have prohibitive computational cost. Several simplified kinetic models have been proposed, such as the relaxation approximation of the Boltzmann equation known as BGK equation \cite{BGK1954}. The most problematic aspect of the BGK equation is that, in the macroscopic limit, it produces the incorrect coefficients for the Navier Stokes equation. To this aim, a corrected version of this model has been proposed, called ellipsoidal statistical BGK equation (ES-BGK). This model has been introduced in \cite{holway1965} and it satisfies the correct conservation laws, yielding the Navier-Stokes approximation via the Chapman-Enskog expansion, see also \cite{andries2000}.\\
One of the most challenging aspects of the numerical approximation of kinetic equations is the ``curse of dimensionality''. These models combine linear transport with nonlinear collision operators in a high dimensional kinetic space, making them significantly more expensive than their fluid-dynamical counterparts from a computational point of view. The other issue is represented by the presence of multiple scales, which requires proper numerical methods to treat the stiff dynamics arising in problems with different regimes \cite{filbetjin2010, jin2022}. The aim of this paper is the construction of an efficient high-order numerical solver for the multiscale ES-BGK equation, developing a high-order in time domain decomposition method, which is able to rapidly compute an accurate solution for all regimes. In particular, the kinetic solver is based on a finite volume scheme combined the IMEX method proposed by Filbet and Jin in \cite{fil-jin-2011}. For the sake of completeness, we mention different techniques developed for the approximation of the ES-BGK equation, e.g. random particle methods \cite{andries2002}, semi-Lagrangian schemes \cite{boscarino2025ESBGK} or micro-macro flux vector splitting \cite{rossmanith2025}.\\
The motivation behind domain decomposition methods is the use of macroscopic equations (e.g. \textit{Euler equations}) in regions where the fluid is close to equilibrium, while restricting kinetic solvers for regions where non-equilibrium effects are dominant. Indeed, Euler equations are not valid when the fluid is out of equilibrium; therefore, kinetic models are required in those regions. These domain decomposition methods blend kinetic and fluid models, by partitioning the domain and using switching criteria based on physical indicators \cite{tiwari1998, tiwariklar1998, dimarcopareschi2008, dimarcomieussens2014}. The core of this strategy lies in replacing the full kinetic equation with its macroscopic closure, thereby reducing the computational cost associated with the numerical approximation. Several papers have developed methods based on this idea, introducing many different types of regime indicators, see \cite{tiwari1998, kol-ars-ari-fro-zab-2007, ala-pup-2011, ala-pup-2012, deg-dim-2012, filbetrey2015, caparello2025}. For the sake of completeness, we mention other recent works exploiting this hybrid domain decomposition strategy for different kinetic equations, such as the Vlasov-BGK or the Vlasov-Poisson-BGK equations \cite{dimarcomieussens2014, lai-2023, lai-rey-2023, laidin2025}.\\
However, in regimes far from equilibrium, a robust and efficient solver for the resolution of the full ES-BGK equation is necessary. Different strategies have been implemented to avoid explicit time integration of general stiff kinetic equations, e.g. IMEX schemes \cite{pareschirusso2005, pieraccini2007, filbetjin2010, dimarcopareschi2012, hu2017,  puppo2019, boscheridimarco2021, boscarinopareschirusso2025} or projective integration methods \cite{lafitte2012, melis2019, bailo2022, tenna2025}. These schemes allow the use of bigger time steps, without any restriction due to the stiffness of the system. The most challenging aspect is related to the combination of such high-order schemes with domain decomposition methods. Indeed, inaccurate coupling strategies between the kinetic and the fluid models may lead to low-order approximation, even using high-order time integrators.\\
In this paper, we propose a numerical method which exploits the advantages of a domain decomposition method for the ES-BGK model to detect regions of equilibrium, employing a new coupled strategy based on IMEX schemes to reach higher order approximations. The approach we propose is not problem-dependent and it can be easily extended to other types of kinetic equations, leading to a greater speedup and enhancing performance of the method, along with a reduced computational cost. Several numerical examples in two space dimensions are included to show the quality of the approximation and the desired properties of the scheme.
\subsubsection*{Structure of the paper.} The rest of the paper is organized as follows. Section \ref{Section_ESBGK} is devoted to the introduction of the ES-BGK model and its Chapman-Enskog expansion. The numerical strategy is presented in Section \ref{Section_Numerical}, where domain decomposition method is described and the coupled IMEX scheme is constructed. In Section \ref{Section_Simulations} we perform several numerical tests in two space dimensions to validate the accuracy of the proposed scheme and its ability to reproduce stiff phenomena.

\section{The ES-BGK Equation}
\label{Section_ESBGK}
Let us consider a given non-negative initial distribution function $f_0(x,v)$ for $x \in \Omega \subset \R^{d_x}$ and $v \in \R^{d_v}$. The open set $\Omega$ is a bounded Lipschitz-continuous domain of $\R^{d_x}$, and appropriate boundary conditions will be described in the following sections. Then, the particle distribution function $f^\eps = f^\eps(t,x,v)$, for $t \geq 0$, evolves according to the ES-BGK equation which, in absence of external forces, has the form
\begin{equation} 
\label{eqCollision}
\left\{ \begin{aligned}
  & \frac{\partial f^\eps}{\partial t} + v \cdot \nabla_x f^\eps \,=\, \frac{1}{\eps}\;\Q(f^\eps), 
  \\
  &  \;
  \\
  & f^\eps(0, x, v) = f_{0}(x,v),
\end{aligned} \right.
\end{equation}
The parameter $\eps >0$ is the \textit{Knudsen number}, defined as the ratio between the mean free path of particles before a collision and the physical length scale of observation. This quantity governs the frequency of collisions, giving an indication of the rarefaction of the gas: if $\eps \ll 1$ the evolution of the particle system is governed by collisions and the gas is \emph{dense} (fluid regime); viceversa, if $\eps \sim 1$ the evolution of the particle system is governed by the advection term and the gas can be considered \emph{rarefied} \cite{cercignani1988}.
In the following we will denote $f$ without the superscript $\eps$ whenever the dependence on $\eps$ is understood.\\

We assume that the collision operator fulfils three fundamental assumptions
\begin{enumerate}[label=\textbf{(H$\bm{_\arabic{*}}$)}, ref=\textbf{(H$\bm{_\arabic{*}}$)}]

\item \label{hypConservations} 
It preserves mass, momentum and kinetic energy. This properties can be, respectively, expressed as
  \begin{equation*}
    \int_{\R^{d_v}} \Q(f)(v) \, dv = 0, \quad  \int_{\R^{d_v}} \Q(f)(v) \, v \, dv = 0, \quad \int_{\R^{d_v}} \Q(f)(v) \, |v|^2 \, dv = 0.
  \end{equation*}
 
\item \label{hypEntropy}
  The Boltzmann entropy is dissipated (H-theorem) \cite{cercignani1988, cercignani1994}
  \begin{equation*}
		\int_{\R^{d_v}} \Q(f)(v) \, \log(f)(v) \, dv \, \leq \, 0.
	\end{equation*}
	
\item \label{hypEquilib} 
The equilibria are given by Maxwellian distributions
\begin{equation*}
\label{eq::Maxwellian}
\mathcal M_{\rho, \bm u, T} := \frac{\rho}{(2 \pi T)^{d_v/2}} \exp \left ( - \frac{|v-\bm u|^2}{2 T} \right ),
\end{equation*}
where the \emph{density}, \emph{velocity} and	\emph{temperature} of the gas $\rho$, $\bm u$ and $T$  are computed as suitable moments of the distribution function $f$,  
\begin{equation*}
\rho = \int_{\R^{d_v}}f(v)\,dv, \quad \bm{u} = \frac{1}{\rho}\int_{\R^{d_v}}v f(v) \, dv, 
\quad T = \frac{1}{d_v\,\rho} \int_{\R^{d_v}} \vert \bm{u} - v \vert^2 f(v) \,dv.
\end{equation*}

\end{enumerate}  
According to assumptions \ref{hypEntropy}-\ref{hypEquilib},when $\eps \to 0$, the distribution $f^\eps$ converges (at least formally) to a Maxwellian distribution, whose moments are solution to the compressible Euler system
\begin{equation}
	\label{eqHydroClosedEuler}
	\left\{ \begin{aligned}
	  & \partial_t \rho  + \diverg_x  (\rho \, \bm{u} ) \,=\, 0, 
	  \\
	  \,
	  \\
	  & \partial_t(\rho \, \bm{u} ) + \diverg_x  \left(\rho \, \bm{u} \otimes \bm{u} \,+\, \rho \, T  \,{\rm\bf I}\right) \, =\, \bm{0}_{\R^{3}}, 
	  \\
	  \,
	  \\
	  & \partial_t E + \diverg_x \left ( \bm{u} \left ( E +\rho \, T\right )  \right ) \,=\, 0,
	\end{aligned} \right.
\end{equation}
where $\mathbf{I}$ is the identity matrix and $E$ represents the total energy $$E=\frac{1}{2}\rho |\mathbf{u}|^2 + \frac{d_v}{2}\rho T.$$ The pressure $P$ is related to the internal energy by $$P=\left(\gamma-1\right)\left(E-\frac{1}{2}\rho |\mathbf{u}|^2 \right),$$
where $\gamma=(d_v+2)/d_v$ is the polytropic constant, thus yielding $P=\rho T$.
This limit is typically referred to as an \emph{hydrodynamic limit}. We refer to \cite{golse2005, saintraymond2009} for a complete review of hydrodynamic limits of Boltzmann-type equations.\\
In the BGK equation, the collision operator consists in a relaxation operator towards the Maxwellian equilibrium state, namely
\begin{equation}
\label{eq::Collision_BGK}
    \mathcal{Q}(f) = \frac{\tau}{\eps} \left( \mathcal{M}_{\rho, \bm u, T} [f] - f \right),
\end{equation}
where $\mathcal{M}_{\rho, \bm u, T}$ is given by \ref{eq::Maxwellian} and $\tau>0$, generally depending on the macroscopic moments \cite{dimarcopareschi2014}. For the sake of simplicity, in the following we will avoid the subscripts in the expression of the Maxwellian distribution function $\mathcal{M}$, whenever the dependence on macroscopic moments is not needed.\\
As $\eps \to 0$, this model has the correct Euler limit \eqref{eqHydroClosedEuler} but, at the first-order correction in $\eps$, the coefficients obtained at the Navier-Stokes level are incorrect. Indeed, despite this model preserves mass, momentum and kinetic energy and provides the correct dissipation of the Boltzmann entropy, it produces the incorrect Prandtl number, as described in Section \ref{Section::CE}.\\
Several kinetic models have been proposed as alternative to the BGK equation, in order to give the correct transport coefficients at the Navier-Stokes level, such as the ellipsoidal statistical (ES-BGK) model \cite{holway1965, andries2000}. It consists in a relaxation operator, where the Maxwellian distribution function $\mathcal{M}[f]$ in \eqref{eq::Collision_BGK} is replaced by an anisotropic Gaussian function $\mathcal{G}[f]$, defined by
\begin{equation}
    \mathcal{G}[f] = \frac{\rho}{\sqrt{\det(2 \pi \mathcal{T})}} \exp \left ( - \frac{(v-\bm u)^\intercal\mathcal{T}^{-1}(v-\bm u)}{2} \right ),
\end{equation}
where the temperature tensor $\mathcal{T}$ is defined as
\begin{equation}
    \mathcal{T} = \left[(1-\beta) T \bm I + \beta \mathbf{\Theta} \right],
\end{equation}
and the opposite of the stress tensor $\mathbf{\Theta}$ is given by
\begin{equation}
    \mathbf{\Theta} = \frac{1}{\rho} \int_{\R^{d_v}} (v - \bm u)\otimes (v-\bm u) f(v) dv.
\end{equation}
Therefore now the collision operator reads as follows
\begin{equation}
\label{eq::Collision_ESBGK}
    \mathcal{Q}(f) = \frac{\tau}{\eps} \left( \mathcal{G}[f] - f \right).
\end{equation}

\subsection{Chapman-Enskog expansion}
\label{Section::CE}
 Let us introduce a perturbation $h$ around the Maxwellian equilibrium $\mathcal{M}$ defined as
 \begin{equation}\label{eqCE_distfun}
     f=\mathcal{M}+\eps h.
 \end{equation}
 In the asymptotic limit $\varepsilon\to 0$, the distribution function $\mathcal{G}$ becomes isotropic, and it converges to the classical Maxwellian $\mathcal{M}$ \cite{fil-jin-2011}. Thus, at order zeroth in the Chapman-Enskog expansion (i.e. $f=\mathcal{M}$), one recovers the Euler equations.\\
 Let us now focus on the first order expansion.\\
 Similarly to Equation \eqref{eqCE_distfun}, the operator $\mathcal{G}$ and $\mathbf{\Theta}$ can be expanded as \cite{fil-jin-2011}
 \begin{equation}
     \begin{cases}
         \mathcal{G} = \mathcal{M} + \varepsilon g, \\
         \mathbf{\Theta} = T \bf{I} + \varepsilon \mathbf{\Theta}_1.\\
     \end{cases}
 \end{equation}
 For the sake of simplicity, we will restrict to the case $d_v=3$.
 The quantity $h$ can be expressed as \cite{fil-jin-2011}
 \begin{equation}
     h = g - \frac{1}{\tau}\frac{\partial \mathcal{M}}{\partial t} - \frac{1}{\tau} v \cdot \nabla_x \mathcal{M}.
 \end{equation}
 From \cite{levermore1996, golse2005, saintraymond2009} we rewrite the following equation
 \begin{equation}
 \label{eq::evolution_Maxwellian}
     \frac{\partial \mathcal{M}}{\partial t} + v \cdot \nabla_x = \mathcal{M}\left[\frac{1}{T}\left(c\otimes c - \frac{1}{3}|c|^2\bf{I}\right):\nabla_x \bm u + \left(\frac{|c^2|}{2T}-\frac{5}{2}\right)c\cdot \frac{\nabla_x T}{T}\right],
 \end{equation}
 where $c=v-\bm u$.\\
 By introducing the following Sonine polynomials \cite{levermore1996, golse2005, lev-mor-nad-1998} 
 \begin{equation}
    \begin{cases}
        \mathbf{A}(v-\bm u) = \displaystyle \frac{v- \bm u}{\sqrt{T}}\otimes \frac{v-\bm u}{\sqrt{T}}-\frac{1}{3}\frac{|v-\bm u|^2}{T}\bf{I},\\[10pt]
        \mathbf{B}(v-\bm u) = \displaystyle \left[\frac{|v-\bm u|^2}{2T}-\frac{5}{2}\right]\frac{v-\bm u}{\sqrt{T}},\\
    \end{cases}
 \end{equation}
 it is possible to rewrite Equation \eqref{eq::evolution_Maxwellian} as
 \begin{equation}
     \frac{\partial \mathcal{M}}{\partial t} + v \cdot \nabla_x = \mathcal{M}\left[\mathbf{A}:\nabla_x \bm u + \mathbf{B}\cdot\frac{\nabla_x T}{\sqrt{T}}\right].
 \end{equation}
 From \cite{fil-jin-2011} the perturbation of the anisotropic Gaussian $\mathcal{G}$
 \begin{equation}
     g = \frac{\beta}{2T^2}\mathcal{M}\left(v-\bm u\right)^\intercal \mathbf{\Theta}_1 \left(v-\bm u\right).
 \end{equation}
 Then combining these last two quantities, one gets that
 \begin{equation}
     h = \frac{\beta}{2T^2}\mathcal{M}\left(v-\bm u\right)^\intercal \mathbf{\Theta}_1 \left(v-\bm u\right) - \frac{\mathcal{M}}{\tau}\left[\mathbf{A}:\nabla_x \bm u + \mathbf{B}\cdot\frac{\nabla_x T}{\sqrt{T}}\right].
 \end{equation}
 Hence, the first order correction to the fluxes in the formal conservation law yields the compressible Navier-Stokes system with $\mathcal{O}(\eps)$ dissipation terms \cite{fil-jin-2011}
 \begin{equation}
     \begin{cases}
      & \partial_t \rho  + \diverg_x  (\rho \, \bm{u} ) \,=\ 0,\\ 
	  & \partial_t(\rho \, \bm{u} ) + \diverg_x  \left(\rho \, \bm{u} \otimes \bm{u} \,+\, \rho \, T  \,{\rm\bf I}\right) \, = \eps \text{div}_x (\mu \mathbf{D}(\bm u)),\\
	  & \partial_t E + \diverg_x \left ( \bm{u} \left ( E +\rho \, T\right )  \right ) \,=\, \eps \text{div}_x \left (\mu \mathbf{D}(\bm u)\bm u + \kappa \nabla_x T \right),
     \end{cases}
 \end{equation}
 where the scalar quantities $\mu$ and $\kappa$ are respectively the \textit{viscosity} and the \textit{thermal conductivity}, while $\mathbf{D}(\bm{u})=\left(\nabla_x \bm{u}\right) + \left( \nabla_x \bm{u}\right)^\intercal - \frac{2}{3}\text{div}_x \left(\bm{u}\right)$ as defined in \cite{golse2005}. In the ES-BGK model, we have \cite{fil-jin-2011} 
 \begin{equation*}
     \mu = \frac{1}{1-\beta} \frac{\rho T}{\tau} \qquad \text{and} \qquad \kappa = \frac{5}{2}\frac{\rho T}{\tau}.
 \end{equation*}
The Prandtl number is defined as \cite{fil-jin-2011}
\begin{equation}
    \text{Pr} = \frac{\gamma}{\gamma-1} \frac{\mu}{\kappa},
\end{equation}
where the polytropic constant $\gamma = (d_v + 2)/d_v$ represents the ratio between specific heat at constant pressure and at constant volume. Then, for the ES-BGK model, it satisfies
\begin{equation}
    \text{Pr} = \frac{1}{1-\beta},
\end{equation}
where $\beta\in\left[-1/2,1 \right)$. In the case of BGK equation $\beta=0$ and one gets the incorrect Prandtl number $\text{Pr}=1$. Since the majority of the gases has $\text{Pr}<1$ (in particular the hard-sphere model of monoatomic gas in the Boltzmann equation leads to a Prandtl number close to $2/3$), this motivates the introduction of the ES-BGK equation.\\

\section{Numerical Method}
\label{Section_Numerical}
\subsection{Domain Decomposition Method}
 The approach used in this paper is inspired by the foundational works of Levermore, Morokoff, and Nadiga \cite{lev-mor-nad-1998} and Filbet and Rey \cite{filbetrey2015}. In particular, we combine a domain decomposition (also called \textit{hybrid}) strategy with a temporal high-order coupling, resulting in a very efficient fully explicit solver.\\ The proposed hybrid method offers two main advantages: it guarantees high accuracy in multiscale regimes compared to classical fluid solvers, since in moderate and highly non-equilibrium areas (where Euler equations are no longer valid) the method employs the full kinetic ES-BGK model, providing accurate resolution of shocks and boundary layers; at the same time, by dynamically updating equilibrium regions, it confines the expensive resolution of the full kinetic equation to small regions, relying on the efficient fluid solver elsewhere. This approach allows to obtain a significant speedup for the approximation of the ES-BGK equation without compromising accuracy and physical fidelity.\\
 Several criteria may be introduced to decompose the domain. In \cite{boy-che-can-1995}, a macroscopic criterion based on the value of the Knudsen number $\eps$ is used to switch between the hydrodynamic description and the kinetic one. We recall that the Knudsen number is defined as the ratio between the mean free path of particles and the characteristic physical length scale of the system. By fixing a threshold for the Knudsen number, the regime indicator is easily constructed.\\
 Following \cite{filbetrey2015, filbetxiong2018, caparello2025}, in this paper two different criteria are employed: one to identify when the hydrodynamic description breaks down, needing a transition to the kinetic framework, and another to determine when it is appropriate to switch from the kinetic regime back to the hydrodynamic description. The key tool for introducing these criteria is the Chapman-Enskog expansion, see Section \ref{Section::CE}.

 \subsubsection{Computation of indicators}
 In this subsection, we define the indicators which will be used in the method. Let us introduce the moment realizability matrix, which captures the internal structure of the distribution function $f$ through its collisional invariants. This matrix encodes how the macroscopic fields defined in \ref{hypEquilib} deviate from those compatible with a realizable non-negative distribution function $f$. Heuristically, it reflects whether the computed moments lie within the physically admissible set generated by $f$. Following the ideas of \cite{lev-mor-nad-1998, filbetrey2015}, the moment realizability matrix provides a natural indicator of hydrodynamic breakdown: when the matrix is not compatible anymore with a positive distribution function, the macroscopic closure becomes unreliable and $f$ is no longer close to equilibrium.\\
 To this aim, let us introduce the following vector
 \begin{equation}
     \mathbf{m} =\displaystyle
     \begin{pmatrix}
      \displaystyle1\\
      \displaystyle\frac{v- \bm u}{\sqrt{T}}\\
      \displaystyle\frac{|v-\bm u|^2}{2T}-\frac{5}{2}
     \end{pmatrix},
 \end{equation}
 and let us define the following moment realizability matrix $M$
 \begin{equation}
 \label{eq:MRM}
     \mathbf{M} = \frac{1}{\rho}\int (\mathbf{m} \otimes \mathbf{m}) f dv.
 \end{equation}
 The integral such defined is a positive semi-definite matrix for any non-negative (and not identically zero a.e.) distribution function \cite{tiw-2000}. The \emph{moment realizability criterion} is based on the following idea. By computing this matrix for different closures of the distribution function $f$ (such as Euler or Navier-Stokes) and by checking whether the matrix $\mathbf{M}$ is positive semi-definite, we can determine whether the closure is physically acceptable or not.\\
 By explicitly computing \eqref{eq:MRM}, we obtain the expression
 \begin{equation}
     \mathbf{M} =
     \begin{pmatrix}
      1  & \mathbf{0}^\intercal            & -1\\
      \mathbf{0}  & \bf{I}+\bf{\bar{A}}    & \bf{\bar{B}}\\
      -1 & \bf{\bar{B}}^\intercal & \bar{C}
     \end{pmatrix},
 \end{equation}
 where we have introduced the following quantities
 \begin{equation}
     \begin{cases}
         \bar{\mathbf{A}} = \displaystyle \frac{1}{\rho}\int \mathbf{A}(v-\bm u) f dv,\\[8pt]
         \bar{\mathbf{B}} = \displaystyle \frac{1}{\rho}\int \mathbf{B}(v-\bm u) f dv,\\[8pt]
         \bar{C} = \displaystyle \frac{1}{\rho}\int\left(\frac{|v-\bm u|^2}{2T}-\frac{5}{2}\right)^2 f dv.\\
     \end{cases}
 \end{equation}
 Now we will compute an important relation for the quantity $\bar{C}$. Let us observe that
 \begin{equation}
     \bar{C}=\frac{1}{4\rho T^2}\int|v-\bm u|^4 f dv - \frac{5}{4}.
 \end{equation}
 At this point it is possible to use the Hölder's inequality 
 \begin{equation}
 \label{eq::Holder}
    \int\left|g(v)h(v)\right|dv < \left(\int|g(v)|^{p} dv\right)^{1/p} \left(\int|h(v)|^{q}dv\right)^{1/q},  
 \end{equation}
with $p=q=2$, $h(v)=\sqrt{f}$, $g(v)=\sqrt{f}|v-\bm u|^2$ to set a lower bound on the fourth moment appearing in the last expression of $\bar{C}$. Please note that the inequality \eqref{eq::Holder} holds strict, since $h$ and $g$ are not linearly dependent \cite[Section 6.2]{folland1999}. \\
 After some computation, one gets
 \begin{equation}
     \int f |v- \bm u|^4 dv>9\rho T^2,
 \end{equation}
which implies $\bar{C}>1$. Since it is a positive quantity, it is possible to compute the Schur complement \cite{horn2012} of $\mathbf{M}$. This operation preserves the eventual positive semi-definiteness of the matrix, in the sense that once the leading block is known to be positive semi-definite, the semi-definiteness of the full matrix is equivalent to the one of the Schur complement $\tilde{\mathcal{S}}$. Then, we obtain
 \begin{equation}
     \tilde{\mathcal{S}} =
     \begin{pmatrix}
      1-\displaystyle \frac{1}{\rule{0pt}{1.2em}\bar{C}}     & \displaystyle \frac{\bar{\mathbf{B}}^\intercal}{\rule{0pt}{1.2em}\bar{C}} \\[8pt]
      \displaystyle \frac{\bar{\mathbf{B}}}{\rule{0pt}{1.2em}\bar{C}} & \mathbf{I}+\bar{\mathbf{A}}- \displaystyle\frac{\bar{\mathbf{B}}\bar{\mathbf{B}}^\intercal}{\rule{0pt}{1.2em}\bar{C}}
     \end{pmatrix}
 \end{equation}

 Since $\bar{C}>1$, then we can compute again the Schur complement $\mathcal{S}$ obtaining
 \begin{equation}\label{eqmatS}
     \mathcal{S}=\mathbf{I}+\bar{\mathbf{A}} - \displaystyle\frac{\bar{\mathbf{B}}\bar{\mathbf{B}}^\intercal}{\rule{0pt}{1.2em}\bar{C}-1}
 \end{equation}
 So the original matrix $\mathbf{M}$ is semidefinite positive if and only if the matrix $\mathcal{S}$ is also semidefinite positive.\\
 
 Let us compute the matrix $\mathcal{S}$ at order zeroth and first in Chapman-Enskog expansion, by using results obtained in Section \ref{Section::CE}. We start by computing $\bar{\mathbf{A}}$, $\bar{\mathbf{B}}$ and $\bar{C}$ at order zeroth and first.

 \paragraph{\textbf{Order 0}}
     \begin{equation}
     \label{eq::A0}
         \bar{\mathbf{A}}_0 = \frac{1}{\rho}\int \mathbf{A}(v-\bm u) \mathcal{M} dv=\mathbf{0},\\
    \end{equation}
        
    \begin{equation}
        \bar{\mathbf{B}}_0 = \frac{1}{\rho}\int \mathbf{B}(v-\bm u) \mathcal{M} dv=\mathbf{0},\\
    \end{equation}

    \begin{equation}
    \label{eq::C0}
        \bar{C}_0 = \frac{1}{\rho}\int\left(\frac{|v-\bm u|^2}{2T}-\frac{5}{2}\right)^2 \mathcal{M} dv=\frac{5}{2}\\
    \end{equation}

 \paragraph{\textbf{Order 1}}
     \begin{equation}
     \label{eq::A1}
         \bar{\mathbf{A}}_1 = \frac{1}{\rho}\int \mathbf{A}(v-\bm u) \left(\mathcal{M}+\varepsilon h\right)dv= \mathbf{0} +\varepsilon \frac{\beta}{T}\mathbf{\Theta}_1 - \varepsilon \frac{\mathbf{D}(\bm u)}{\tau}=-\varepsilon \mathbf{D}(\bm u)\frac{\mu}{\rho T},\\
    \end{equation}
    where we used the fact that the tensor $\mathbf{\Theta_1}$ is traceless, $\mathbf{\Theta}_1 = \mathbf{\Theta}_1^\intercal$, $\mathbf{\Theta}_1=-\frac{ T}{(1-\beta)}\frac{\mathbf{D}(\bm u)}{\tau}$ and that the viscosity $\mu=\frac{\rho T}{\tau(1-\beta)}$ \cite{fil-jin-2011}.
        
    \begin{equation}
        \bar{\mathbf{B}}_1 = \frac{1}{\rho}\int \mathbf{B}(v-\bm u) \left(\mathcal{M}+\varepsilon h\right) dv=\mathbf{0}-\varepsilon \frac{\kappa}{\rho T^{\frac{3}{2}}}\nabla_x T,\\
    \end{equation}
    where we recall that the thermal conductivity $\kappa=\frac{5}{2}\frac{\rho T}{\tau}$ \cite{fil-jin-2011}.

    \begin{equation}
    \label{eq::C1}
        \bar{C}_1 = \frac{1}{\rho}\int\left(\frac{|v-\bm u|^2}{2T}-\frac{5}{2}\right)^2 \left(\mathcal{M}+\varepsilon h\right) dv=\frac{5}{2} + \varepsilon 0\\
    \end{equation}
    
At this point we can now compute the matrix $\mathcal{S}$ at order zeroth and first of the Chapman-Enskog Expansion, by simply plugging \eqref{eq::A0}-\eqref{eq::C0} and \eqref{eq::A1}-\eqref{eq::C1} into equation \eqref{eqmatS}.\\
\paragraph{\textbf{Order 0}}
 At order zeroth we have
 \begin{equation}
     \mathcal{S}_0 = \mathbf{I}.
 \end{equation}
\paragraph{\textbf{Order 1}}
 At order first we have
 \begin{equation}
     \mathcal{S}_1 = \mathbf{I} - \varepsilon \mathbf{D}(\bm u)\frac{\mu}{\rho T} - \frac{2}{3} \left(\varepsilon \frac{\kappa}{\rho T^{\frac{3}{2}}}\right)^2 \left(\nabla_x T\right) \left(\nabla_x T\right)^\intercal.
 \end{equation} 
 
 Then, the idea is to monitor the evolution of the matrix $\mathcal{S}$ \cite{tiw-2000, filbetrey2015, lev-mor-nad-1998}. Indeed, the moment realizability criterion states that the fluid dynamic description breaks down when the perturbation is too large such that $\mathcal{S}$, and thus $\bm{M}$, is no longer positive semi-definite. More details on how to build the hybrid solver are discussed in the following Section.

\subsubsection{Breakdown criteria}
The goal of this subsection is defining a criterion to determine the appropriate model -- fluid or kinetic -- to be used \cite{tiw-2000, lev-mor-nad-1998, filbetrey2015, filbetxiong2018, caparello2025}. 

The underlying idea is the following: we evaluate the quantity $\mathcal{S}$ defined in Equation \eqref{eqmatS} at two successive orders of the Chapman-Enskog expansion in $\eps$. Since higher order expansions yield increased model accuracy, each order is naturally associated with a different model. In particular, the zero-th order in the Chapman-Enskog expansion corresponds to the Euler equations (representing the lowest accuracy), while the first-order expansion corresponds to the ES-BGK equation (representing the highest accuracy).\\
Therefore, we obtain two different matrices corresponding to the different order of the Chapman-Enskog expansion in $\eps$: $\mathcal{S}_{0}$, associated with the zero-th order expansion and the Euler equations, and $\mathcal{S}_{1}$, associated to first order expansion and the ES-BGK equation. In each cell of the spatial domain, both matrices can be evaluated, and consequently compared in order to determine whether a transition between levels is required. If, in a given cell, the matrix $\mathcal{S}_{1}$ differs \textit{significantly} from $\mathcal{S}_{0}$, this indicates that the zero-th order layer (Euler equations) is no longer sufficiently accurate, and the first-order layer (kinetic, ES-BGK) must be employed instead. To quantify this difference, we compare the eigenvalues of these two matrices \cite{filbetrey2015, filbetxiong2018, caparello2025}.\\ 
From Section \ref{Section::CE} we already know that at order zeroth
\begin{equation}
     \mathcal{S}_0 = \mathbf{I},
 \end{equation}
while at first order
 \begin{equation}
     \mathcal{S}_1 = \mathbf{I} - \varepsilon \mathbf{D}(\bm u)\frac{\mu}{\rho T} - \frac{2}{3}\left(\varepsilon \frac{\kappa}{\rho T^{\frac{3}{2}}}\right)^2 \left(\nabla_x T\right) \left(\nabla_x T\right)^\intercal.
 \end{equation} 
Hence, we claim that the zero-th order layer (Euler equations) is valid when the matrix $ \mathcal{S}_{1}$ behaves similarly to the reference matrix $\mathcal{S}_{0}=\mathbf{I}$, meaning that it remains positive definite and \emph{its eigenvalues are close to unity}. Conversely, the zero-th order layer description will become inadequate when the discrepancy between the eigenvalues $\lambda_{1}$ of $\mathcal{S}_{1}$ and the eigenvalues $\lambda_{0}=1$ of $\mathcal{S}_{0}$ becomes \emph{too} large. Therefore, our criterion for switching from the Euler description (Layer 0) to ES-BGK description (Layer 1) is given by the following condition   
\begin{equation}
  \label{crit:CompEuler}
    \exists \lambda_{1} \in\, {\rm Sp}(\mathcal{S}_{1}) \qquad \text{such that} \,\, \left| \lambda_{A} \right|=\left |\lambda_{1}-1\right | > \eta,
\end{equation}
where $\eta > 0$ is a small parameter (problem-dependent quantity). Indeed, if this conditions is satisfied, it means that the eigenvalues associated to the first order expansion deviate excessively from those corresponding to the zero-th order expansion, implying that the Layer-0 description is no longer accurate. Therefore, we adopt this criterion in our numerical scheme to automatically switch between the Euler equations (Layer 0) and the ES-BGK equation (Layer 1).

\begin{rem}
    This criterion must also be verified when transitioning from Layer 1 (kinetic, ES-BGK) to Layer 0 (Euler equations). In particular, the switch from Layer 1 to Layer 0 may be carried out only when the following condition is satisfied
    \begin{equation}
    |\lambda_A|=|\lambda_{1}(t, \mathbf{x})-1|\leq\eta \qquad \forall \,\, \lambda_{1} \in\, {\rm Sp}(\mathcal{S}_{1}).
    \end{equation}
    This condition is necessary, but not sufficient, to justify a transition from Layer 1 to Layer 0, as we will discuss below.
\end{rem}
	    
Whenever the full kinetic description of a gas is available (Layer 1), several methods exist to assess how far the gas is from thermal equilibrium, see for instance \cite{saintraymond2009}. \\
In our numerical scheme we compare for each point $(t,x)$, the distribution function $f$ obtained as solution of the ES-BGK equation, with the Maxwellian $\mathcal{M}$, whose moments coincide with those of $f$. Hence, our numerical scheme can automatically switch from Layer 1 (kinetic, ES-BGK) to Layer 1 (Euler equations) if the following condition is met
\begin{equation}
    \label{eqKintoFluid1}
    \left \| f(t,x,\cdot) - \mathcal M_{\rho,\bm{u},T}(t,x,\cdot) \right \|_{L^1_v} \leq \delta,
\end{equation}
where $\delta$ is a small parameter (which is a problem dependent quantity). Both conditions \eqref{crit:CompEuler} and \eqref{eqKintoFluid1} must be satisfied in order for the transition to Layer 0 o occur.

\subsection{Velocity and Spatial Discretization}
\label{Section:spatial_discretization}
Let us consider a general $2$D spatial cartesian grid $[0,L_x] \times [0,Ly]$, discretized using $N_x$ and $N_y$ cells in the $x$- and $y$-direction, respectively. Then, the spatial cartesian grid is discretized into volumes of size $\Delta x\Delta y$. Each cell is given by $$\mathcal{C}_{ij} = [x_{i-1/2}, x_{i+1/2}] \times [y_{j-1/2}, y_{j+1/2}], \qquad i=1,\dots,N_x, \,\, j= 1,\dots,N_y,$$ 
where $x_{i \pm 1/2} = x_i \pm \Delta x /2$ and $y_{j \pm 1/2} = y_j \pm \Delta y /2$.\\
For the kinetic model, we introduce a $2$D discrete velocity grid, representing the set of admissible velocities. Let us truncate the velocity space by fixing some bounds and setting a square grid $[-v_\text{max}, v_\text{max}]^2$ in the velocity space with a total of $N_v$ points and $\Delta v$ as the grid step. For the sake of simplicity and without loss of generality, we consider the same number of points in each direction. The continuous distribution function $f(t,x,v)$ is then replaced by a vector  
\begin{equation*}
    f_k(x,t) \approx f(t,x,v_k), \qquad k=1,\dots,N_v.
\end{equation*}
Using this approximation, the discrete kinetic model consists of a set of $N_v$ evolution equations in the velocity space of the form
\begin{equation}
    \partial_t f_k + v_k \cdot \nabla_x(f_k) = \frac{1}{\eps} \mathcal{Q}(f_k), \qquad 1 \leq k \leq N_v.
\end{equation}
The velocity and spatial discretizations lead to an approximation of the form
\begin{equation}
    f_{i,k}(t) \approx f(t,\bm x_i,v_k), \qquad \text{where} \,\, \bm x_i \in [0,L_x] \times [0,Ly], \,\, v_k \in [-v_\text{max}, v_\text{max}]^2.
\end{equation}
For the sake of completeness, we report here the Euler equations \eqref{eqHydroClosedEuler} in 2D written in a compact form as 
\begin{equation}
\label{eq::Euler_Compact}
    \partial_t U + \diverg_x \mathbf{F}(U) = 0,
\end{equation}
where
\begin{equation}
    U=
    \begin{pmatrix}
        \rho \\
        \rho u_x\\
        \rho u_y\\
        E\\
    \end{pmatrix}, \ \ 
    \mathbf{F}_x(U)=
    \begin{pmatrix}
        \rho u_x \\
        \rho u_x^2 + P\\
        \rho u_x u_y\\
        u_x (E+P)\\
    \end{pmatrix}, \ \ 
    \mathbf{F}_y(U)=
        \begin{pmatrix}
            \rho u_y \\
            \rho u_y u_x\\
            \rho u_y^2 + P\\
            u_y (E+P)\\
        \end{pmatrix},
\end{equation}
and $P=\left(\gamma - 1 \right)\left(E-\frac{1}{2}\rho\left(u_x^2 + u_y^2\right)\right)$ is the pressure, $\rho$, $u_x$, $u_y$, $E$ and $\gamma$ are, respectively, the density, velocity along $x$-direction, velocity along $y$-direction, energy and specific heat ratio.\\
In our numerical scheme, we employ a CWENO3 scheme with a Lax-Friedrichs flux reconstruction \cite{pupporusso1999, pupposemplice2023}, both for the discretization of the flux $\mathbf{F}$ of the Euler equations and for the discretization of the advection term of the ES-BGK model.\\ 
Since the focus of this work is on time integration, we shall use $f$ and $U$ to denote their spatially discrete counterparts (for $f$ including the velocity discretization), whenever no ambiguity may arise.

\subsection{Construction of the coupling model}
\label{Section::Coupling_EulHyb}

It is important to note that the coupled strategy presented in this section can be adapted to any kind of coupled system.\\
Let us describe the evolution of the domain decomposition method, by focusing on the transition between macroscopic and kinetic regime in a cell $\mathcal{C}_{ij}$ at time $t^n$. We define $\Omega^n_{\text{macro}}$ the set of macroscopic cells and $\Omega^n_{\text{kin}}$ the set of kinetic cells at time $t^n$, such that the entire domain $\Omega^n = \Omega^n_{\text{macro}} \cup \Omega^n_{\text{kin}}$. Assuming a regime transition for the cell $\mathcal{C}_{ij}$, two different situations may occur:

\begin{itemize}
    \item \emph{Transition from kinetic to fluid regime.} If $\mathcal{C}_{ij} \subset \Omega^n_{\text{macro}}$, we can reconstruct the hydrodynamic fields as
	      \begin{equation}
	        \label{eq:Coupling_KinFlu}
             U_{ij}^n := \left (\rho_{ij}^n, \bm{u}_{ij}^n, T_{ij}^n\right ) \simeq \left (\rho(t^n, x_i, y_j), \bm{u}(t^n, x_i, y_j), T(t^n, x_i, y_j)\right ), 
	      \end{equation}
	 computed as suitable moments of the distribution function $f_{ij}^n(v)$, using \ref{hypEquilib}.
    \item \emph{Transition from fluid to kinetic regime.} If $\mathcal{C}_{ij} \subset \Omega^n_{\text{kin}}$, we can reconstruct the particle distribution function $f_{ij}^n(v)$ as a Maxwellian distribution function
	    \begin{equation}
         \label{eq:Coupling_FluKin}
	      f_{ij}^{n}(v_k) = \mathcal{M}_{\rho_{ij}^n, \bm{u}_{ij}^n, T_{ij}^n} (v_k), \quad \forall \, v_k,
	    \end{equation}
    where $(\rho_{ij}^n, \bm{u}_{ij}^n, T_{ij}^n)$ is the numerical solution to the Euler equation \eqref{eq::Euler_Compact}. 
\end{itemize}

\begin{rem}
\label{rem_coupling}
    The procedure described above is not valid only to describe a regime transition from $t^{n-1}$ to $t^{n}$, but should be applied at each time step to reconstruct the numerical fluxes in buffer regions. Indeed, the numerical approximation of $U$ and $f$ in the cell $\mathcal{C}_{ij}$ requires information from the neighbor cells, according to the stencil of the method. To clarify this idea, let us consider an example where the cell $\mathcal{C}_{ij}$ is kinetic and the cell $\mathcal{C}_{i-1\,j}$ is macroscopic. In this setting, we need to reconstruct $f_{i-1\,j}$, in order to evaluate the transport and evolve $f_{ij}$. Conversely, this procedure applies to macroscopic cells having kinetic cells in their stencil.
\end{rem}

\subsection{Time Integration}
\label{Section::Time_Discretization}
Let us consider the following semidiscrete version of the macroscopic Euler system
\begin{equation}
\label{eq::semidiscrete_Euler}
    \partial_t U = -\Phi(U),
\end{equation}
where $\Phi$ is a proper discretization of $\diverg_x \mathbf{F}$, as described in Section  \ref{Section:spatial_discretization}.\\
Analogously, let us introduce the semidiscrete version of the kinetic ES-BGK model as
\begin{equation}
\label{eq::semidiscrete_ESBGK}
    \partial_t f = -\mathcal{L}(f) + \frac{1}{\eps} \mathcal{Q}(f),
\end{equation}
where $\mathcal{L}$ represents a suitable discretization in the phase space of the advection operator $v \cdot \nabla_x f$, as described in Section \ref{Section:spatial_discretization}.\\

\subsubsection{High-Order Runge-Kutta Schemes}
\label{Section_RK}
Let us recall the main properties of Runge-Kutta (RK) methods, which will be employed for the time discretization of the Euler equations. We consider an explicit $S$-stage RK method for the Equation \eqref{eq::semidiscrete_Euler}.
Let us consider a time step discretization given by $\Delta t$, then  we have
\begin{equation}
    \begin{cases}
        U^{(s)} = U^n + \Delta t \displaystyle \sum_{l=1}^{s-1} a_{s,l} \mathbf{k}_l,\\[10pt]
        \mathbf{k}_s = - \Phi \left(t^n+c_s\,\Delta t, U^{(s)} \right),\\[10pt]
    \end{cases} \quad 1\leq s \leq S,
\end{equation}
and finally compute
\begin{equation}
    U^{n+1} = U^n + \Delta t \displaystyle \sum_{l=1}^{S} b_l \mathbf{k}_l.
\end{equation}
The matrix $A=(a_{s,l})_{s,l=1}^S$ is the so-called \textit{Runge-Kutta matrix}, the vector $b=(b_s)_{s=1}^S$ and $c=(c_s)_{s=1}^S$ are the \textit{Runge-Kutta weights} and the \textit{Runge-Kutta nodes}, respectively. The weights $\mathbf{b}_s$ and $\mathbf{c}_s$ are chosen to ensure consistency of the method, requiring $0 \leq b_s \leq 1$ and $0 \leq c_s \leq 1$ and
\begin{equation}
\label{Assumptions_RK}
    \sum_{s=1}^S b_s = 1, \sum_{s=1}^{S-1} a_{s,l} = c_s, \quad 1 \leq s \leq S.
\end{equation}
We recall that each Runge-Kutta scheme can be expressed in term of its Butcher's tableau as
\begin{equation*}
    \begin{tabular}
{c|c}
$c$ &
$A$\\[3pt]
\hline
\rule{0pt}{3ex} 
& $b^T$
\end{tabular}.
\end{equation*}

\subsubsection{Coupling Conditions}
\label{Section:Coupling_Conditions}
The construction of classical Implicit-Explicit (IMEX) Runge-Kutta schemes is based on the application of an implicit discretization to the source terms and an explicit one to the non-stiff term \cite{ascher1997, pareschirusso2005, pieraccini2007}. We refer to \cite{boscarinopareschirusso2025} for a complete review on the topic. This approach applied to system \eqref{eq::semidiscrete_ESBGK} takes the form
\begin{eqnarray}
 f^{(i)} & = & f^n - \Delta t \sum_{j=1}^{i-1} \ta_{ij} \mathcal{L}(f^{(j)}) +
                 \Delta t \sum_{j=1}^i a_{ij} \frac{1}{\varepsilon} \mathcal{Q}(f^{(j)}),
                                                                 \label{eq:RKEI1}\\
 f^{n+1} & = & f^n - \Delta t \sum_{i=1}^{s} \tb_{i} \mathcal{L}(f^{(i)})  +
                 \Delta t \sum_{i=1}^s b_{i} \frac{1}{\varepsilon} \mathcal{Q}(f^{(i)}).
                                                                 \label{eq:RKEI2}
\end{eqnarray}
The discretization is performed by choosing $\Delta t$, which is the time step.\\
The matrices $\tilde{A} = (\tilde{a}_{ij})$, $\tilde{a}_{ij}=0$ for $j \geq i$ and $A = (a_{ij})$ are the $s \times s$ Runge-Kutta matrices chosen such that the resulting scheme is explicit in $\mathcal{L}$ and implicit in $\mathcal{Q}$. For the sake of clarity, the reader can think about $A$ as a lower-triangular matrix (typically diagonally-implicit method) and $\tilde{A}$ as a strongly lower-triangular matrix. The corresponding weight (or coefficient) vectors are denoted by $\tilde{b}= (\tilde{b}_1,\ldots,\tilde{b}_s)^T$, $b =(b_1,\ldots,b_s)^T$.
Several techniques have been proposed to accelerate the computation, without increasing the complexity of the system.
The order conditions derived in \cite{ascher1997, pareschirusso2005} for IMEX schemes are based on the Taylor expansion of the exact and numerical solutions. More precisely, conditions for schemes of order $p$ are obtained by imposing that the Taylor expansion of the solution to system \eqref{eq::semidiscrete_ESBGK} at time $t^n + \Delta t$ matches the Taylor expansion of the Runge-Kutta scheme, up to $\mathcal{O}(\Delta t^{p+1})$.\\
For coupled problems, e.g. domain decomposition methods, new conditions arise because of the coupling of different schemes. In the particular setting of our domain decomposition strategy, a standard Runge-Kutta method is used for the discretization of the Euler equations and an IMEX scheme \cite{fil-jin-2011} is employed for the ES-BGK equation.\\
This approach applied to systems \eqref{eq::semidiscrete_Euler}-\eqref{eq::semidiscrete_ESBGK} leads to \eqref{eq:RKEI1}-\eqref{eq:RKEI2} for the ES-BGK model, coupled with the following additional equation for the Euler system 
\begin{eqnarray}
U^{(i)} & = & U^n - \Delta t \sum_{j=1}^{i-1} \ba_{ij} \Phi(U^{(j)}), \label{eq:RKEuler}\\
U^{n+1} & = & U^n - \Delta t \sum_{i=1}^{s} \bb_i \Phi(U^{(i)}), \label{eq:RKEuler_2}
\end{eqnarray}
where the matrix $\bA$ can be thought as a strongly lower-triangular matrix. This Runge-Kutta scheme is also characterized by the coefficient weights $\bb$.
In this case, the coupling involves all the Equations \eqref{eq:RKEI1}-\eqref{eq:RKEuler_2}, because at each stage $i=1,\dots,S$, we must solve the two systems in a fully coupled manner. At every stage, this requires both reconstructing the Maxwellian distribution functions or computing the corresponding moments in the buffer regions, as described in Section \ref{Section::Coupling_EulHyb}. This procedure is essential to preserve the overall high-order accuracy of the scheme, since each intermediate Runge-Kutta stage must incorporate the updated information provided by the other system, as we will show in Section \ref{Section_Simulations}. As a consequence, no system can be advanced independently of the other and coupling conditions must be provided.\\
The order conditions up to $p=3$ are given by the following relations.\\

\noindent \textit{First order}.\\
\begin{equation}
    \sum_{i=1} \bb_i = 1, \quad \sum_{i=1} \tb_i = 1, \quad \sum_{i=1} b_i = 1.
   \label{eq:cond1}
\end{equation}
\noindent \textit{Second order.}\\
\begin{equation}
    \sum_i \bb_i \bc_i = 1/2, \quad \sum_i \tb_i \tc_i = 1/2, \quad \sum_i b_i c_i = 1/2,
   \label{eq:cond2}
\end{equation}

\begin{equation}
    \begin{aligned}
   \sum_i \bb_i c_i = 1/2, \quad \sum_i \bb_i \tc_i = 1/2, \quad \sum_i \tb_i c_i = 1/2, \\
   \sum_i \tb_i \bc_i = 1/2, \quad \sum_i b_i \bc_i = 1/2. \quad \sum_i b_i \tc_i = 1/2.
   \label{eq:mixed2}
   \end{aligned}
\end{equation}

\noindent \textit{Third order.}\\
\begin{equation}
\begin{aligned}
   \sum_{ij} \bb_i \ba_{ij} \bc_j = 1/6, \quad
   \sum_{i}  \bb_i \bc_{i}  \bc_i = 1/3, \\
   \sum_{ij} \tb_i \ta_{ij} \tc_j = 1/6, \quad
   \sum_{i}  \tb_i \tc_{i}  \tc_i = 1/3, \\
   \sum_{ij}   b_i   a_{ij}   c_j = 1/6, \quad
   \sum_{i}    b_i   c_{i}    c_i = 1/3,
\label{eq:cond3}
\end{aligned}
\end{equation}
\begin{equation}
\begin{aligned}
&\sum_{i,j} \beta_i \, \alpha_{ij} \, \gamma_j = \frac{1}{6}, 
   && \forall (\beta, \alpha, \gamma) \in \{b, \bb, \tb\} \times \{a, \ba, \ta\} \times \{c, \bc, \tc\},\\[2mm]
&\sum_i \beta_i \, \gamma_i \, \gamma'_i = \frac{1}{3}, 
   && \forall (\beta, \gamma, \gamma') \in \{b, \bb, \tb\} \times \{c, \bc, \tc\} \times \{c, \bc, \tc\}.
   \end{aligned}                                     \label{eq:mixed3}
\end{equation}
The appearance of additional conditions naturally arises from the need to couple all three schemes simultaneously. Indeed, conditions \eqref{eq:cond1}, \eqref{eq:cond2}, \eqref{eq:cond3} are
the standard order conditions for the three {\em tableau\/}, each of them taken separately. Conditions \eqref{eq:mixed2} and \eqref{eq:mixed3} are new conditions that arise because of the coupling of the three schemes.\\
In our setting, the order conditions will simplify a lot, since we can set $\bA=\tA$, $\bb=\tb$ and $\bc=\tc$. Indeed, the Butcher's tableau used for the approximation of the explicit term in the ES-BGK model can also be employed for the time integration of the Euler equations. This automatically reduces the coupling conditions to the ones described in \cite{pareschirusso2005} for IMEX schemes.
In our code we implement the temporal third order IMEX DIRK scheme ARS(2,3,3) \cite{ascher1997}, which is characterized by the following two Butcher tableaus
\begin{center}
    \begin{tabular}{c|c}
        $\tilde{c}$ & $\tilde{A}$\\                
        \hline
        \\
                    & $\tilde{b}^\intercal$\\
    \end{tabular}, \quad  \quad   
    \begin{tabular}{c|c}
        $c$ & $A$\\                
        \hline
        \\
            & $b^\intercal$\\
    \end{tabular}.
\end{center}
given by
\begin{center}
    \begin{tabular}{c|c c c}
        $0$        & $0$        & $0$         & $0$\\   
        $\alpha$   & $\alpha$   & $0$         & $0$\\ 
        $1-\alpha$ & $\alpha-1$ & $2-2\alpha$ & $0$\\
        \hline
                   & $0$          & $1/2$       & $1/2$\\
    \end{tabular}, \quad  \quad  
    \begin{tabular}{c|c c c}
        $0$        & $0$ & $0$         & $0$\\   
        $\alpha$   & $0$ & $\alpha$    & $0$\\ 
        $1-\alpha$ & $0$ & $1-2\alpha$ & $\alpha$\\
        \hline
                   & $0$          & $1/2$       & $1/2$\\
    \end{tabular},
\end{center}
where $\alpha = (3+\sqrt3)/6$. The tableau on the left is used to integrate in time the explicit term, while the one on the right is used for the implicit integration.\\
Different IMEX-RK methods can be applied but their discussion lies beyond the scope of this work. For further details, we refer the reader to \cite{hu2017, boscarino2025, boscarinopareschirusso2025} and references therein. 

\begin{rem}
As already mentioned in the Introduction, different schemes can be introduced to solve the kinetic equation. For instance, projective integration methods \cite{lafitte2012, melis2019} have been proposed to deal with more complicated collision operators for which IMEX schemes are not easily constructed (e.g. the full Boltzmann equation). The coupling strategy for domain decomposition methods can be easily extended to these time integrators, reaching high-order of accuracy in time.
\end{rem}

\subsection{Numerical Algorithm}
\label{sec_temporal_algorithm}
Let us describe the numerical approach to implement this strategy, which is summarized in Algorithm \ref{algoritmo1}. The key aspect is the coupling between the two models \eqref{eq:Coupling_KinFlu}-\eqref{eq:Coupling_FluKin}, which must be performed at each stage of the Runge-Kutta scheme, see also Remark \ref{rem_coupling}. This coupling is fundamental to maintain the overall high-order of accuracy of the method. If the coupling is only performed at each time step $t^n$ (not at each stage), the order of accuracy reduces to $1$. In the numerical tests of Section \ref{Section_Numerical} we will refer to this latter approach as low-order coupling (LOC) hybrid method.\\
\begin{algorithm}
\caption{High-order Domain Decomposition Method for the ES-BGK Equation}
\label{algoritmo1}
\begin{algorithmic}[1]
\REQUIRE Initial condition $\left(f^0, U^0\right)$, time step $\Delta t$, Butcher tableaux $(c_i, a_{ij}, b_i)$, $(\tc_i, \ta_{ij}, \tb_i)$ and $(\bc_i, \ba_{ij}, \bb_i)$
\WHILE{$t<T$}
    \STATE Regime update according to \eqref{crit:CompEuler}
    \FOR{$i = 1$ to $s$} 
    \STATE \{Compute stages\}
        \STATE Macroscopic solution: $U^{(i)} = U^n - \Delta t \displaystyle \sum_{j=1}^{i-1} \ba_{ij} \Phi(U^{(j)})$\\[7pt]
        \STATE Kinetic solution: $f^{(i)} =  f^n - \Delta t \displaystyle \sum_{j=1}^{i-1} \ta_{ij} \mathcal{L}(f^{(j)}) + \Delta t \sum_{j=1}^i a_{ij} \frac{1}{\varepsilon} \mathcal{Q}(f^{(j)})$\\[7pt]
        \STATE Coupling through \eqref{eq:Coupling_KinFlu}-\eqref{eq:Coupling_FluKin}
    \ENDFOR
    \STATE Update macroscopic solution: $U^{n+1} = U^n - \Delta t \displaystyle \sum_{i=1}^{s} \bb_i \Phi(U^{(i)})$\\[7pt]
    \STATE Update kinetic solution: $f^{n+1} = f^n - \Delta t \displaystyle \sum_{i=1}^{s} \tb_i \mathcal{L}(f^{(i)}) + \Delta t \sum_{i=1}^s b_i \frac{1}{\varepsilon} \mathcal{Q}(f^{(i)})$
    \STATE Coupling through \eqref{eq:Coupling_KinFlu}-\eqref{eq:Coupling_FluKin}
    \STATE Update time $t=t+\Delta t$
\ENDWHILE
\RETURN $(f^N, U^N)$
\end{algorithmic}
\end{algorithm}
Here it is important to remark that the coupling operation on the seventh row of the Algorithm \ref{algoritmo1} is what allows the scheme to maintain the overall temporal order of accuracy. Indeed, as numerically demonstrated in Section \ref{Section_Simulations}, if the coupling procedure is executed only after each timestep (eleventh row in the Algorithm \ref{algoritmo1}) the overall temporal order of accuracy reduces to one, while our approach is able to preserve the temporal third order of accuracy of the ARS233 numerical scheme.

\subsection{Properties of the scheme}
We observe that for $\eps=0$ the regime indicators naturally lead to a fully hydrodynamic domain, meaning that we automatically solve the Euler equations. Indeed, we recall the expression of $\mathcal{S}_1$ associated to the first order of the Chapman-Enskog expansion is 
 \begin{equation}
     \mathcal{S}_1 = \mathbf{I} - \varepsilon \mathbf{D}(\mathbf{u})\frac{\mu}{\rho T} -\frac{2}{3} \left(\varepsilon \frac{\kappa}{\rho T^{\frac{3}{2}}}\right)^2 \left(\nabla_x T\right) \left(\nabla_x T\right)^\intercal.
 \end{equation}     
For $\eps=0$ we have $\mathcal{S}_1 = \mathcal{S}_0$ and $\lambda_{A}=0$, which implies that the inequality \eqref{crit:CompEuler}, namely
\begin{equation}
    \exists \lambda_{1} \in\, {\rm Sp}(\mathcal{S}_{1}) \qquad \left| \lambda_{A} \right| > \eta,
\end{equation}
is never satisfied for any value of $\eta > 0$. This means that for $\eps=0$, the regime indicators identify the entire domain as hydrodynamic and the scheme reduces to solving the Euler equations.\\
Let us now focus on the Navier-Stokes limit and let us assume $\eps$ is chosen such that the regime indicators lead to a full kinetic domain, meaning that we never solve the Euler equations. In this regime, the preservation of the asymptotic behavior is automatically inherited from the AP properties of the IMEX scheme for the ES-BGK equation, as shown in \cite{fil-jin-2011}. This means that the combination of a domain decomposition method with an asymptotic-preserving scheme for the kinetic equation ensures high efficiency and accuracy in all regimes \cite{jin2022}. For the sake of completeness, we mention a very recent paper on the asymptotic analysis of IMEX-RK schemes for the ES-BGK equation at the Navier-Stokes level \cite{boscarino2025}. \\
Concerning conservation properties of the scheme, we observe that the same spatial discretization is employed for both the Euler and the ES-BGK equations. This choice ensures consistency between the two descriptions and, in particular, guarantees the global conservation properties of the overall numerical method. By relying on a unified discretization framework, the scheme preserves the fundamental conservation structure of both models across regime transitions, as already observed in \cite{dimarcomieussens2014}.

\section{Numerical Simulations}
In this Section we perform some numerical tests in order to validate our numerical scheme. A first numerical test is performed to show the novelty of this work, i.e. a coupling time integrator able to preserve the overall time order accuracy. Then we perform some tests 2D in space and 2D in velocity to show the ability of the scheme to deal with long-time simulations and handle complex geometries, like a flow over a cylinder or a rectangular obstacle. We will show that in all the simulations, the kinetic regime is always well confined in space.\\
It is important to remark that the ES-BGK equation is integrated in time by extending the first order IMEX AP scheme proposed by \cite{fil-jin-2011} to third-order using the ARS233 scheme \cite{ascher1997, pareschirusso2005}, whose Butcher's tableau is illustrated in Section \ref{Section:Coupling_Conditions}.\\
All simulations have been conducted on a dedicated computer equipped with two \emph{Intel(R) Xeon(R) Gold 5418Y} CPUs and 128 GB of RAM, kindly provided by the Mathematics and Computer Science Department of the University of Ferrara. In order to speed up the code execution, a CPU multithread approach has been implemented: each thread operates on a subdomain of the entire spatial domain. In all simulations we implement a value $\beta=-0.5$.\\
\label{Section_Simulations}
\subsection{One-dimensional case}
\subsubsection{Test 1: Order of accuracy}
\begin{figure}
    \centering
    \begin{subfigure}{0.45\textwidth}
    \centering
    \includegraphics[width=\textwidth, trim={1cm 0cm 1cm 0cm}]{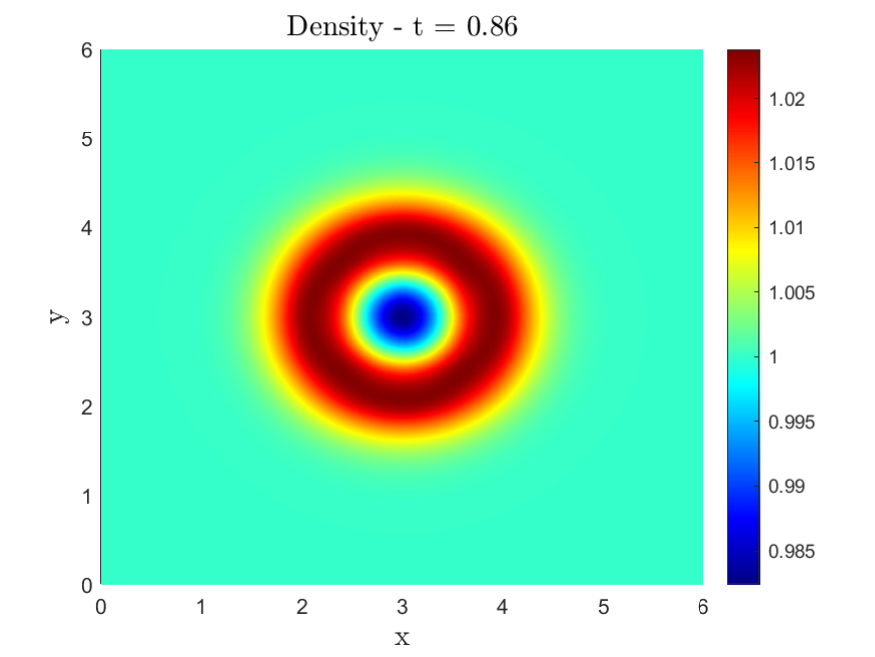}
    \end{subfigure}
    \begin{subfigure}{0.45\textwidth}
    \centering
    \includegraphics[width=\textwidth]{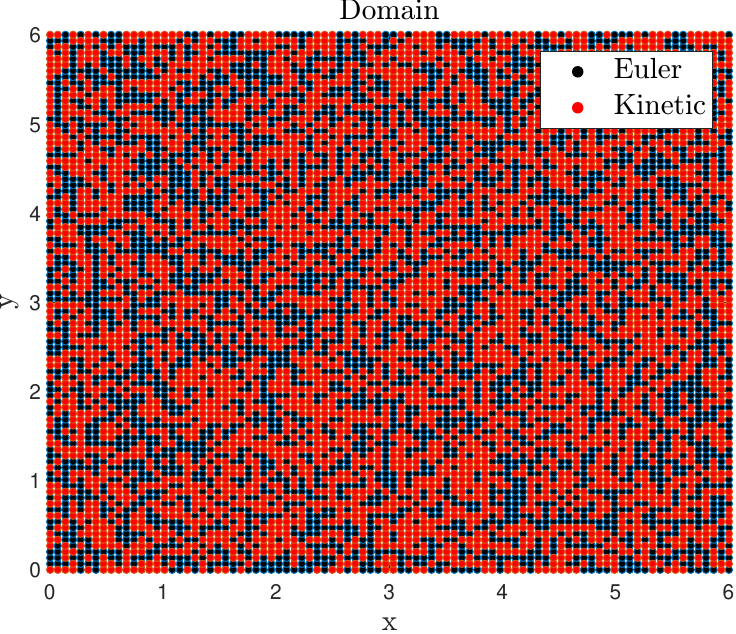}
    \end{subfigure}
    \hfill
    \caption{On the left, density at final time $t_{end}=0.86$ obtained with the smallest time step and Knudsen number $\eps=10^{-6}$, using the hybrid scheme with high-order coupling. On the right, domain decomposition used in all simulations for testing the order.}
    \label{figDensityDomainOrder}
\end{figure}

In this subsection, we test the order of accuracy of the time integrator for the domain decomposition method combined with the strategy described in Section \ref{Section::Time_Discretization}. Since we would like to compare solutions computed with different $\Delta t$, we assume that the domain partition is constant in time, arbitrary chosen at $t_0 = 0$.\\
As initial condition, we consider the following space-dependent moments
\begin{equation}
    \begin{cases}
            \rho(x)=1+0.1\exp(-\lambda d^2),\\            T(x)=0.3\left(1+0.1\exp(-\lambda d^2)\right),\\
            u_x(x)=0,\\
            u_y(x)=0,
    \end{cases}
\end{equation}
where $\lambda=3$ and $d=\sqrt{(x-x_c)^2+(y-y_c)^2}$, where $(x_c, y_c)$ are the coordinates of the center of the domain. For the kinetic region, the distribution function $f$ is initialized as a Maxwellian having those macroscopic moments.\\ 
The final time is $t_{end}=0.86$.
The space domain is $[0, L_x]\times [0, L_y]$, where $L_x=L_y=6$, discretized using $N_x=N_y=90$ cells in both directions. The velocity space $\Omega_v=[-5, 5]^2$ is discretized using $N_v=32$ points along both directions. In order to stress as much as possible the coupling approach, we randomly initialize the regime in each cell of the domain (macroscopic or kinetic), as shown in Figure \ref{figDensityDomainOrder}.\\
The temporal order of the method is proven using an approach based on Richardson extrapolation \cite{richardson1911, tenna2025}, where two solutions computed with different time step $\frac{\Delta t}{k}$, $k \in \mathbb{N}$ are compared.\\
Let us assume a time discretization of order $p$ (e.g. IMEX-RK of order $p$), then the solution $u^{n+1}_{\Delta t/k}$ is deviating from the exact solution $u^{n+1}_*$ according to 
\begin{equation}
    u^{n+1}_{\Delta t / k} = u^{n+1}_* + k\, C_n \, \left( \frac{\Delta t}{k} \right)^{p+1} + \mathcal{O} \left( \Delta t^{p+2} \right),
\end{equation}
where $C_n$ is the error constant.
This leads to the following relation 
\begin{equation}
    u^{n+1}_{\Delta t}-u^{n+1}_* = \frac{2^p \left( u^{n+1}_{\Delta t}-u^{n+1}_{\Delta t/2} \right)}{2^p-1}.
\end{equation}
In other words, it is possible to compute the order of convergence studying the ratio between two ``consecutive'' errors, i.e. errors computed with $\Delta t/k$ and $\Delta t/(2k)$. Therefore, let us define the relative error at final time $t_{end}$ as follows
\begin{equation}
\label{error_temporal}
    \text{err}_{\Delta t} = || u_{\Delta t}(t_{end}) - u_{\Delta t/2}(t_{end})||_2,
\end{equation}
where $||\cdot||_2$ is the discrete $L^2$-norm.

\begin{figure}
    \centering
    \begin{subfigure}{0.325\textwidth}
    \centering
    \includegraphics[width=\textwidth, trim={0cm 0cm 0cm 0cm}]{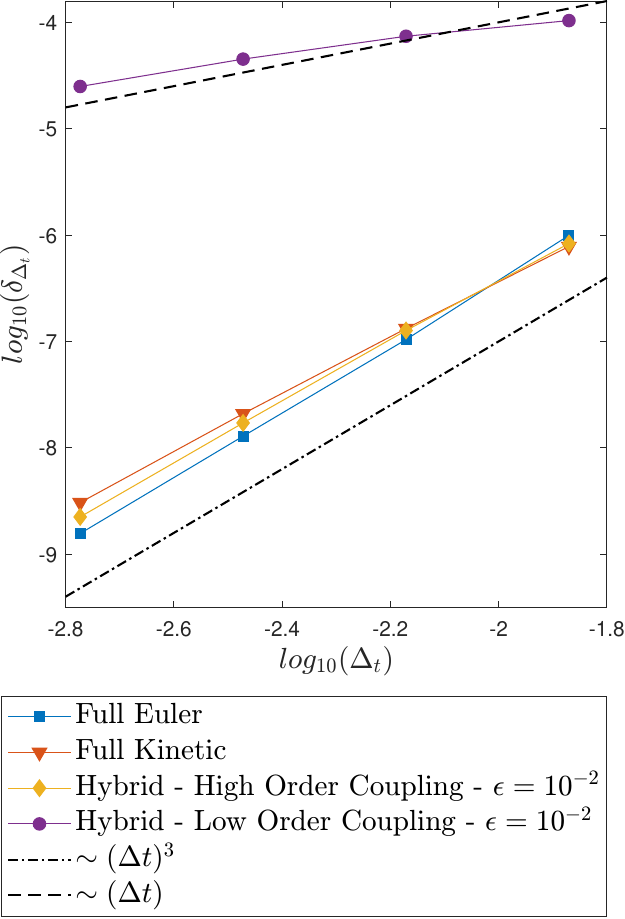}
    \end{subfigure}
    \begin{subfigure}{0.325\textwidth}
    \centering
    \includegraphics[width=\textwidth, trim={0cm 0cm 0cm 0cm}]{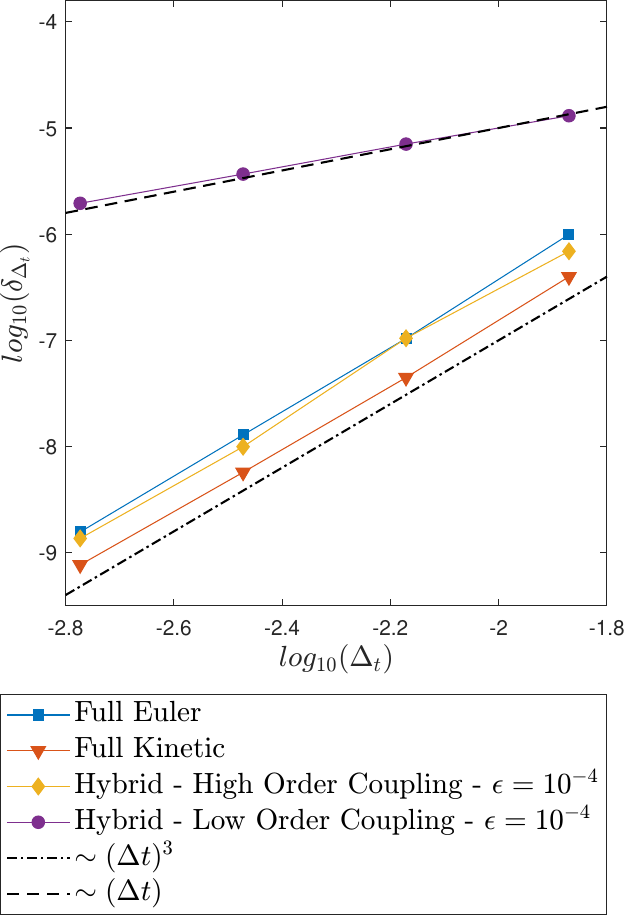}
    \end{subfigure}
    \begin{subfigure}{0.325\textwidth}
    \centering
    \includegraphics[width=\textwidth]{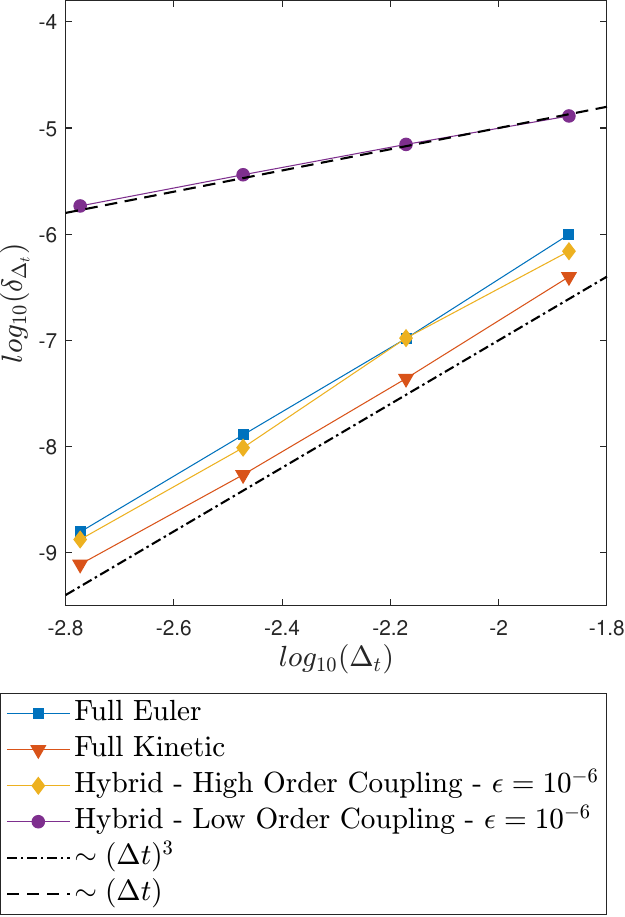}
    \end{subfigure}
    \hfill
    \caption{Convergence curves obtained using the full Euler Solver, the full Kinetic Solver and a comparison between the high-order coupling and low-order coupling, for three different values of the Knudsen number $\eps=10^{-2}$ (on the left), $\eps=10^{-4}$ (in the middle) and $\eps=10^{-6}$ (on the right).}
    \label{figErrorCurves}
\end{figure}

In the following we numerically prove that our coupling strategy between Euler and Kinetic solver, executed at each stage of the IMEX-RK is able to preserve the overall temporal order of accuracy of the scheme (see Section \ref{sec_temporal_algorithm}). We also demonstrate that the coupling between solvers executed only at the end of each timestep degrades the temporal accuracy to order 1 (see Section \ref{sec_temporal_algorithm}). We refer to this latter approach as Low-Order Coupling hybrid method (LOC-hybrid).\\ 
In Figure \ref{figErrorCurves} we plot the convergence curves of the full ES-BGK model, full Euler model, hybrid model \textbf{with} high-order coupling and hybrid model \textbf{without} high-order coupling (LOC-hybrid). The error curve is obtained using \eqref{error_temporal}.\\ We also report the order of accuracy for the hybrid method \textbf{with} and \textbf{without} the high-order coupling for $\eps=10^{-2}$ in Table \ref{tab:order_accuracy_2}, for $\eps=10^{-4}$ in Table \ref{tab:order_accuracy_4} and for $\eps=10^{-6}$ in Table \ref{tab:order_accuracy_6}.\\
In Figure \ref{figDensityDomainOrder} we show the final density (at time $t_{end}=0.86$) obtained with the smallest time step and $\eps=10^{-6}$, using the hybrid scheme with high-order coupling.\\

\begin{table}[ht!]
    \centering
    \begin{tabular}{|c|c c|c c|}
        \hline
        & \multicolumn{2}{c|}{\textbf{Hybrid ($\eps=10^{-2}$)}} & \multicolumn{2}{c|}{\textbf{LOC-Hybrid ($\eps=10^{-2}$)}}\\
        \hline
         $\mathbf{\Delta t}$ & \textbf{err}$_{\Delta t}$ & \textbf{order} & \textbf{err}$_{\Delta t}$ & \textbf{order}\\
        \hline
        0.0135 & $0.8319 \cdot 10^{-6}$ & // & $0.1038 \cdot 10^{-3}$ & // \\
        0.0067 & $0.1265 \cdot 10^{-6}$ & $2.72$ & $0.0740 \cdot 10^{-3}$ & $0.49$ \\
        0.0034 & $0.0171 \cdot 10^{-6}$ & $2.89$ & $0.0451 \cdot 10^{-3}$ & $0.71$ \\
        0.0017 & $0.0022 \cdot 10^{-6}$ & $2.96$ & $0.0249 \cdot 10^{-3}$ & $0.86$ \\
        \hline
    \end{tabular}
    \caption{Temporal error and order of convergence for the Hybrid domain decomposition method applied to the ES-BGK model for $\eps=10^{-2}$ \textbf{with} high-order coupling (on the left) and \textbf{without} high-order coupling (on the right).}
    \label{tab:order_accuracy_2}
\end{table}

\begin{table}[ht!]
    \centering
    \begin{tabular}{|c|c c|c c|}
        \hline
        & \multicolumn{2}{c|}{\textbf{Hybrid ($\eps=10^{-4}$)}} & \multicolumn{2}{c|}{\textbf{LOC-Hybrid ($\eps=10^{-4}$)}}\\
        \hline
         $\mathbf{\Delta t}$ & \textbf{err}$_{\Delta t}$ & \textbf{order} & \textbf{err}$_{\Delta t}$ & \textbf{order}\\
        \hline
        0.0135 & $0.6897 \cdot 10^{-6}$ & // & $0.1305 \cdot 10^{-4}$ & //\\
        0.0067 & $0.1045 \cdot 10^{-6}$ & $2.72$ & $0.0706 \cdot 10^{-4}$ & $0.89$ \\
        0.0034 & $0.0099 \cdot 10^{-6}$ & $3.39$ & $0.0368 \cdot 10^{-4}$ & $0.93$ \\
        0.0017 & $0.0014 \cdot 10^{-6}$ & $2.87$ & $0.0195 \cdot 10^{-4}$ & $0.92$ \\
        \hline
    \end{tabular}
    \caption{Temporal error and order of convergence for the Hybrid domain decomposition method applied to the ES-BGK model for $\eps=10^{-4}$ \textbf{with} high-order coupling (on the left) and \textbf{without} high-order coupling (on the right).}
    \label{tab:order_accuracy_4}
\end{table}

\begin{table}[ht!]
    \centering
    \begin{tabular}{|c|c c|c c|}
        \hline
        & \multicolumn{2}{c|}{\textbf{Hybrid ($\eps=10^{-6}$)}} & \multicolumn{2}{c|}{\textbf{LOC-Hybrid ($\eps=10^{-6}$)}}\\
        \hline
         $\mathbf{\Delta t}$ & \textbf{err}$_{\Delta t}$ & \textbf{order} & \textbf{err}$_{\Delta t}$ & \textbf{order}\\
        \hline
        0.0135 & $0.6903 \cdot 10^{-6}$ & // & $0.1298 \cdot 10^{-4}$ & // \\
        0.0067 & $0.1047 \cdot 10^{-6}$ & $2.72$ & $0.0700 \cdot 10^{-4}$ & $0.89$ \\
        0.0034 & $0.0098 \cdot 10^{-6}$ & $3.43$ & $0.0363 \cdot 10^{-4}$ & $0.95$ \\
        0.0017 & $0.0013 \cdot 10^{-6}$ & $2.87$ & $0.0184 \cdot 10^{-4}$ & $0.98$ \\
        \hline
    \end{tabular}
    \caption{Temporal error and order of convergence for the Hybrid domain decomposition method applied to the ES-BGK model for $\eps=10^{-6}$ \textbf{with} high-order coupling (on the left) and \textbf{without} high-order coupling (on the right).}
    \label{tab:order_accuracy_6}
\end{table}

\subsubsection{Test 2: Sod Shock Tube}
\begin{figure}
    \centering
    
    \begin{subfigure}{0.32\textwidth}
    \centering
    \includegraphics[width=\textwidth]{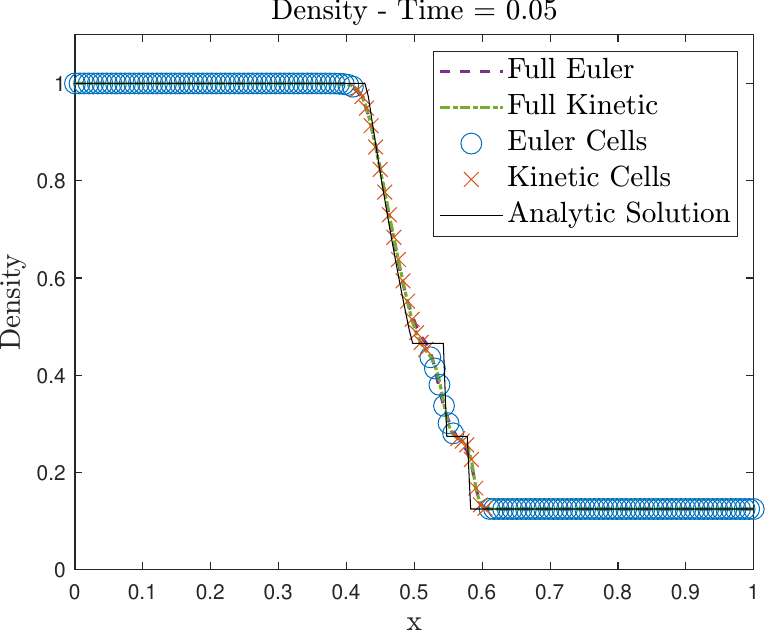}
    \end{subfigure}
    \begin{subfigure}{0.32\textwidth}
    \centering
    \includegraphics[width=\textwidth]{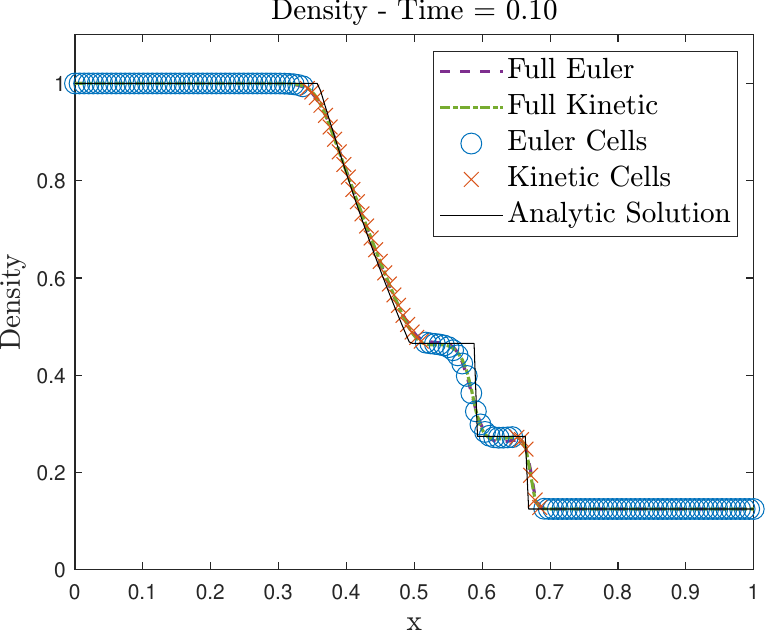}
    \end{subfigure}
    \begin{subfigure}{0.32\textwidth}
    \centering
    \includegraphics[width=\textwidth]{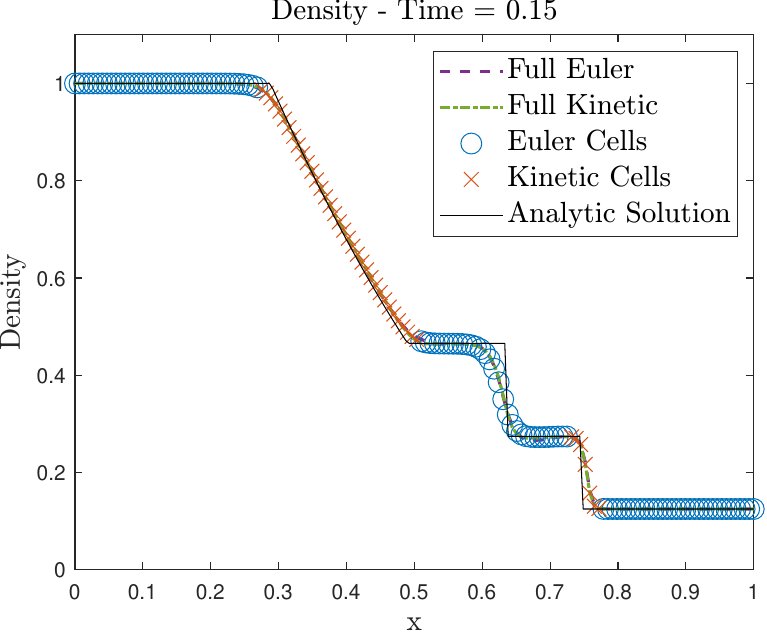}
    \end{subfigure}
    \hfill
    
    \begin{subfigure}{0.32\textwidth}
    \centering
    \includegraphics[width=\textwidth]{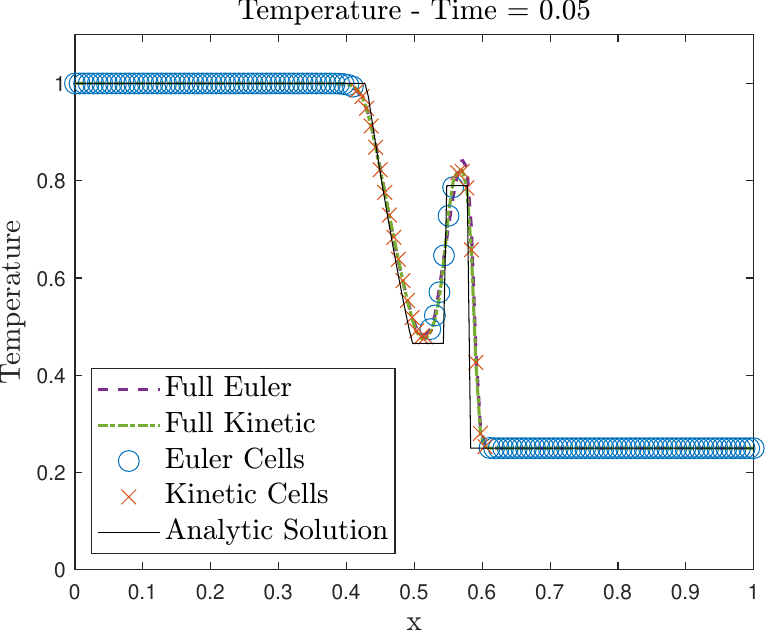}
    \end{subfigure}
    \begin{subfigure}{0.32\textwidth}
    \centering
    \includegraphics[width=\textwidth]{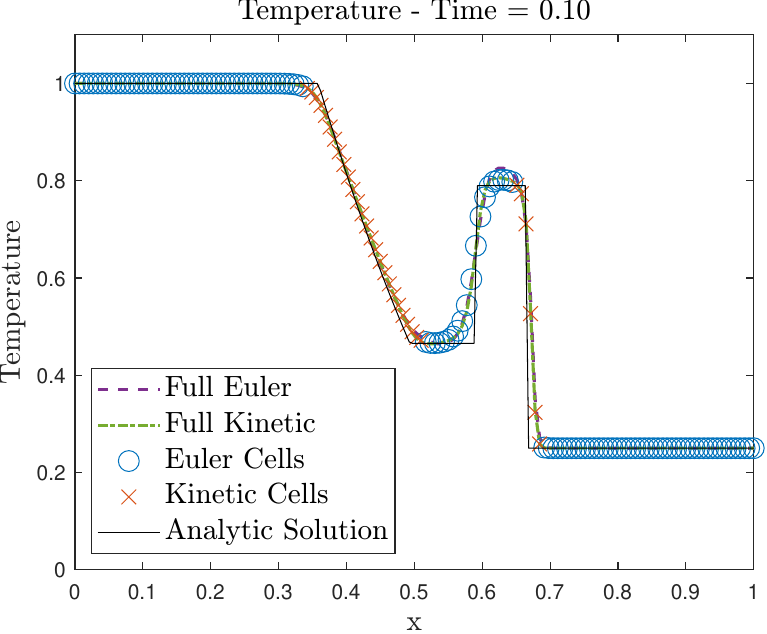}
    \end{subfigure}
    \begin{subfigure}{0.32\textwidth}
    \centering
    \includegraphics[width=\textwidth]{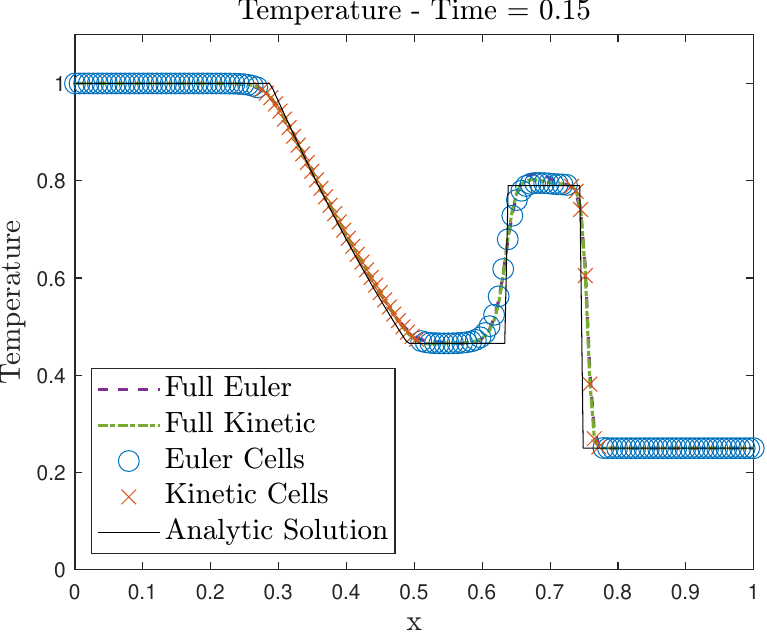}
    \end{subfigure}
    
    \begin{subfigure}{0.32\textwidth}
    \centering
    \includegraphics[width=\textwidth]{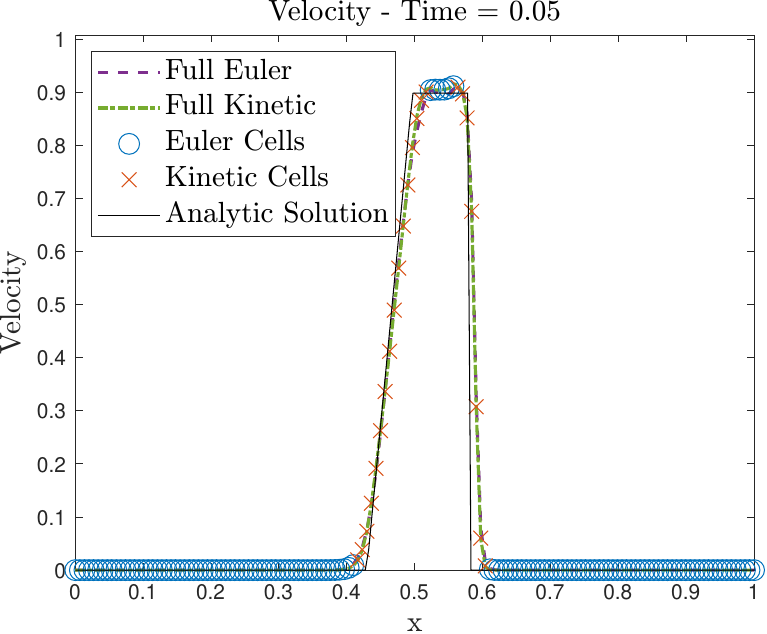}
    \end{subfigure}
    \begin{subfigure}{0.32\textwidth}
    \centering
    \includegraphics[width=\textwidth]{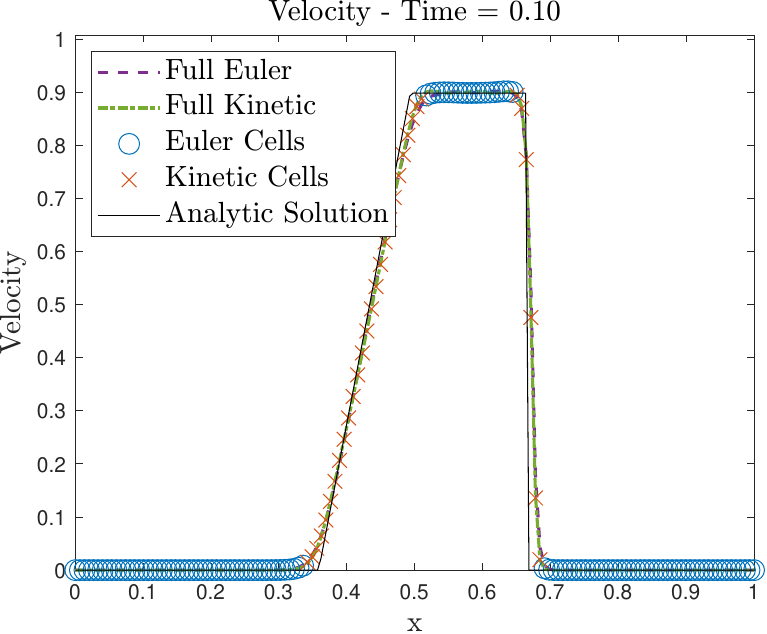}
    \end{subfigure}
    \begin{subfigure}{0.32\textwidth}
    \centering
    \includegraphics[width=\textwidth]{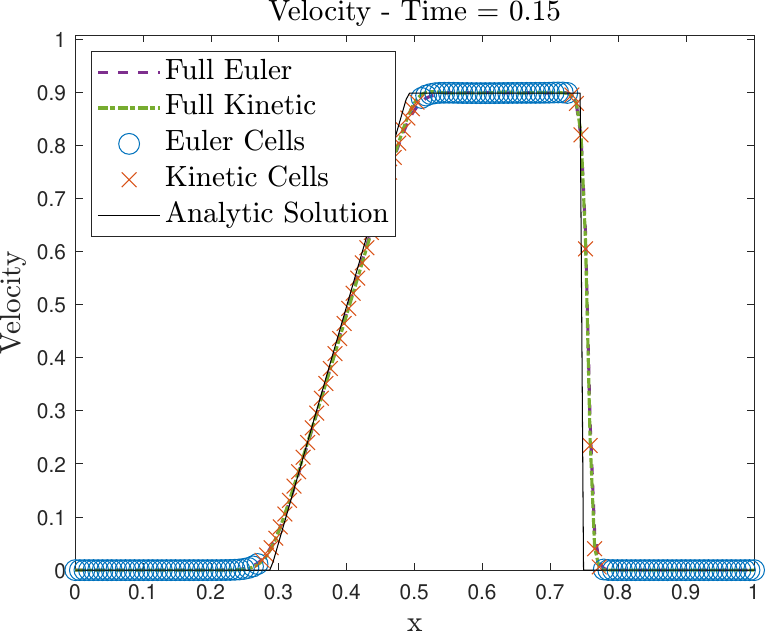}
    \end{subfigure}
    \caption{Density, Temperature and Velocity evolution at $t=0.05$, $t=0.10$ and $t=0.15$ for the Sod shock tube with initial condition \eqref{eq::initial_Sod}.}
    \label{figSod}
\end{figure}

Let us consider a classical benchmark test for the one-dimensional Euler equations, given by the following initial Riemann problem in $[0,1]$, known as Sod shock tube test \cite{sod1978}
\begin{equation}
\label{eq::initial_Sod}
	\begin{cases}
		\begin{aligned}
		&\rho= \rho_L, \qquad &u=u_L, \qquad &p=p_L \qquad &\text{ if } y \leq 0.5,\\
		&\rho= \rho_R, \qquad &u=u_R, \qquad &p=p_R \qquad &\text{ if } y > 0.5.
		\end{aligned}
	\end{cases}
\end{equation}
The left state is prescribed as $U_L =(\rho_L, u_L, p_L) = (1, 0, 1)$, whereas the right state as $U_R =(\rho_R, u_R, p_R) = (2^{-3}, 0, 2^{-5})$. Since our numerical implementation is 2D in space, we initialize the moments uniformly along the $x$-direction. The solution to the Euler equations for this particular configuration consists of a rarefaction wave, a contact discontinuity and a shock wave. This macroscopic solution (obtained by analytically solving the Riemann problem) is the reference solution for the comparison to the numerical one computed approximating the ES-BGK equation in the limit $\eps \to 0$.\\
The physical domain is given by $[0, L_x] \times [0, L_y]$, where $L_x=0.1$, and $L_y=1$, discretized with $N_x=21$ cells, and $N_y=150$ cells, along the $x-$ and $y-$direction respectively. The velocity domain $\Omega_v=[-L,L]^2$, where $L=8$ is discretized using $N_v=32$ points in both directions.\\
The time step used is $\Delta t = 0.8\min{(\Delta x, \Delta y)} / L$ and the final time of integration is $t_{end}=0.15$.\\
We initialize the distribution function everywhere as a Maxwellian, whose moments coincide with those of the macroscopic states $U_L$ and $U_R$. Hence, the computational domain is initially treated as fully macroscopic, after which the regime of each cell is updated at every time step.\\
In Figure \ref{figSod} we compare density, temperature and velocity obtained using the hybrid scheme, the full Euler solver, the full Kinetic (ES-BGK) solver, along with the analytical solution obtained by solving the associated Riemann problem.\\
Let us consider a Knudsen number $\varepsilon=10^{-6}$. We would like to highlight that the adoption of an Asymptotic-Preserving IMEX scheme removes the dependence of the time step on the Knudsen number. As a result, the time step is  constrained only by the hyperbolic CFL condition associated with the transport term of the ES-BGK equation and not by $\eps$.\\
The thresholds used in switching criteria (Equations \eqref{crit:CompEuler} and \eqref{eqKintoFluid1}) for this test are $\eta=\eps$, $\delta=10^{-3}$.\\
This test allows to check that all schemes are converging towards the correct (analytical) solution, when $\varepsilon \to 0$, and that the coupling is correctly executed. Indeed it is clear that the kinetic regime is engaged only within regions which lie close to non-equilibrium conditions. This means that only when significant non-equilibrium effects arise (e.g. steep gradients, shocks or boundary layers), the method employs a kinetic description, whereas the remainder of the flow is captured by the macroscopic model. 

\subsection{Two-dimensional case}
\subsubsection{Test 3: Kelvin-Helmholtz Instability}
\begin{figure}
    \centering
    
    \begin{subfigure}{0.48\textwidth}
    \centering
    \includegraphics[width=\textwidth, trim={0.2cm 0cm 0cm 0cm}]{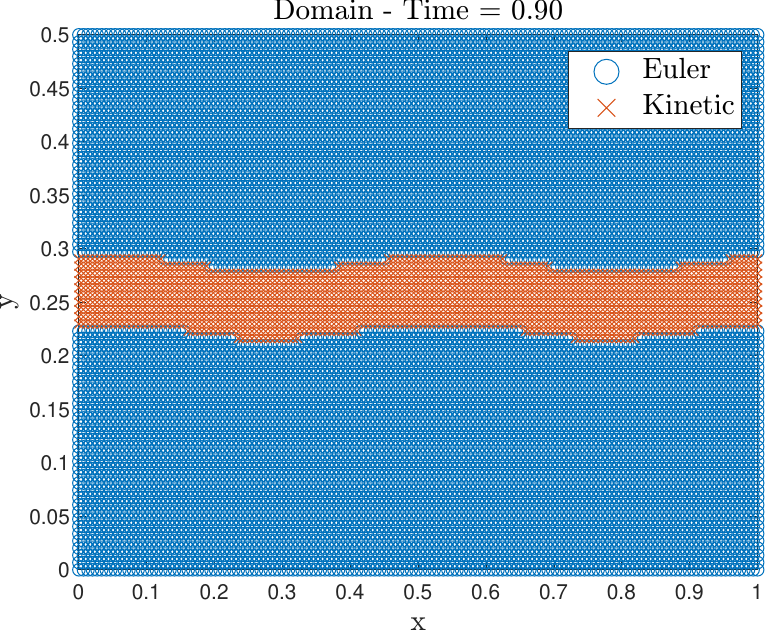}
    \end{subfigure}
    \begin{subfigure}{0.48\textwidth}
    \centering
    \includegraphics[width=\textwidth, trim={0.2cm 0cm 0cm 0cm}]{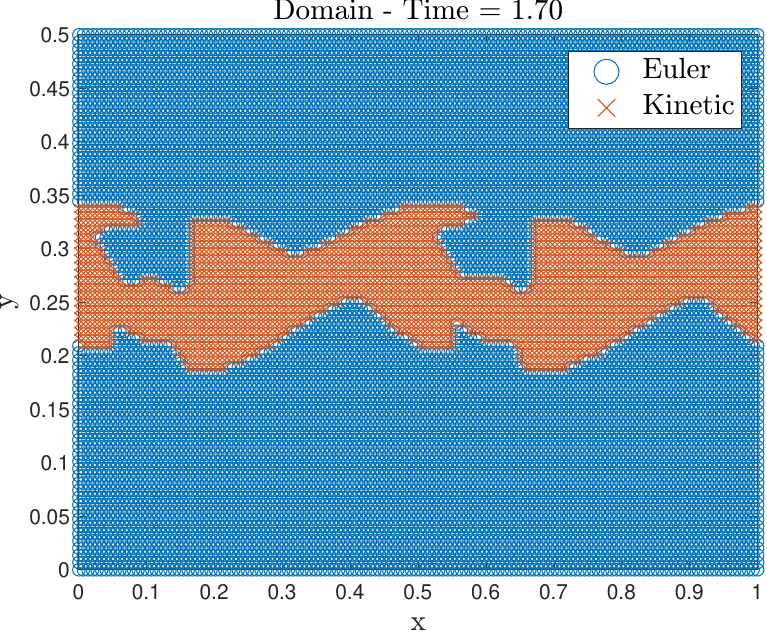}
    \end{subfigure}
    \hfill
    
    \begin{subfigure}{0.48\textwidth}
    \centering
    \includegraphics[width=\textwidth, trim={1cm 0cm 1cm 0cm}]{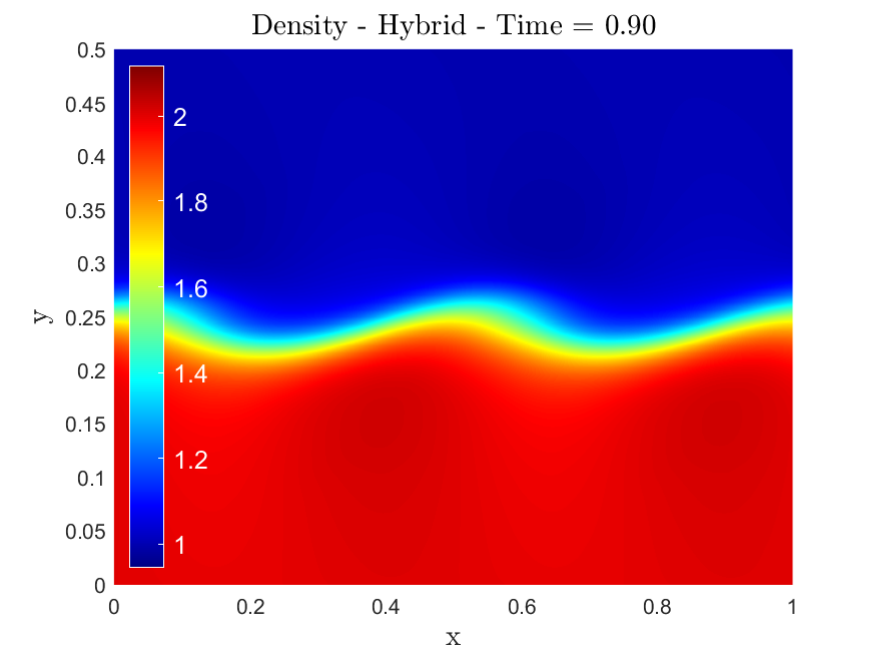}
    \end{subfigure}
    \begin{subfigure}{0.48\textwidth}
    \centering
    \includegraphics[width=\textwidth, trim={1cm 0cm 1cm 0cm}]{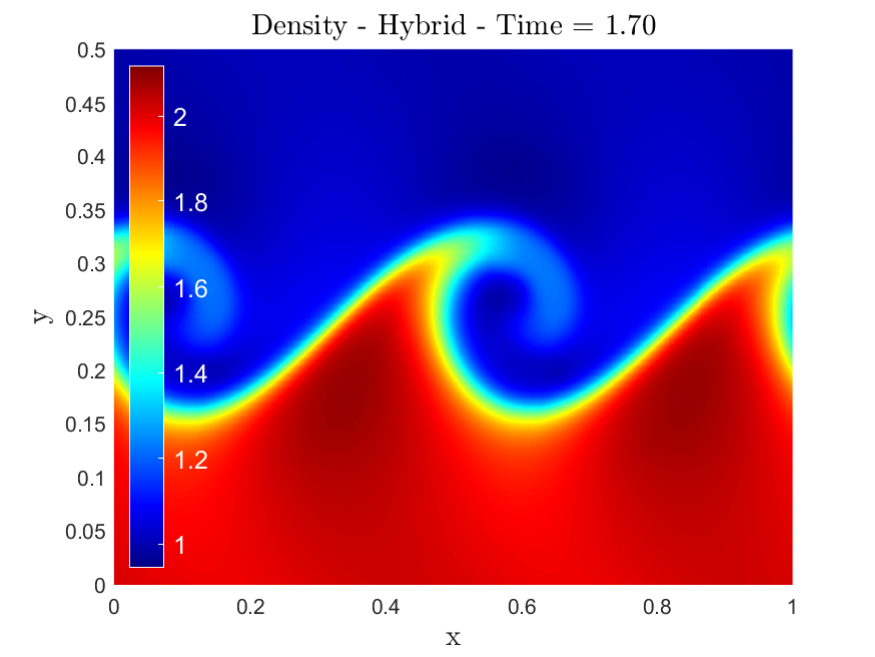}
    \end{subfigure}
    
    \begin{subfigure}{0.48\textwidth}
    \centering
    \includegraphics[width=\textwidth, trim={1cm 0cm 1cm 0cm}]{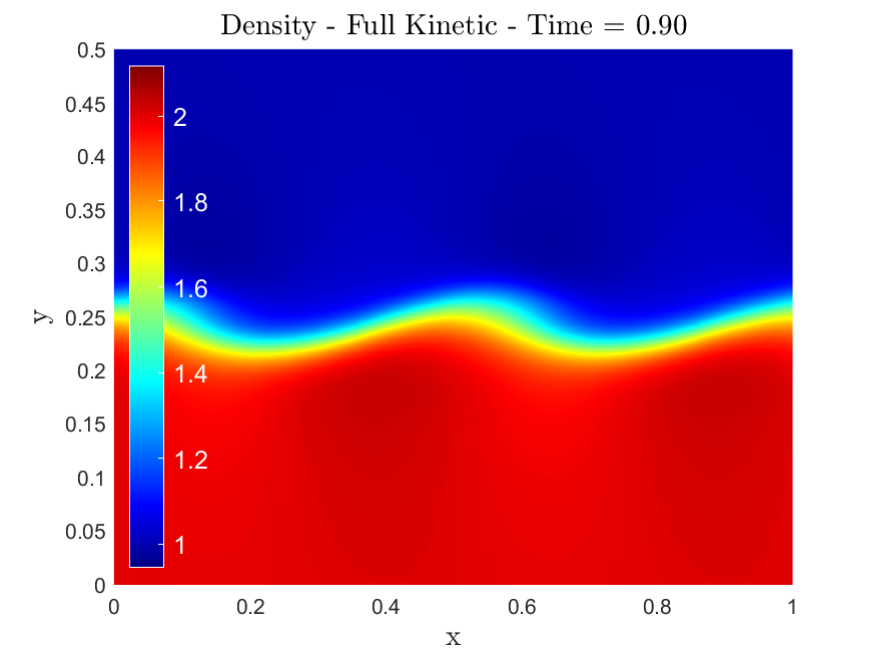}
    \end{subfigure}
    \begin{subfigure}{0.48\textwidth}
    \centering
    \includegraphics[width=\textwidth, trim={1cm 0cm 1cm 0cm}]{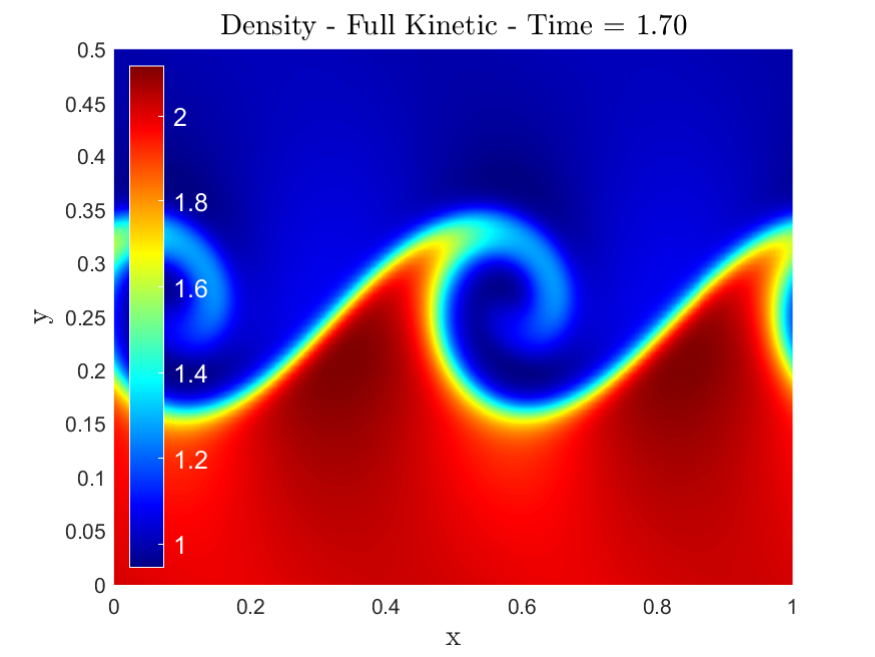}
    \end{subfigure}
    
    \begin{subfigure}{0.48\textwidth}
    \centering
    \includegraphics[width=\textwidth, trim={1cm 0cm 1cm 0cm}]{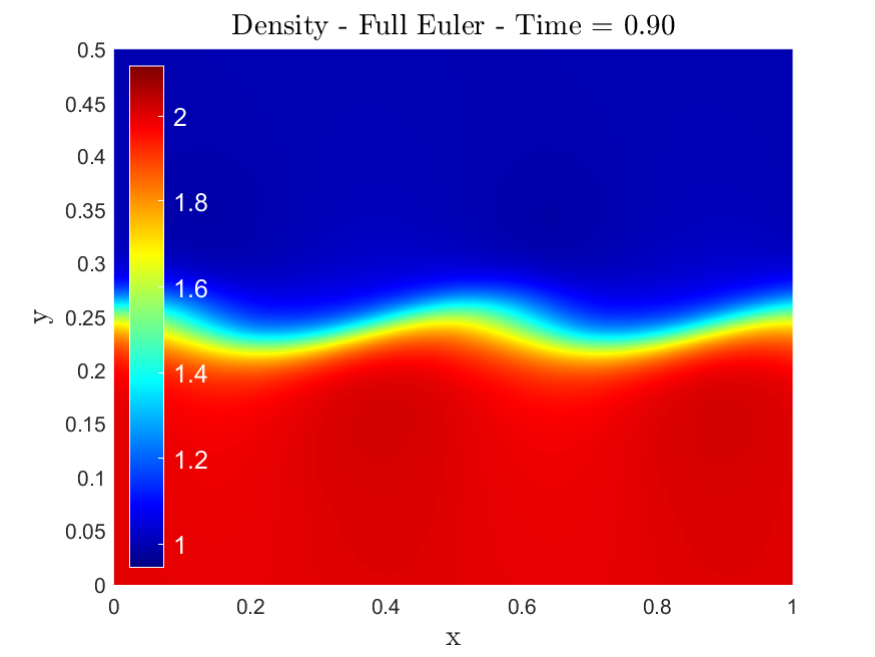}
    \end{subfigure}
    \begin{subfigure}{0.48\textwidth}
    \centering
    \includegraphics[width=\textwidth, trim={1cm 0cm 1cm 0cm}]{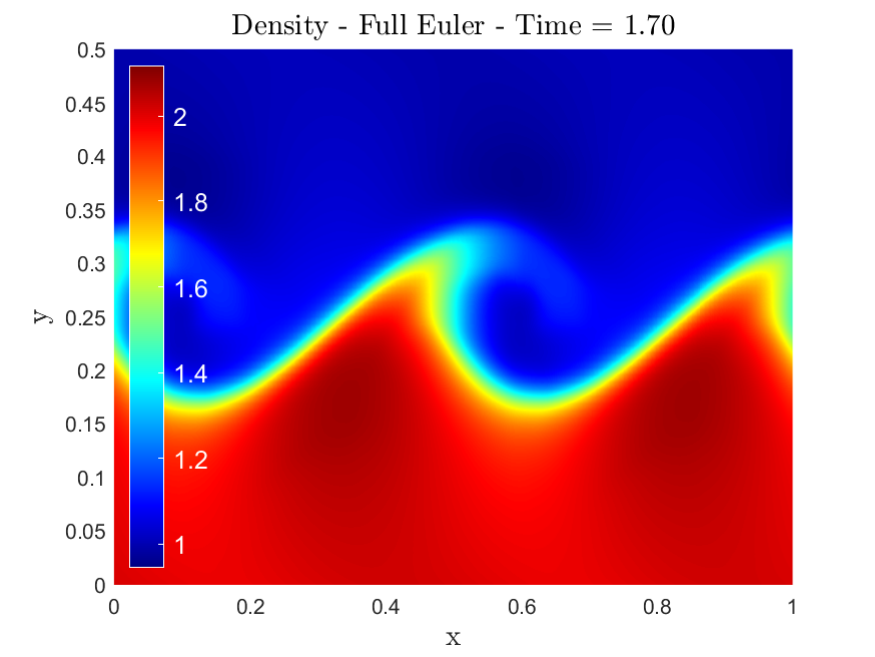}
    \end{subfigure}
    \caption{In the first row domain adaptation is shown, whereas the other rows present the densities obtained using the hybrid (second row), the full Kinetic (third row) and full Euler (fourth row) solvers, respectively, for the Kelvin-Helmholtz instability with initial condition \eqref{eq::initial_Kelvin}. The solutions are displayed at time $t=0.9$ and $t=1.7$.}
    \label{figKH}
\end{figure}

In this test, we investigate the capability of this scheme to capture instabilities arising from particular initial conditions, inspired by \cite{melis2019, bailo2022}.\\
The physical domain is given by $[0, L_x] \times [0, L_y]$, where $L_x=1$, and $L_y=0.5$, discretized with $N_x=150$ cells, and $N_y=75$ cells, along the $x-$ and $y-$direction respectively. The velocity domain $\Omega_v=[-L,L]^2$, where $L=8$ is discretized using $N_v=32$ points in both directions.\\
The time step used is $\Delta t = 0.5\min{(\Delta x, \Delta y)} / L$ and we consider as final time of integration $t_{end}=1.7$.
We set the Knudsen number $\varepsilon=10^{-6}$, thereby requiring an appropriate numerical method to treat the stiffness.\\
The fluid is initialized with the following moments
\begin{equation}
\label{eq::initial_Kelvin}
    \begin{cases}
    \begin{aligned}
		&\rho= \rho_D, \quad &u_x=u_{x_D}, \quad &u_y=0.01\sin{(4\pi x)}, \quad &p=p_D \qquad &\text{ if } y \leq L_y/2,\\
		&\rho= \rho_U, \quad &u_x=u_{x_U}, \quad &u_y=0.01\sin{(4\pi x)}, \quad &p=p_U \qquad &\text{ if } y > L_y/2,
		\end{aligned}
    \end{cases}
\end{equation}
where $\rho_D=2$, $\rho_U=1$, $u_{x_D}=-0.5$, $u_{x_U}=+0.5$, $p_D=1$ and $p_U=1$. Periodic boundary conditions are used along the $x$-direction.\\
The distribution function is initialized everywhere as a Maxwellian, whose moments coincide with those of the macroscopic states precised above in \eqref{eq::initial_Kelvin}.\\
For this test, we consider as thresholds $\eta=4\varepsilon$ and $\delta=10^{-3}$ in switching criteria (Equations \eqref{crit:CompEuler} and \eqref{eqKintoFluid1}) .\\
In Figure \ref{figKH} we compare the density obtained using the hybrid scheme, the full Euler solver and the full Kinetic (ES-BGK) solver. In the first row it is also reported the domain adaptation. It appears evident that the solutions are strongly consistent with one another, and that the kinetic regime is automatically and dynamically engaged only close to the discontinuities.

\subsubsection{Test 4: Shock Bubble Interaction}
\begin{figure}
    \centering
    \begin{subfigure}{0.32\textwidth}
    \centering
    \includegraphics[width=\textwidth, trim={0.2cm 0cm 0cm 0cm}]{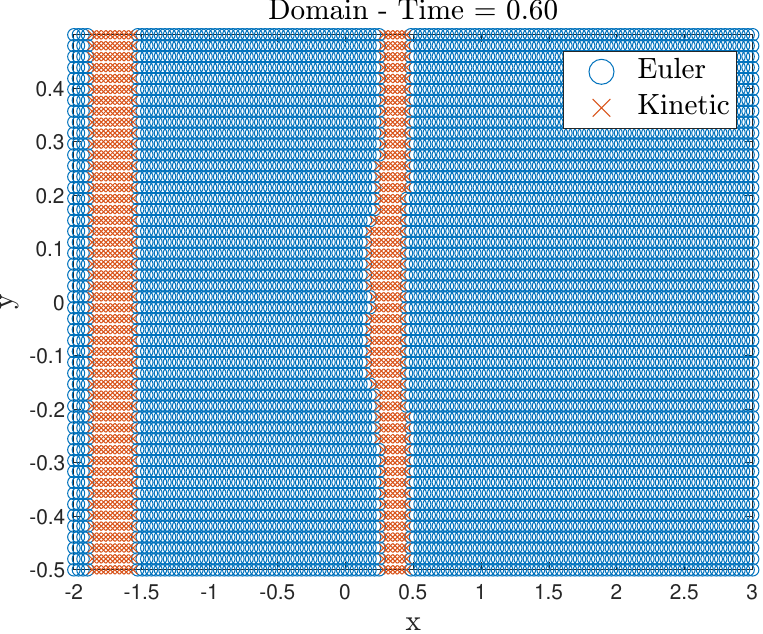}
    \end{subfigure}
    \begin{subfigure}{0.32\textwidth}
    \centering
    \includegraphics[width=\textwidth, trim={0.2cm 0cm 0cm 0cm}]{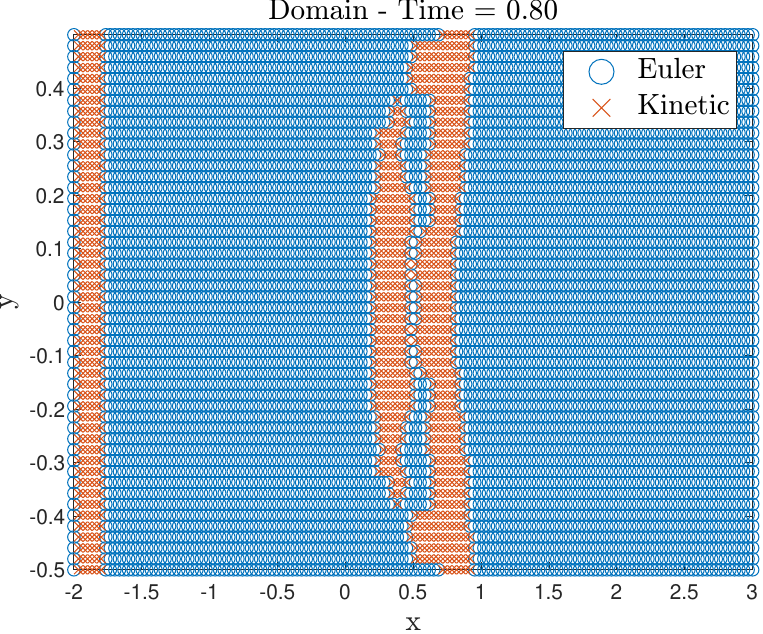}
    \end{subfigure}
    \begin{subfigure}{0.32\textwidth}
    \centering
    \includegraphics[width=\textwidth, trim={0.2cm 0cm 0cm 0cm}]{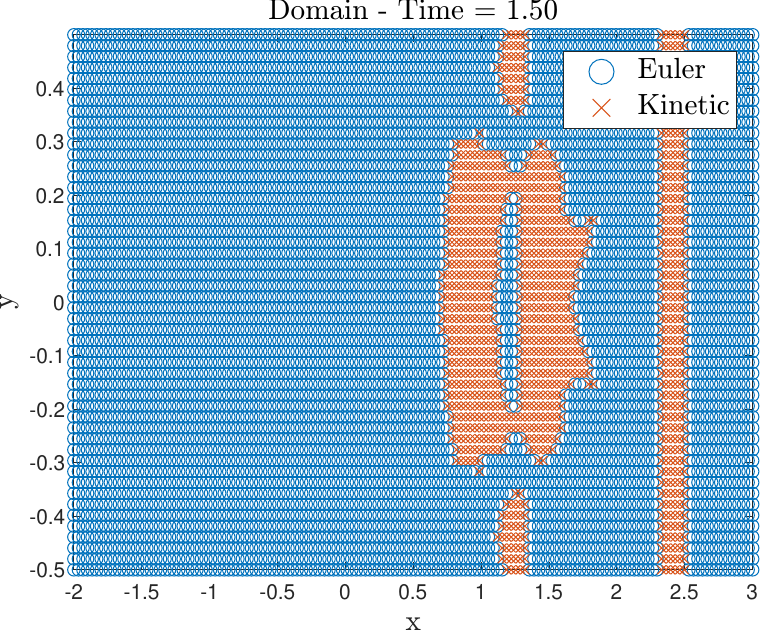}
    \end{subfigure}
    \hfill
    
    \centering
    \begin{subfigure}{0.32\textwidth}
    \centering
    \includegraphics[width=\textwidth, trim={1cm 0cm 1.8cm 0cm}]{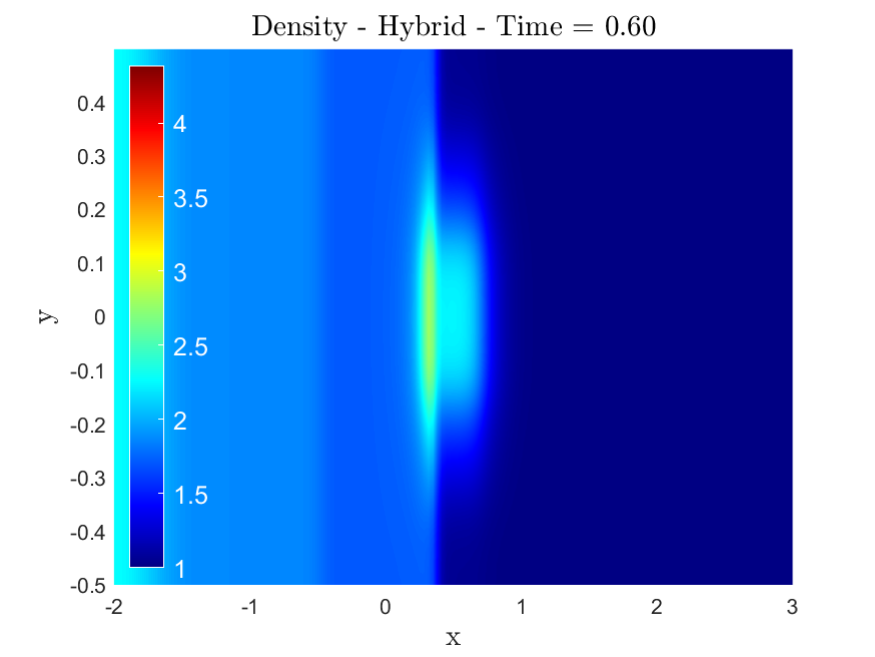}
    \end{subfigure}
    \begin{subfigure}{0.32\textwidth}
    \centering
    \includegraphics[width=\textwidth, trim={1cm 0cm 1.8cm 0cm}]{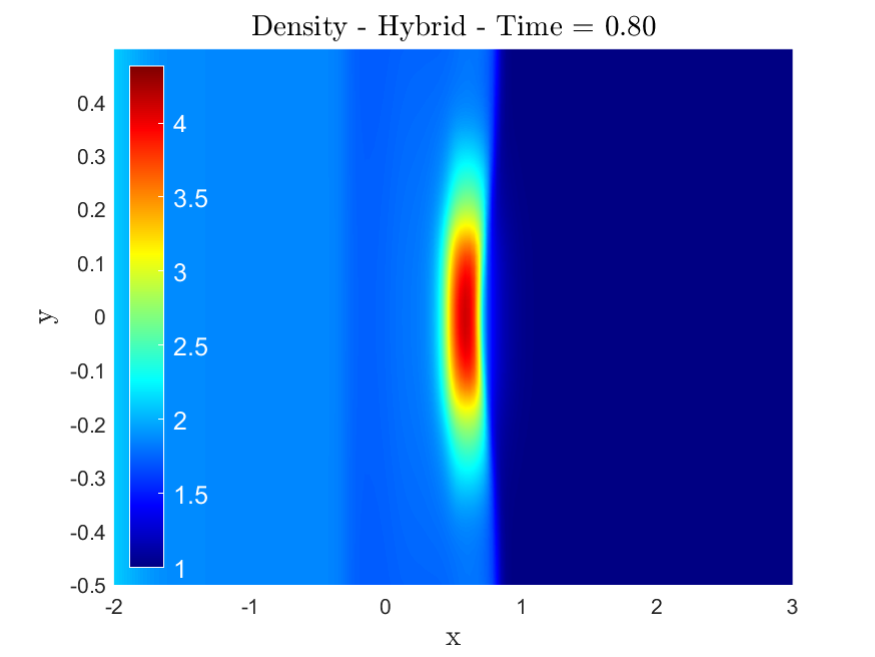}
    \end{subfigure}
    \begin{subfigure}{0.32\textwidth}
    \centering
    \includegraphics[width=\textwidth, trim={1cm 0cm 1.8cm 0cm}]{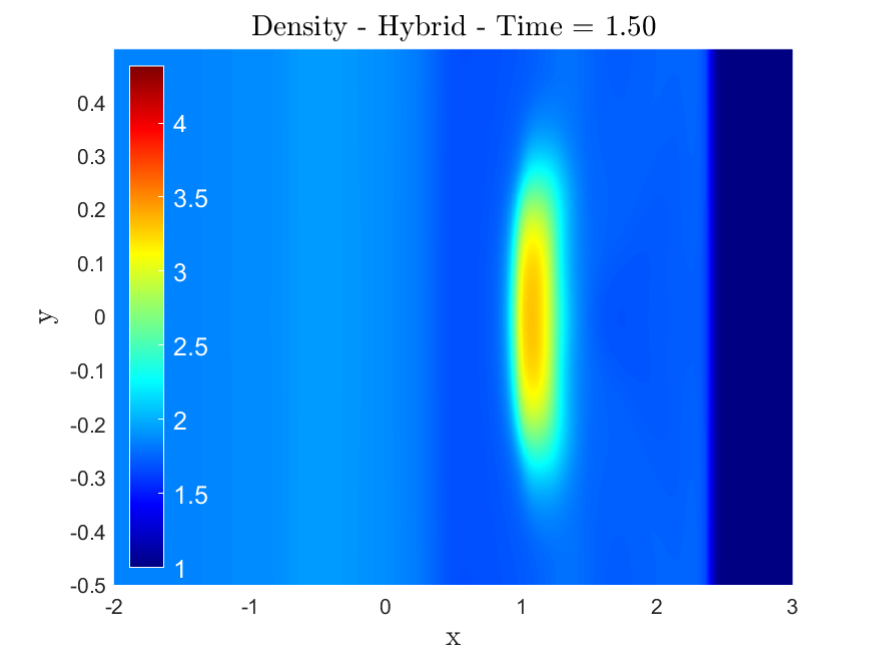}
    \end{subfigure}
    \hfill
    
    \centering
    \begin{subfigure}{0.32\textwidth}
    \centering
    \includegraphics[width=\textwidth, trim={1cm 0cm 1.8cm 0cm}]{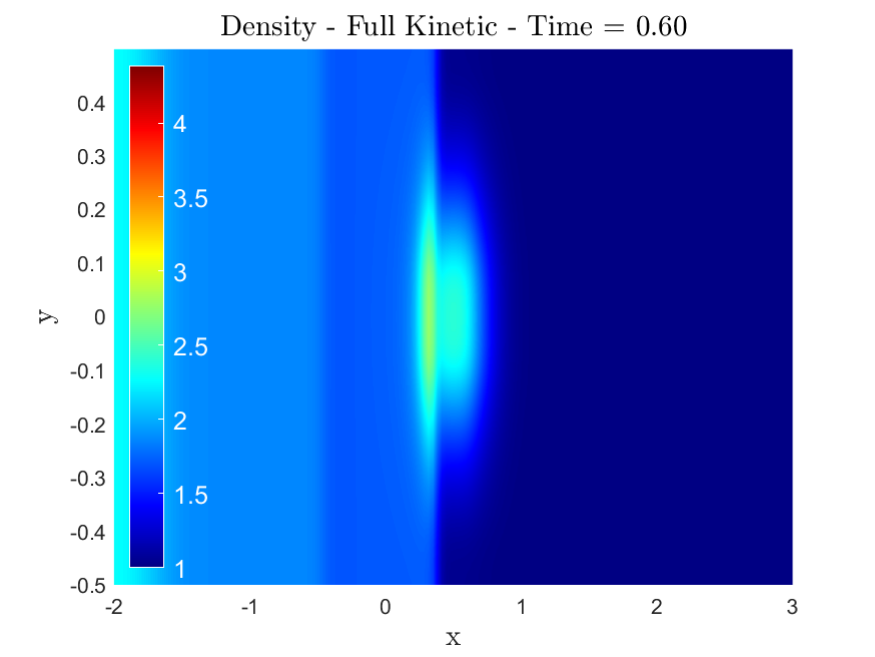}
    \end{subfigure}
    \begin{subfigure}{0.32\textwidth}
    \centering
    \includegraphics[width=\textwidth, trim={1cm 0cm 1.8cm 0cm}]{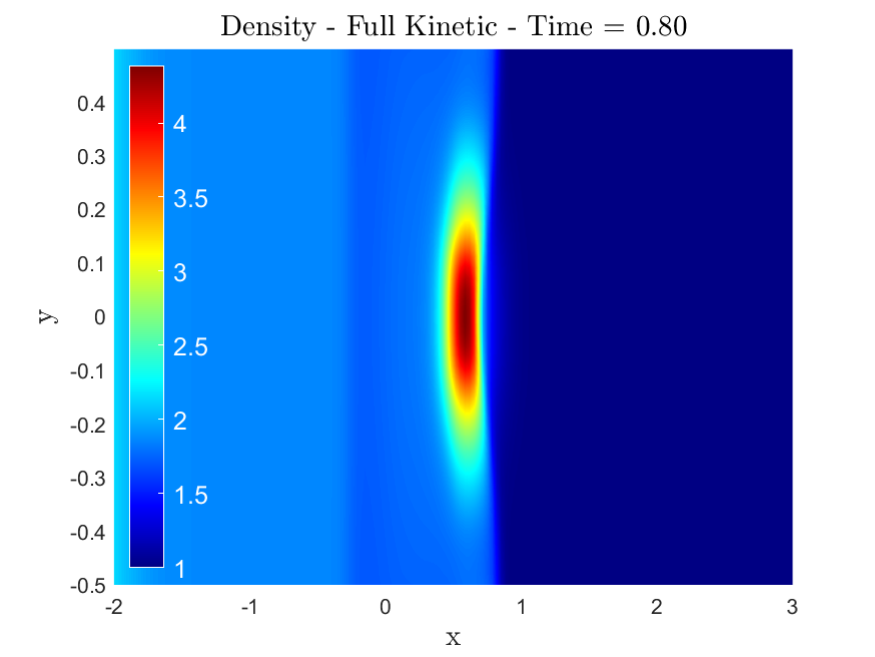}
    \end{subfigure}
    \begin{subfigure}{0.32\textwidth}
    \centering
    \includegraphics[width=\textwidth, trim={1cm 0cm 1.8cm 0cm}]{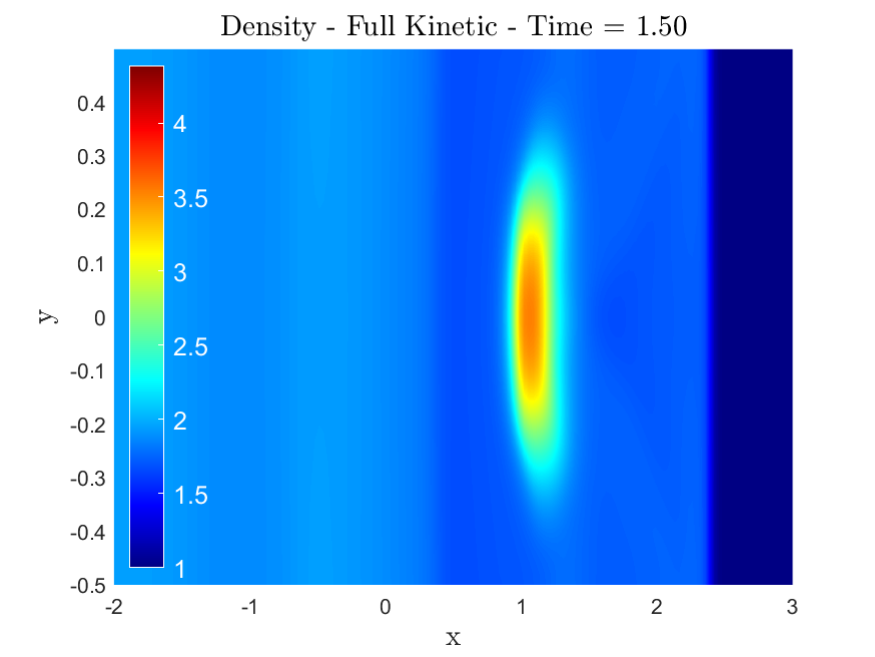}
    \end{subfigure}
    \hfill
    
    \centering
    \begin{subfigure}{0.32\textwidth}
    \centering
    \includegraphics[width=\textwidth, trim={1cm 0cm 1.8cm 0cm}]{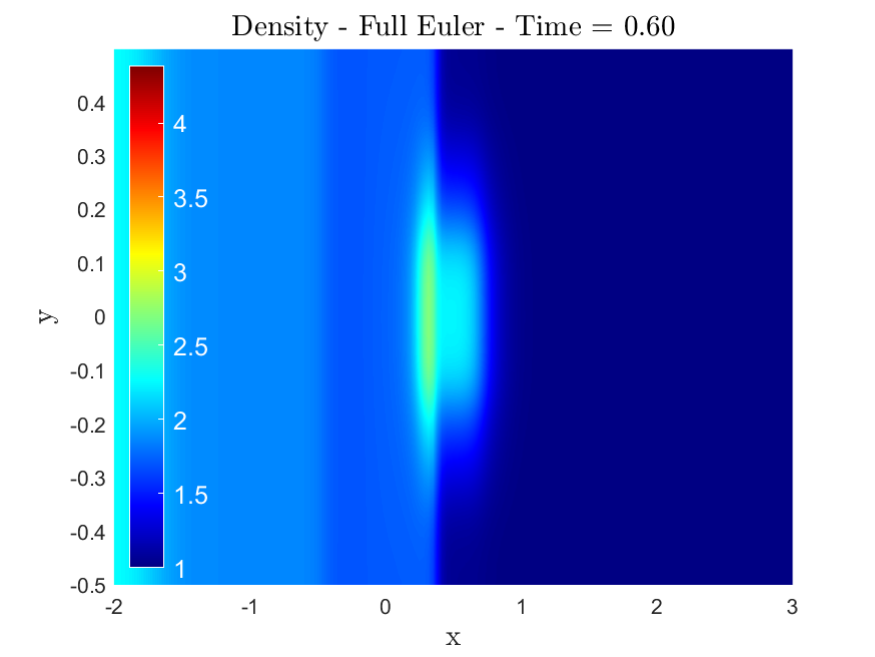}
    \end{subfigure}
    \begin{subfigure}{0.32\textwidth}
    \centering
    \includegraphics[width=\textwidth, trim={1cm 0cm 1.8cm 0cm}]{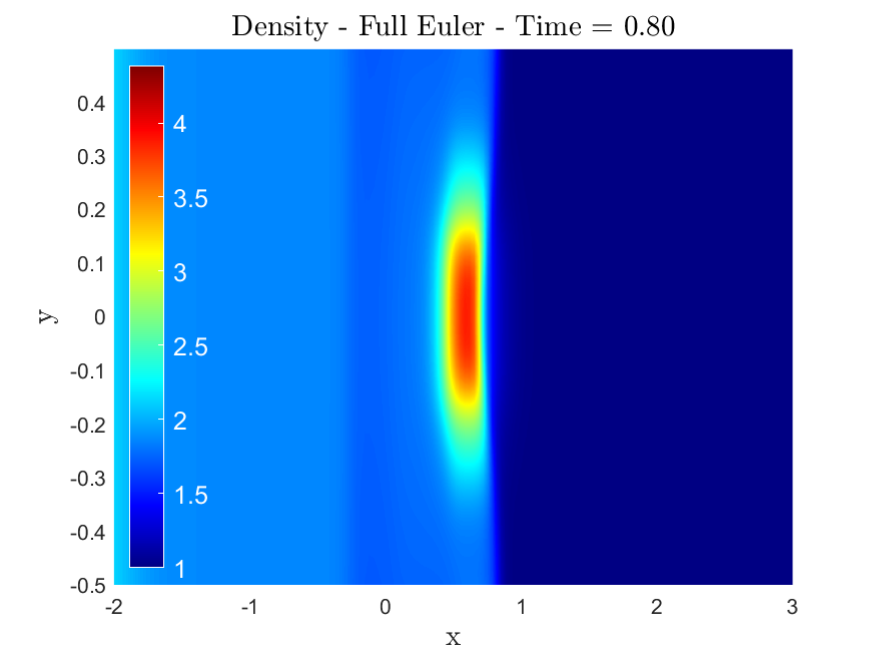}
    \end{subfigure}
    \begin{subfigure}{0.32\textwidth}
    \centering
    \includegraphics[width=\textwidth, trim={1cm 0cm 1.8cm 0cm}]{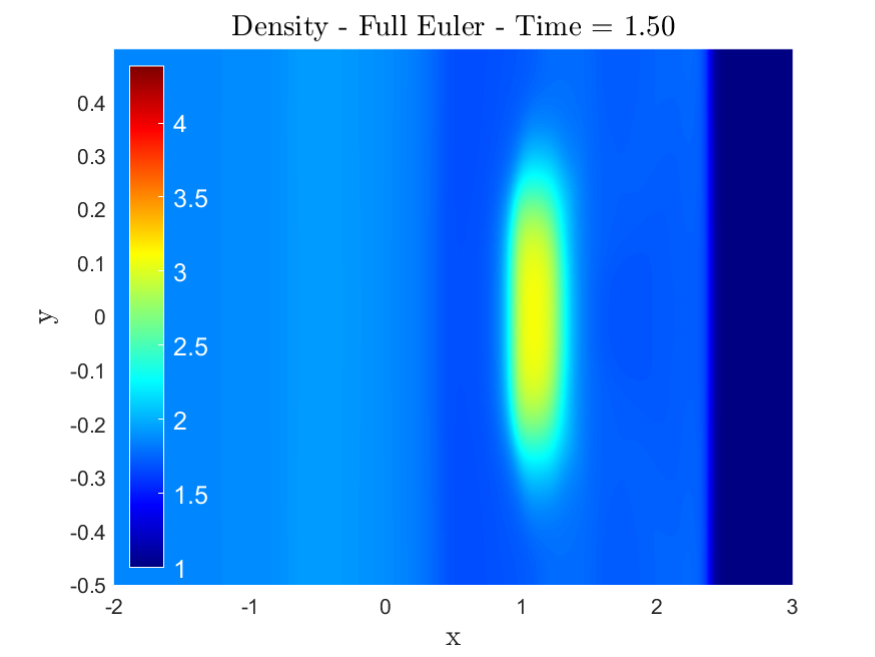}
    \end{subfigure}
    \hfill
    
    \caption{In the first row domain adaptation is shown, whereas the other rows present the densities obtained using the hybrid (second row), the full Kinetic (third row) and full Euler (fourth row) solvers, respectively, for the shock bubble interaction with initial condition \eqref{eq::initial_Bubble}. The solutions are displayed at time $t=0.60$, $t=0.80$ and $t=1.5$.}
    \label{figSB}
\end{figure}

Let us consider the Shock Bubble test to investigate the interaction between a travelling wave and a stationary bubble, inspired by \cite{torrilhon2006, melis2019}.\\
The physical domain is given by $[-2, 3] \times [-0.5, 0.5]$, discretized with $N_x=150$ cells, and $N_y=50$ cells, along the $x-$ and $y-$direction respectively. The velocity domain $\Omega_v=[-L,L]^2$, where $L=8$ is discretized using $N_v=32$ points in both directions.\\
The time step used is $\Delta t = 0.5\min{(\Delta x, \Delta y)} / L$ and we consider as final time of integration $t_{end}=1.5$.\\
Let us consider a Knudsen number $\varepsilon=10^{-6}$, which immediately places the simulation in a stiff regime.\\
The fluid is initialized with the following moments
\begin{equation}
\label{eq::initial_Bubble}
    \begin{cases}
    \begin{aligned}
		&\rho= \rho_L, \qquad &u_x=u_{x_L}, \qquad &u_y=0, \qquad &p=p_L \qquad &\text{ if } x \leq -1,\\
		
		&\rho= \rho_R, \qquad &u_x=u_{x_R}, \qquad &u_y=0, \qquad &p=p_R \qquad &\text{ if } x > -1,
		\end{aligned}
    \end{cases}
\end{equation}
where $\rho_L=16/7$, $\rho_R=1+1.5\exp{(-16d^2)}$, $u_{x_L}=\sqrt{5/3}(7/16)$, $u_{x_R}=0$, $p_L=4.75$ and $p_R=1$, and where $d=\sqrt{(x-x_0)^2+(y-y_0)^2}$ is the distance of each cell from the bubble centered in $(x_0, y_0)=(0.5, 0)$.\\
The distribution function is initialized everywhere as a Maxwellian, whose moments coincide with those of the macroscopic states precised above in \eqref{eq::initial_Bubble}.\\
The thresholds used in switching criteria (Equations \eqref{crit:CompEuler} and \eqref{eqKintoFluid1}) are $\eta=\varepsilon/2$ and $\delta=10^{-3}$.\\
In Figure \ref{figSB} we compare the density obtained using the hybrid scheme, the full Euler solver and the full Kinetic (ES-BGK) solver. In the first row we also report the domain adaptation. Similarly, in this test, the results clearly demonstrate strong agreement between the solutions, while the kinetic regime is automatically and dynamically activated only close to the discontinuity, as expected.

\subsubsection{Test 5: Fluid flux towards a cylinder}
\label{Test5}
The aim of this test is to check that the scheme is able to handle complex geometries, like a flux around a fixed cylinder, as done in \cite{fil-jin-2011}.\\
The physical domain is given by $[0, L_x] \times [0, L_y]$, where $L_x=1.0$, and $L_y=1.0$, discretized with $N_x=78$ cells, and $N_y=78$ cells, along the $x-$ and $y-$direction respectively. The velocity domain $\Omega_v=[-L,L]^2$, where $L=8$ is discretized using $N_v=32$ points in both directions.\\
The time step used is $\Delta t = 0.5\min{(\Delta x, \Delta y)} / L$ and we consider as final time of integration $t_{end}=0.15$.\\ 
We adopt again a Knudsen number $\varepsilon=10^{-6}$, which motivates the need for a stiff efficient solver.\\
The round shape, centered around $(x_c, y_c)=(2L_x/5, L_y/2)$  and of radius $R=L_x/12$ is chosen such that the integrators (Kinetic and Macroscopic) are never executed in it, and it implements specular reflective boundary conditions \cite{cercignani1994} on its surface. The fluid is initizialized elsewhere in the domain with the following macroscopic moments
\begin{equation}
\label{eq::initial_Cylinder}
	\rho= 1, \qquad u_x=2, \qquad u_y=0, \qquad p=1 \qquad \text{ if } r > R,
\end{equation}
where $r=\sqrt{(x-x_c)^2+(y-y_c)^2}$ is the distance of the point $(x, y)$ from the center of the round shape $(x_c, y_c)$.\\
The thresholds used are $\eta=10\varepsilon$ and $\delta=10^{-3}$ in switching criteria (Equations \eqref{crit:CompEuler} and \eqref{eqKintoFluid1}) .\\

\begin{figure}
    \centering
    
    \begin{subfigure}{0.32\textwidth}
    \centering
    \includegraphics[width=\textwidth, trim={0.2cm 0cm 0cm 0cm}]{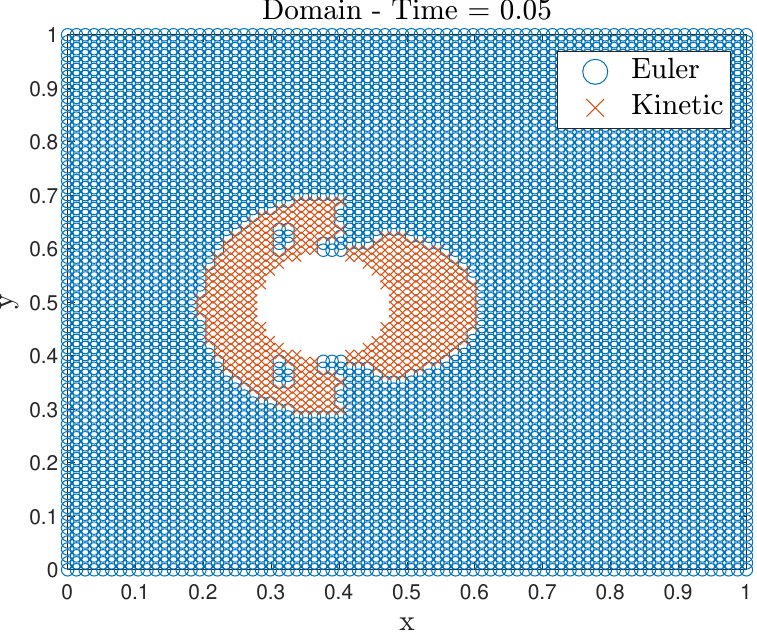}
    \end{subfigure}
    \begin{subfigure}{0.32\textwidth}
    \centering
    \includegraphics[width=\textwidth]{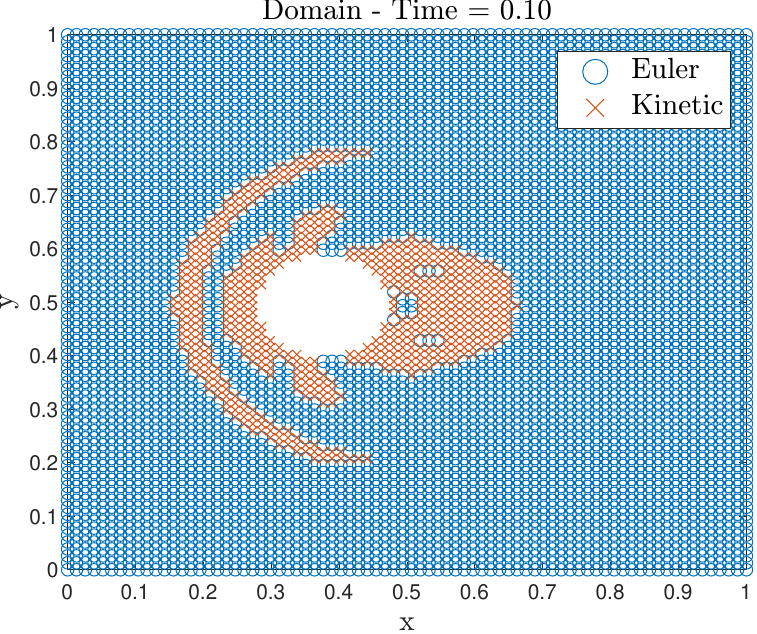}
    \end{subfigure}
    \begin{subfigure}{0.32\textwidth}
    \centering
    \includegraphics[width=\textwidth]{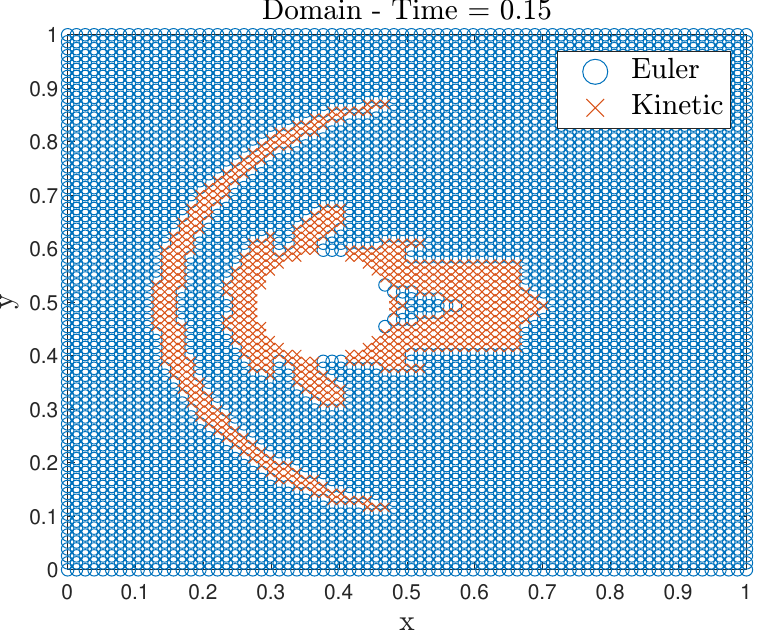}
    \end{subfigure}
    \hfill
    
    \begin{subfigure}{0.32\textwidth}
    \centering
    \includegraphics[width=\textwidth, trim={1cm 0cm 1.8cm 0cm}]{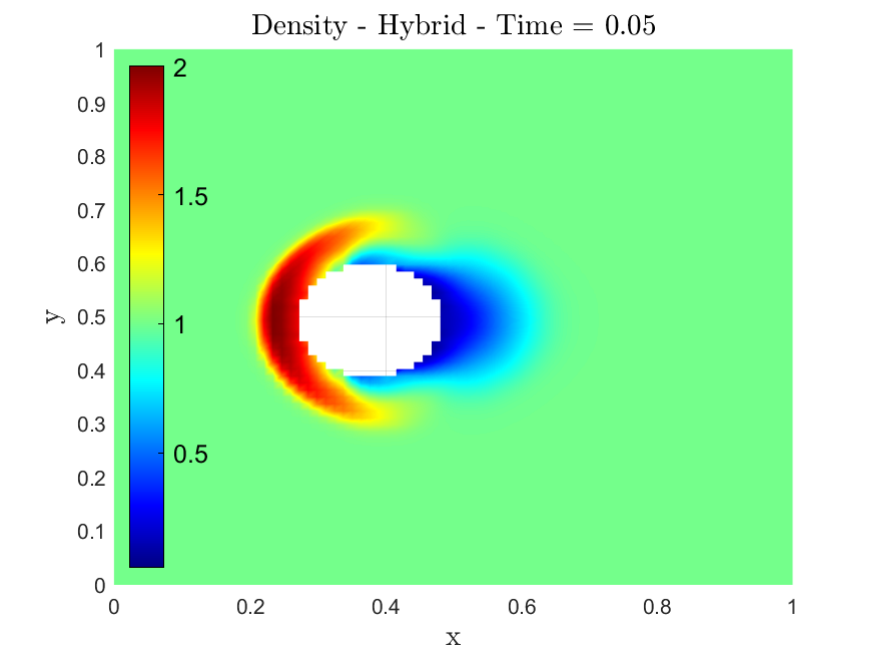}
    \end{subfigure}
    \begin{subfigure}{0.32\textwidth}
    \centering
    \includegraphics[width=\textwidth, trim={1cm 0cm 1.8cm 0cm}]{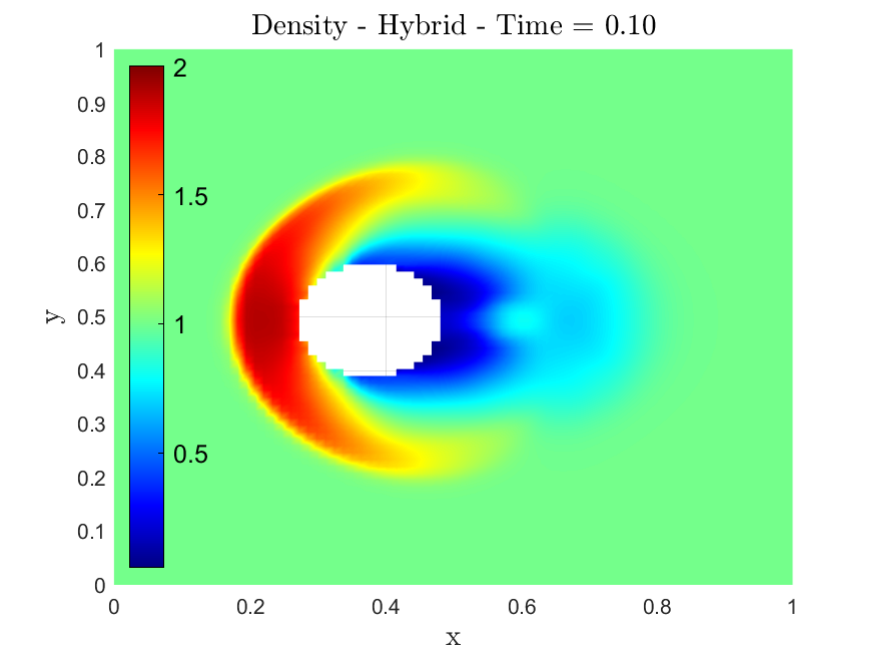}
    \end{subfigure}
    \begin{subfigure}{0.32\textwidth}
    \centering
    \includegraphics[width=\textwidth, trim={1cm 0cm 1.8cm 0cm}]{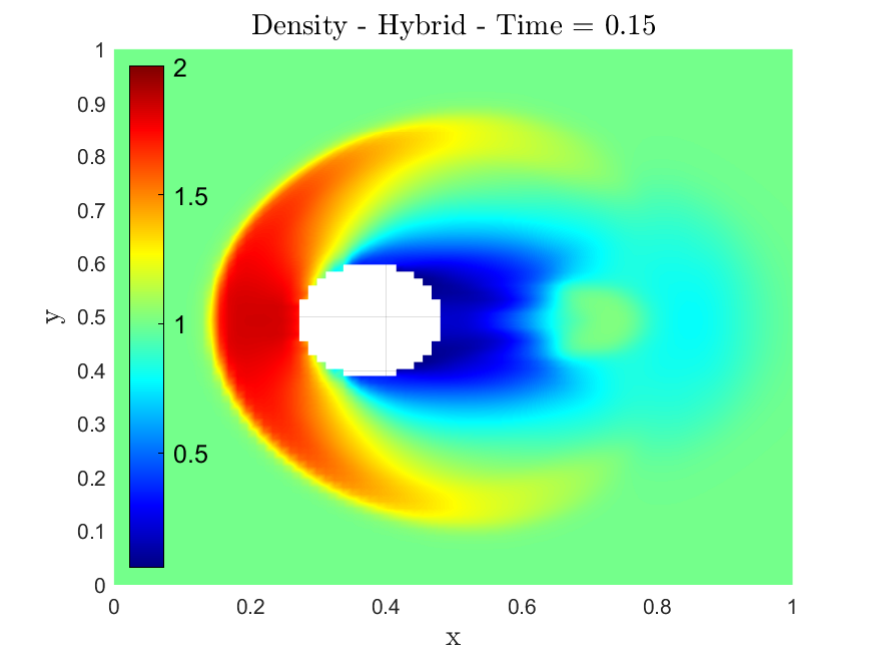}
    \end{subfigure}
    
    \begin{subfigure}{0.32\textwidth}
    \centering
    \includegraphics[width=\textwidth, trim={1cm 0cm 1.8cm 0cm}]{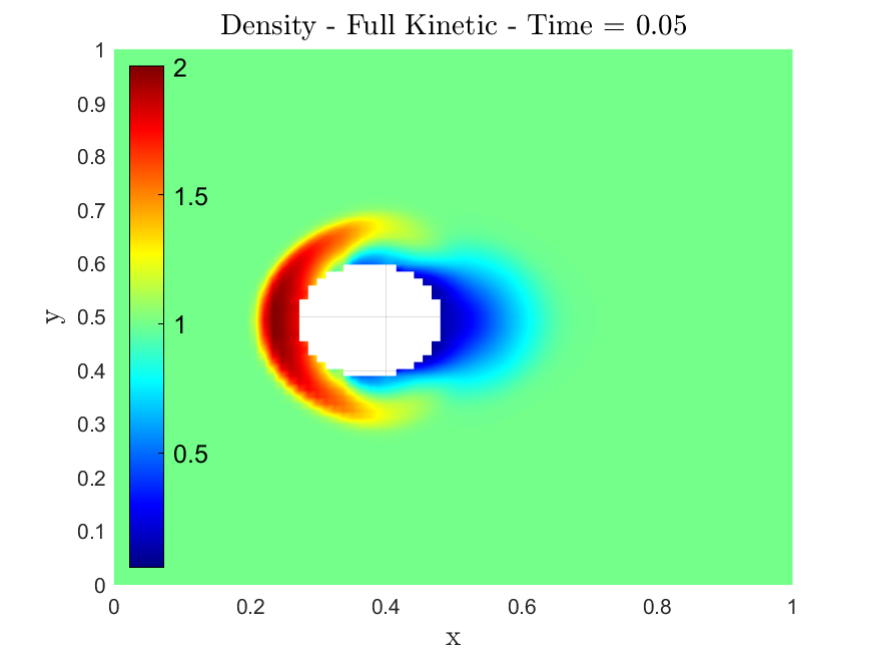}
    \end{subfigure}
    \begin{subfigure}{0.32\textwidth}
    \centering
    \includegraphics[width=\textwidth, trim={1cm 0cm 1.8cm 0cm}]{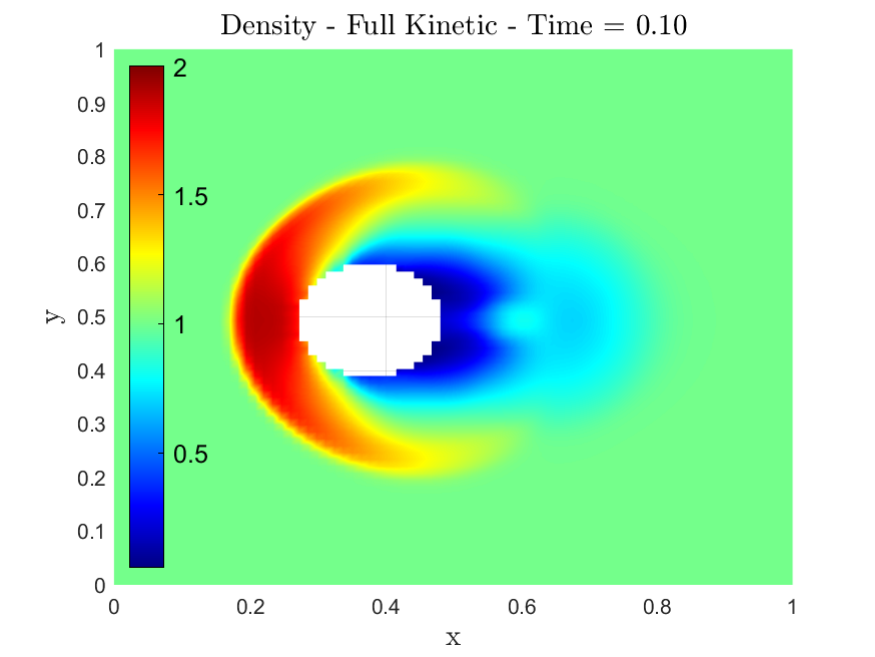}
    \end{subfigure}
    \begin{subfigure}{0.32\textwidth}
    \centering
    \includegraphics[width=\textwidth, trim={1cm 0cm 1.8cm 0cm}]{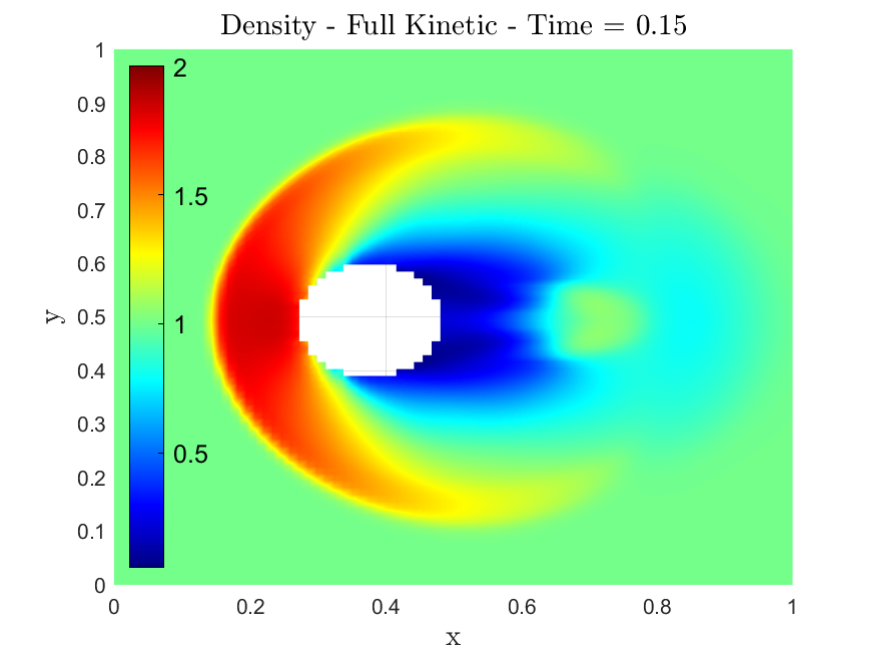}
    \end{subfigure}
    
    \begin{subfigure}{0.32\textwidth}
    \centering
    \includegraphics[width=\textwidth, trim={1cm 0cm 1.8cm 0cm}]{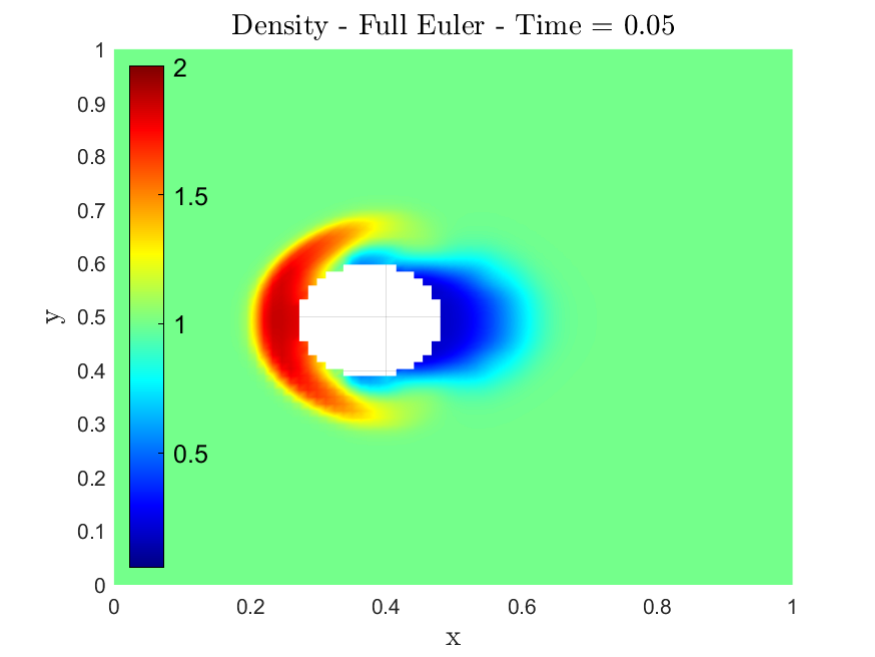}
    \end{subfigure}
    \begin{subfigure}{0.32\textwidth}
    \centering
    \includegraphics[width=\textwidth, trim={1cm 0cm 1.8cm 0cm}]{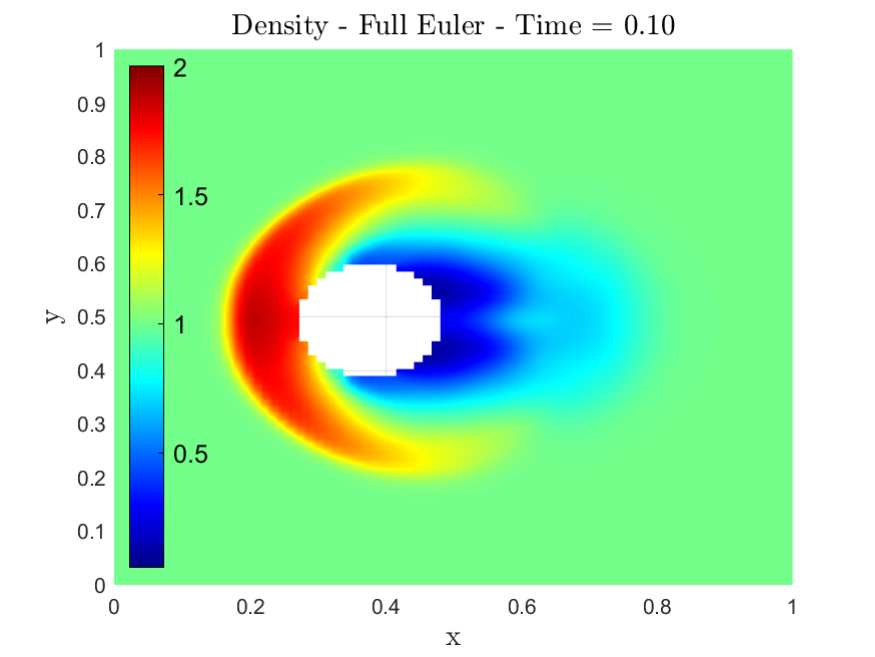}
    \end{subfigure}
    \begin{subfigure}{0.32\textwidth}
    \centering
    \includegraphics[width=\textwidth, trim={1cm 0cm 1.8cm 0cm}]{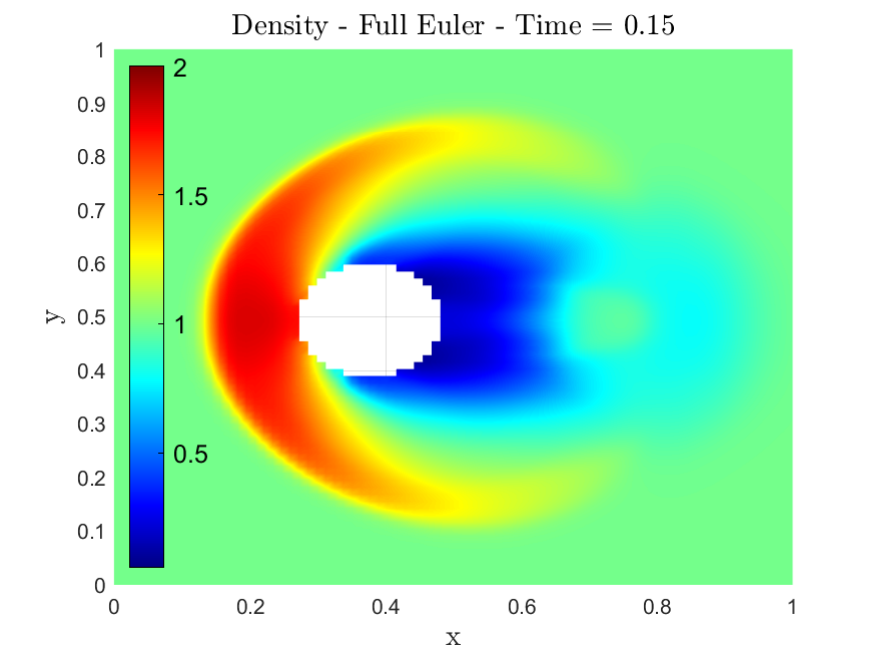}
    \end{subfigure}
    \caption{In the first row domain adaptation is shown, whereas the other rows present the densities obtained using the hybrid (second row), the full Kinetic (third row) and full Euler (fourth row) solvers, respectively, for the fluid flux towards a cylinder with initial condition \eqref{eq::initial_Cylinder}. The solutions are displayed at time $t=0.05$, $t=0.10$ and $t=0.15$.}
    \label{figCyl_density}
\end{figure}

\begin{figure}
    \centering
    
    \begin{subfigure}{0.32\textwidth}
    \centering
    \includegraphics[width=\textwidth, trim={1cm 0cm 1.8cm 0cm}]{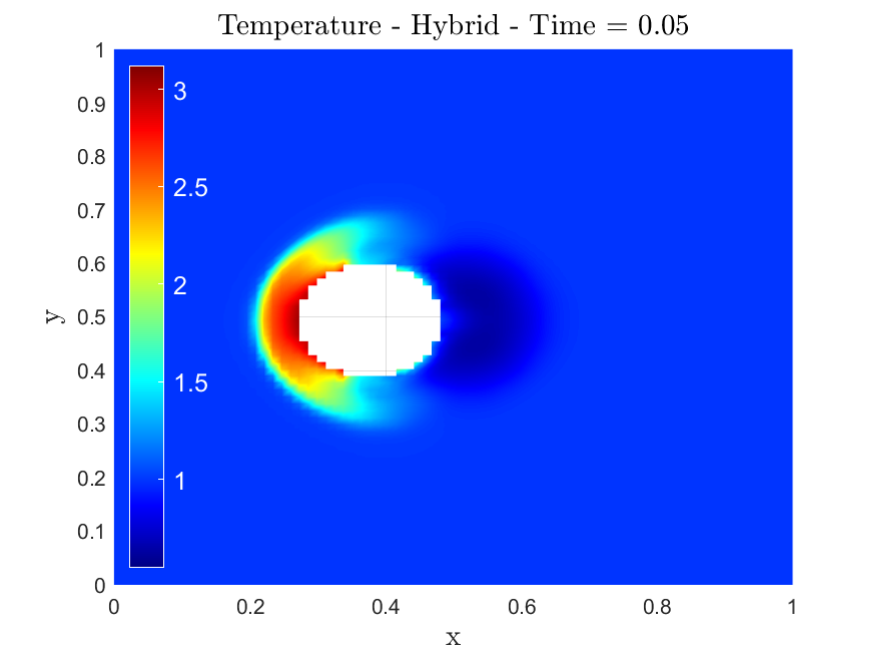}
    \end{subfigure}
    \begin{subfigure}{0.32\textwidth}
    \centering
    \includegraphics[width=\textwidth, trim={1cm 0cm 1.8cm 0cm}]{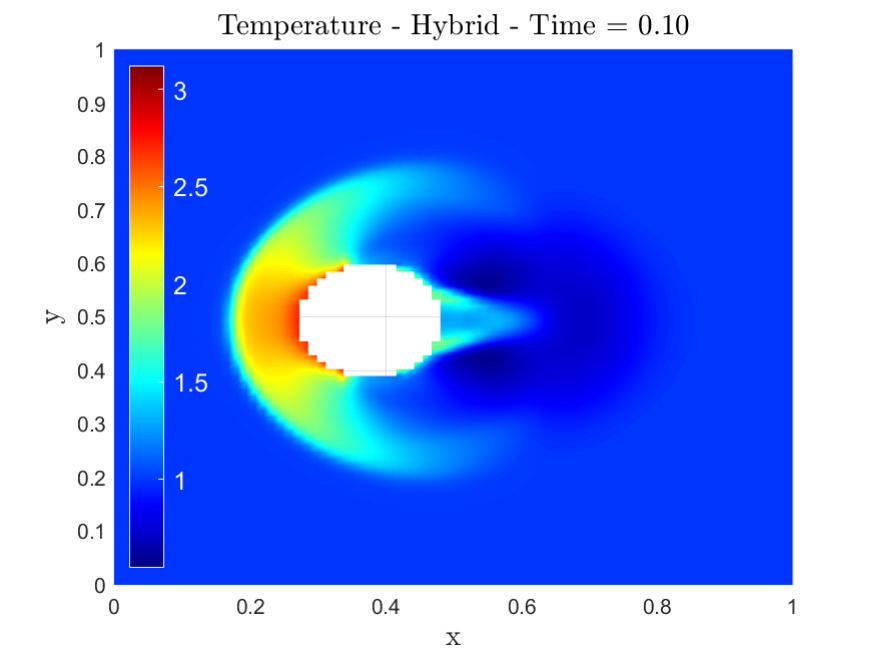}
    \end{subfigure}
    \begin{subfigure}{0.32\textwidth}
    \centering
    \includegraphics[width=\textwidth, trim={1cm 0cm 1.8cm 0cm}]{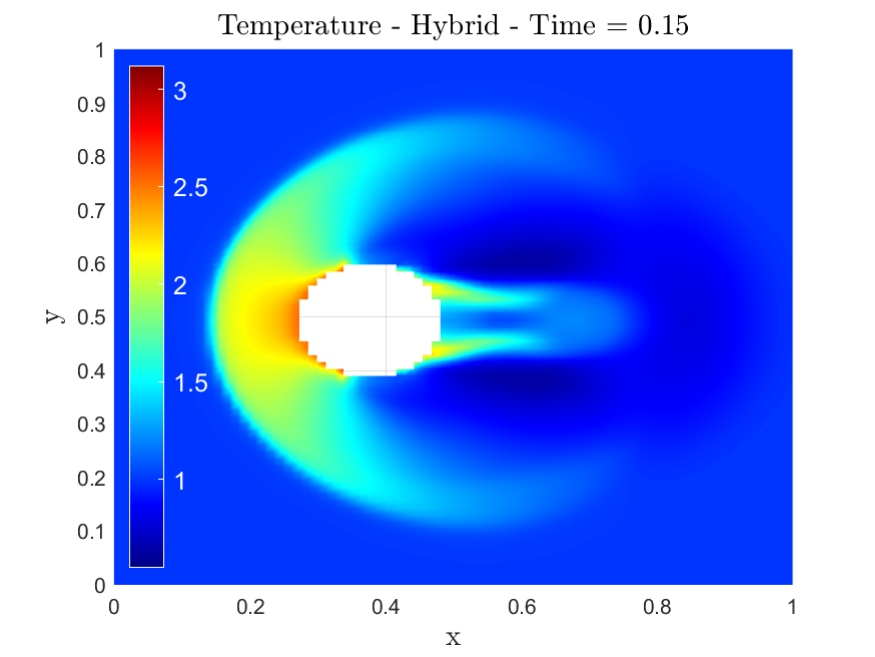}
    \end{subfigure}
    
    \begin{subfigure}{0.32\textwidth}
    \centering
    \includegraphics[width=\textwidth, trim={1cm 0cm 1.8cm 0cm}]{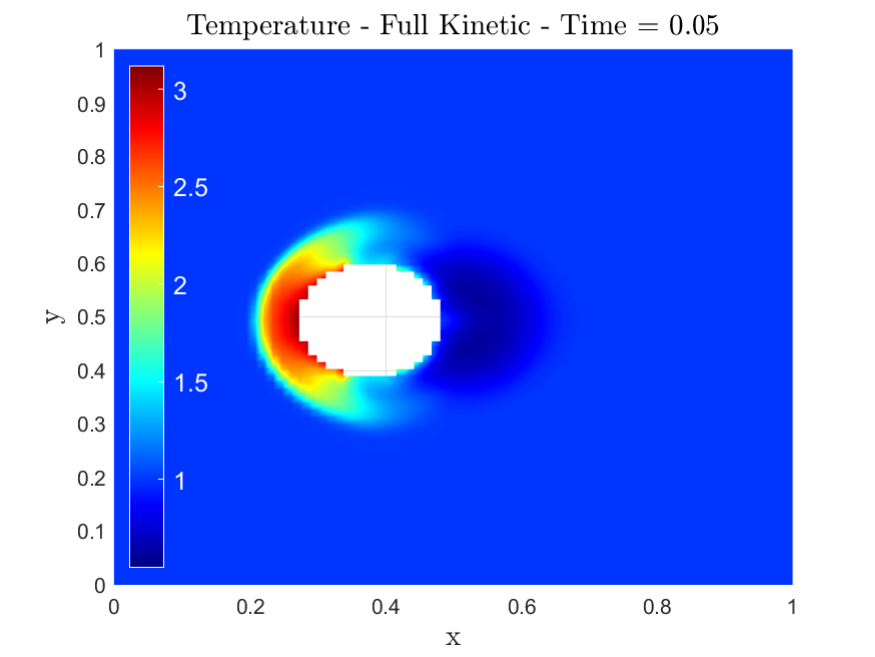}
    \end{subfigure}
    \begin{subfigure}{0.32\textwidth}
    \centering
    \includegraphics[width=\textwidth, trim={1cm 0cm 1.8cm 0cm}]{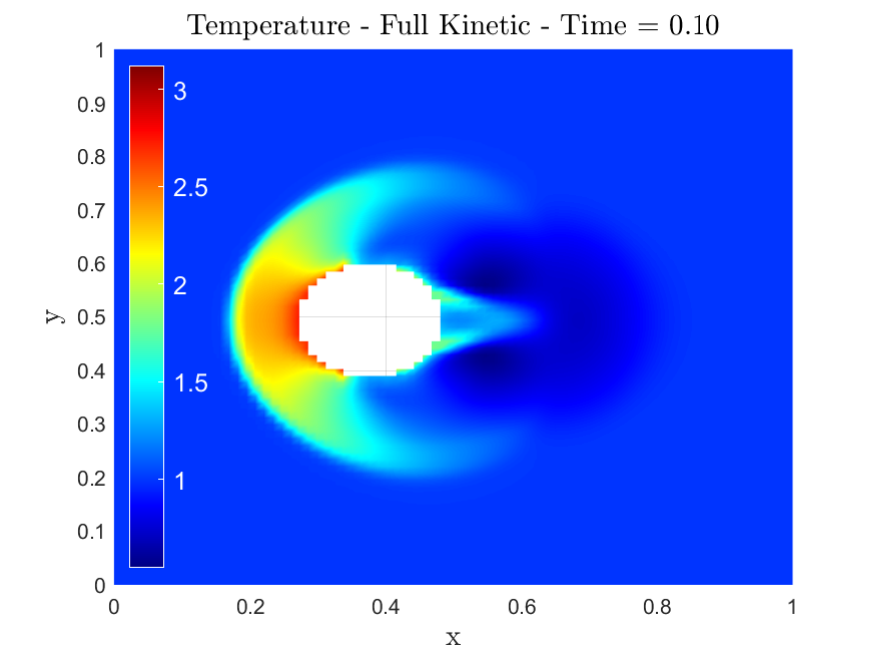}
    \end{subfigure}
    \begin{subfigure}{0.32\textwidth}
    \centering
    \includegraphics[width=\textwidth, trim={1cm 0cm 1.8cm 0cm}]{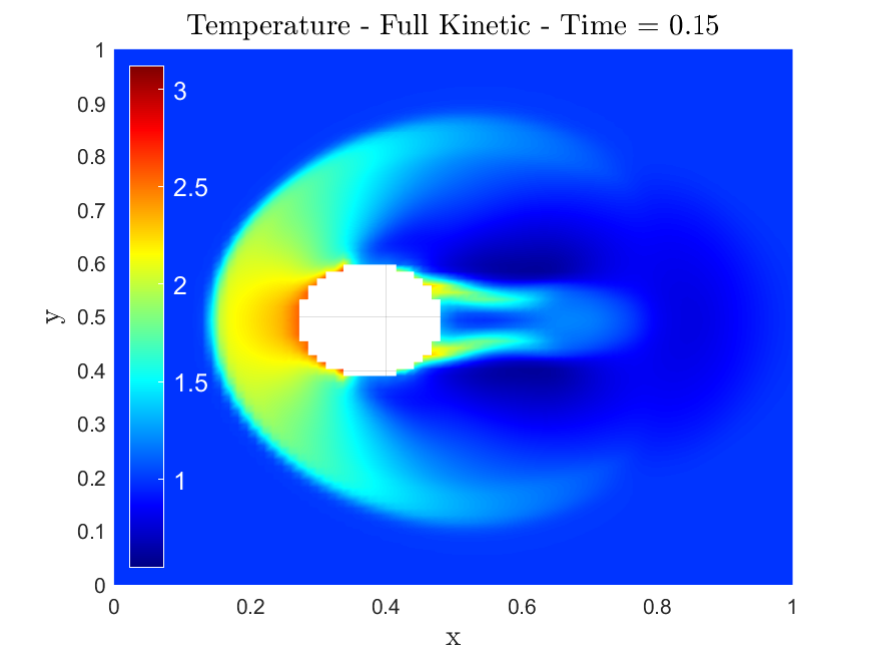}
    \end{subfigure}
    
    \begin{subfigure}{0.32\textwidth}
    \centering
    \includegraphics[width=\textwidth, trim={1cm 0cm 1.8cm 0cm}]{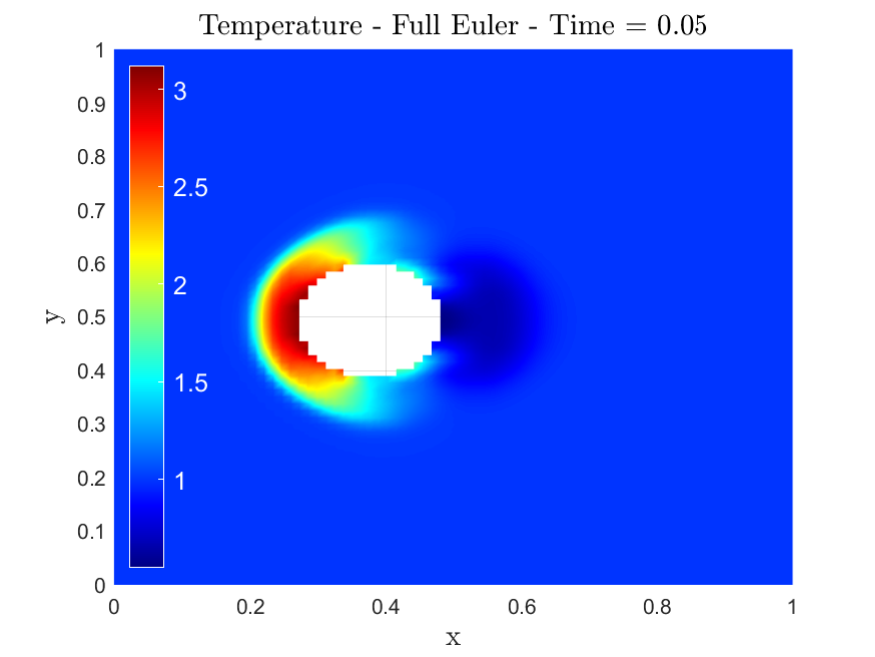}
    \end{subfigure}
    \begin{subfigure}{0.32\textwidth}
    \centering
    \includegraphics[width=\textwidth, trim={1cm 0cm 1.8cm 0cm}]{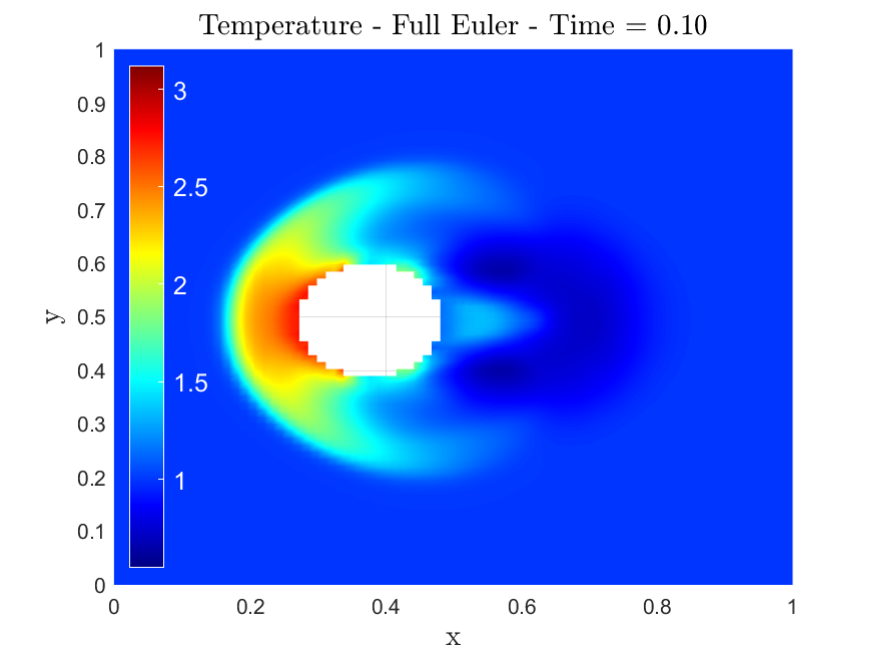}
    \end{subfigure}
    \begin{subfigure}{0.32\textwidth}
    \centering
    \includegraphics[width=\textwidth, trim={1cm 0cm 1.8cm 0cm}]{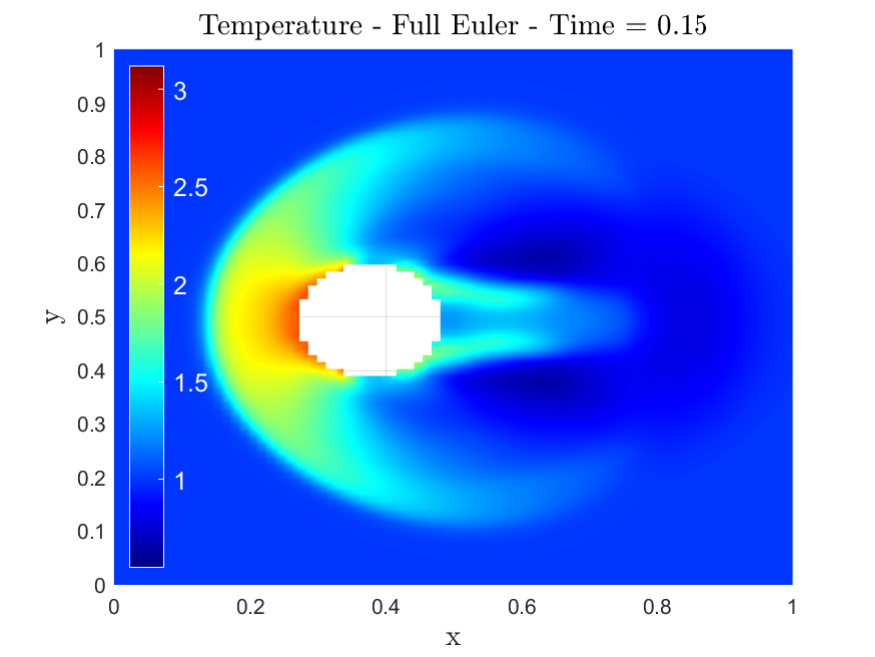}
    \end{subfigure}
    \caption{Temperatures obtained using the hybrid (first row), the full Kinetic (second row) and full Euler (third row) solvers, are shown, for the fluid flux towards a cylinder with initial condition \eqref{eq::initial_Cylinder}. The solutions are displayed at time $t=0.05$, $t=0.10$ and $t=0.15$.}
    \label{figCyl_temperature}
\end{figure}

In Figure \ref{figCyl_density} we compare the density obtained using the hybrid scheme, the full Euler solver and the full Kinetic (ES-BGK) solver. In the first row we show again the domain adaptation. The solutions are clearly strongly consistent with each other also for this test, indicating that the method accurately captures the expected behavior across the domain. The kinetic regime is automatically and dynamically triggered only close to the discontinuity, which in this case is generated by the presence of the obstacle in the flow. Around the cylinder, non-equilibrium effects appear, requiring the kinetic model, while the rest of the domain is still well represented by the Euler macroscopic model. In Figure \ref{figCyl_temperature} are reported the same quantities but for the temperature.

\subsubsection{Test 6: Fluid flux towards a rectangular obstacle}
The aim of this test is to check that the scheme is able to handle complex geometries, like a flux around a fixed rectangular object \cite{caparello2025}.\\
The physical domain is given by $[0, L_x] \times [0, L_y]$, where $L_x=1.0$, and $L_y=1.0$, discretized with $N_x=78$ cells, and $N_y=78$ cells, along the $x-$ and $y-$direction respectively. The velocity domain $\Omega_v=[-L,L]^2$, where $L=8$ is discretized using 32 points in both directions.\\
The time step used is $\Delta t = 0.5\min{(\Delta x, \Delta y)} / L$ and the final time of integration is $t_{end}=0.15$.\\ 
For this simulation, we consider a Knudsen number $\varepsilon=10^{-6}$.\\
The rectangular object, centered in $(x_c, y_c)=(2L_x/5, L_y/2)$  of dimension $10\Delta x \times 4\Delta y$ is chosen such that the integrators (Kinetic and Macroscopic) are never executed in it, and it implements specular reflective boundary conditions \cite{cercignani1994} on its surface. The fluid is initialized elsewhere in the domain (outside of the rectangular object) with the following moments
\begin{equation}
\label{eq::initial_Rectangular}
	\rho= 1, \qquad u_x=2, \qquad u_y=0, \qquad p=1.
\end{equation}
The thresholds used in switching criteria (Equations \eqref{crit:CompEuler} and \eqref{eqKintoFluid1}) are $\eta=10\varepsilon$ and $\delta=10^{-3}$.\\

\begin{figure}
    \centering
    
    \begin{subfigure}{0.32\textwidth}
    \centering
    \includegraphics[width=\textwidth, trim={0.2cm 0cm 0cm 0cm}]{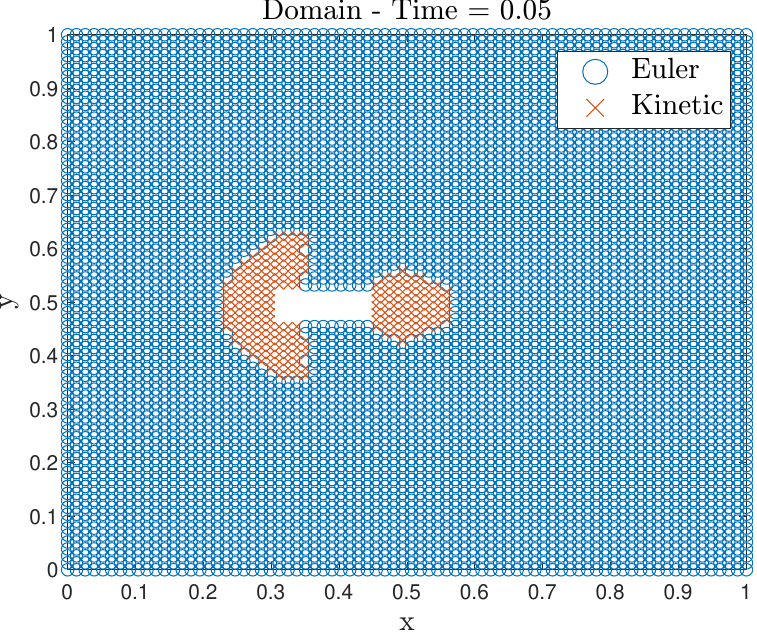}
    \end{subfigure}
    \begin{subfigure}{0.32\textwidth}
    \centering
    \includegraphics[width=\textwidth, trim={0.2cm 0cm 0cm 0cm}]{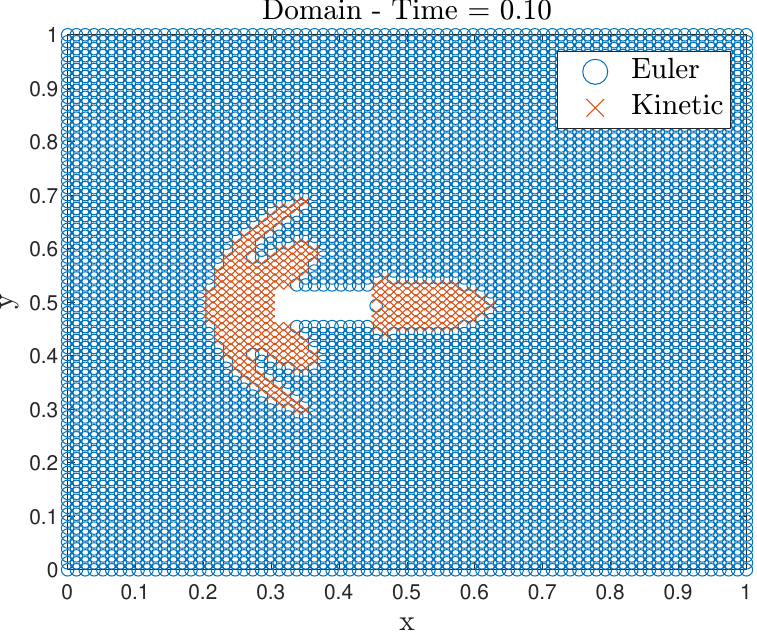}
    \end{subfigure}
    \begin{subfigure}{0.32\textwidth}
    \centering
    \includegraphics[width=\textwidth, trim={0.2cm 0cm 0cm 0cm}]{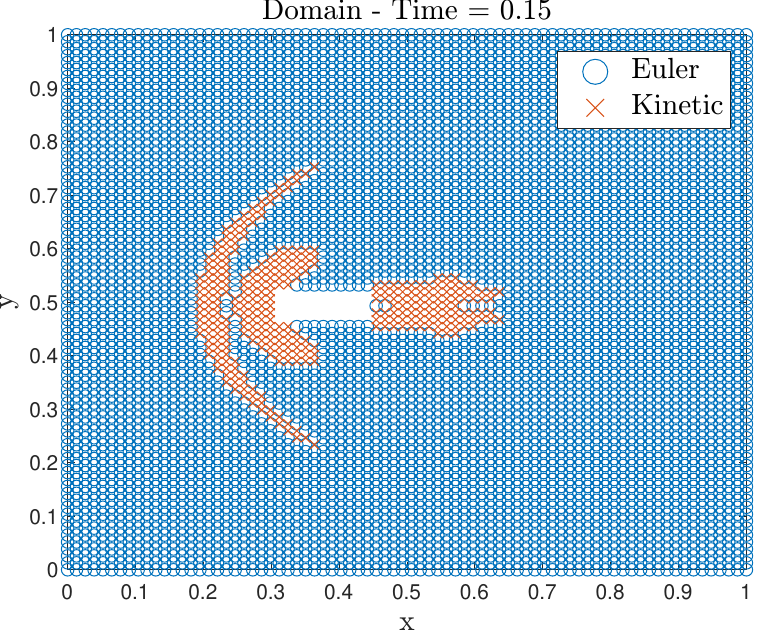}
    \end{subfigure}
    \hfill
    
    \begin{subfigure}{0.32\textwidth}
    \centering
    \includegraphics[width=\textwidth, trim={1cm 0cm 1.8cm 0cm}]{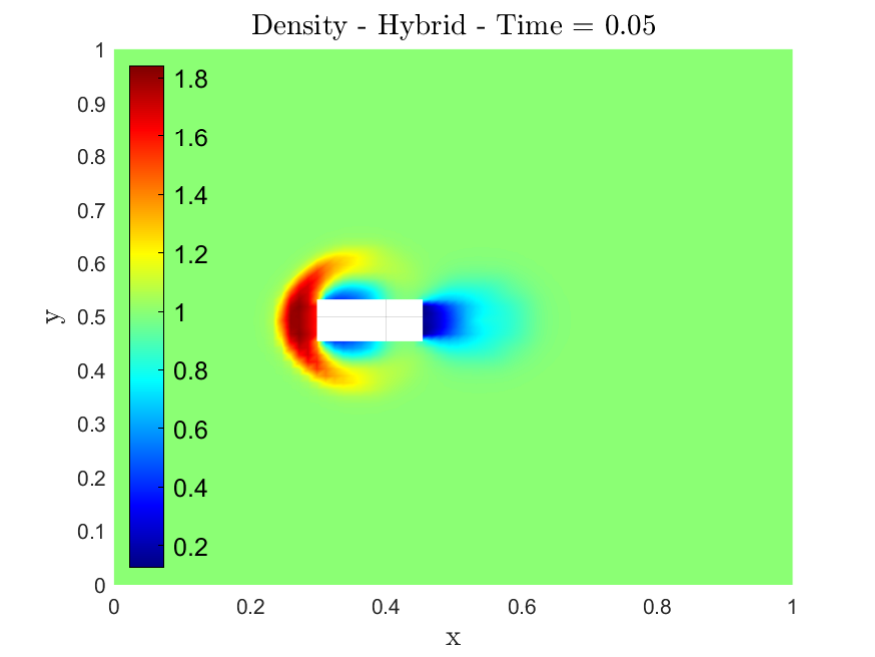}
    \end{subfigure}
    \begin{subfigure}{0.32\textwidth}
    \centering
    \includegraphics[width=\textwidth, trim={1cm 0cm 1.8cm 0cm}]{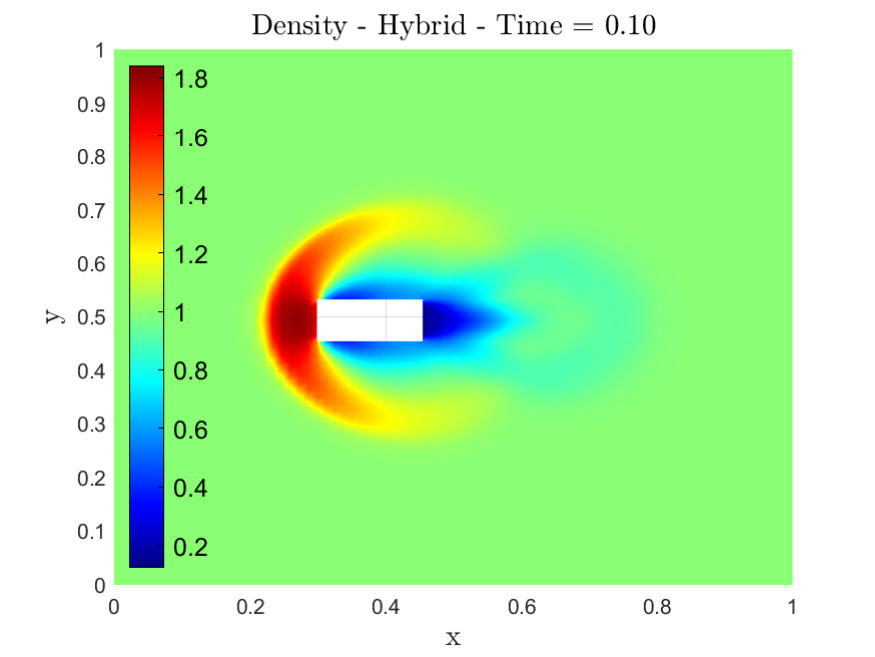}
    \end{subfigure}
    \begin{subfigure}{0.32\textwidth}
    \centering
    \includegraphics[width=\textwidth, trim={1cm 0cm 1.8cm 0cm}]{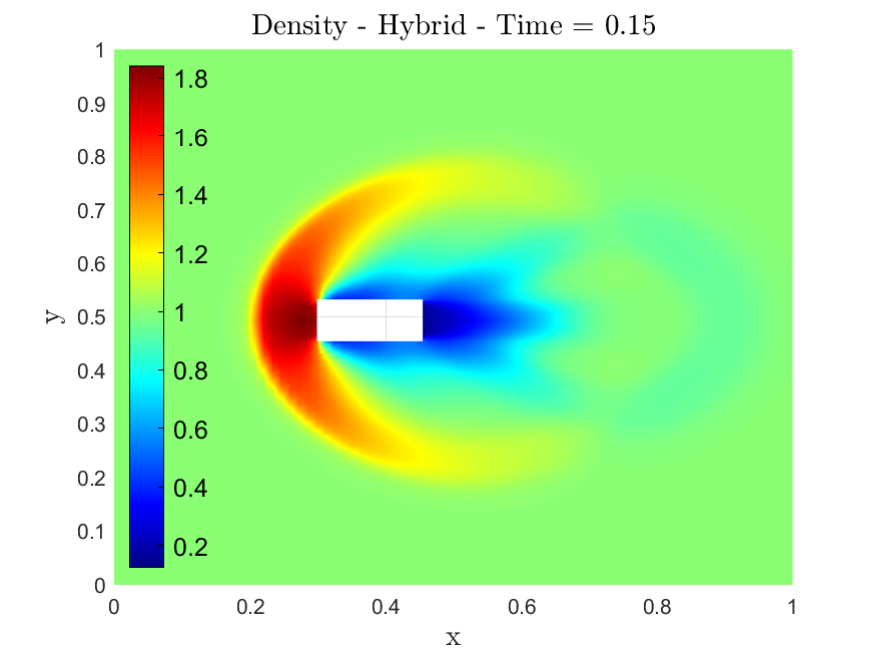}
    \end{subfigure}
    
    \begin{subfigure}{0.32\textwidth}
    \centering
    \includegraphics[width=\textwidth, trim={1cm 0cm 1.8cm 0cm}]{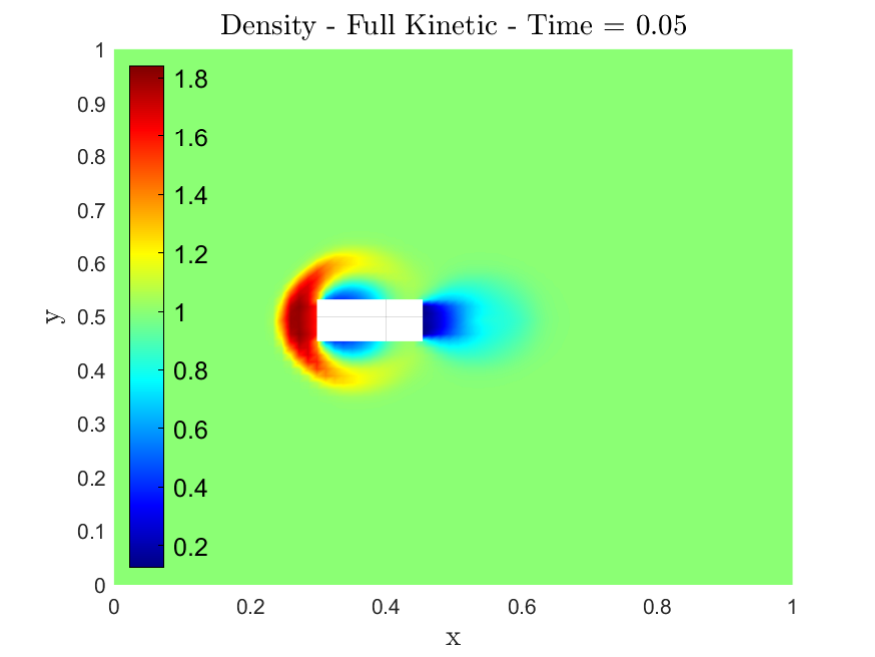}
    \end{subfigure}
    \begin{subfigure}{0.32\textwidth}
    \centering
    \includegraphics[width=\textwidth, trim={1cm 0cm 1.8cm 0cm}]{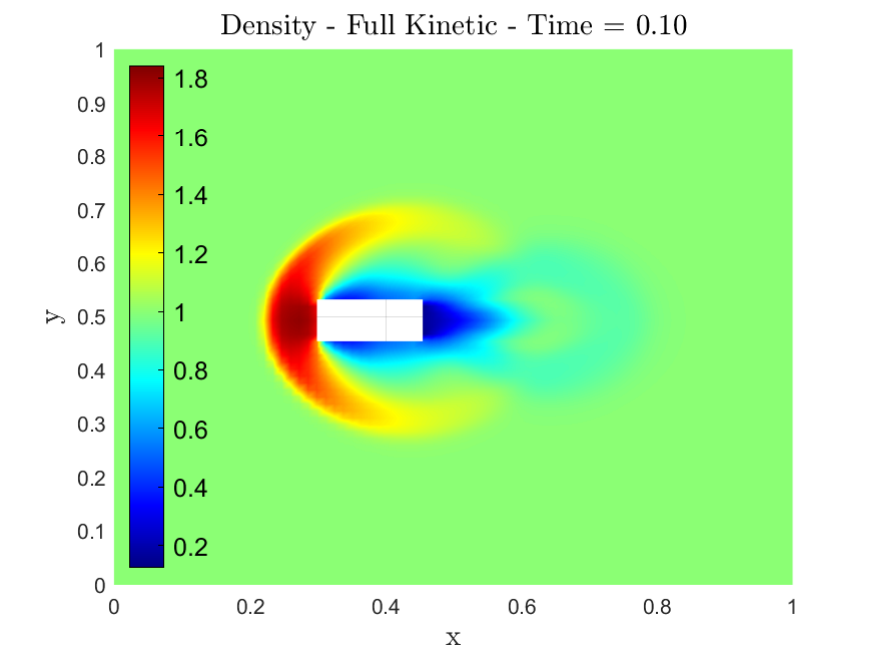}
    \end{subfigure}
    \begin{subfigure}{0.32\textwidth}
    \centering
    \includegraphics[width=\textwidth, trim={1cm 0cm 1.8cm 0cm}]{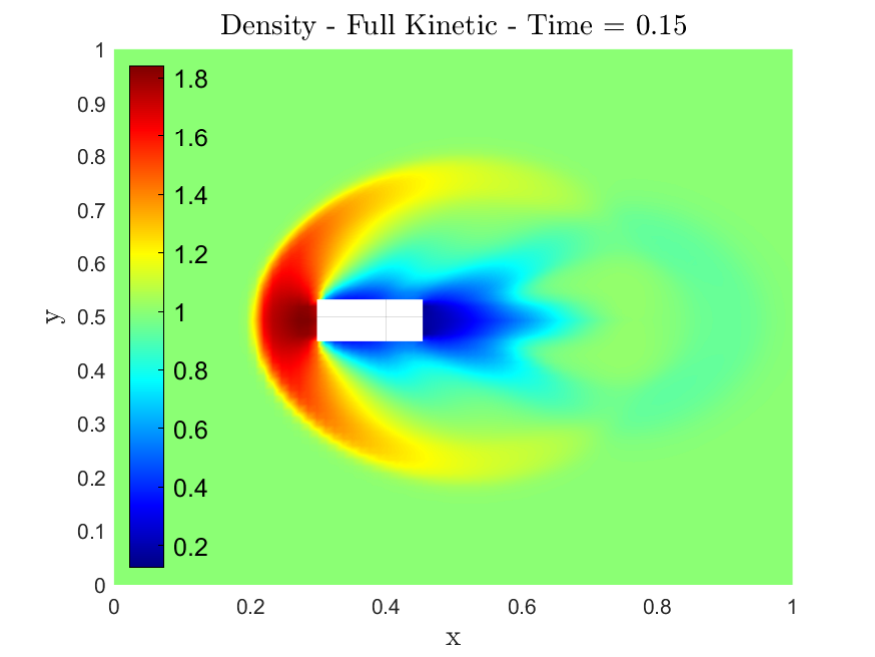}
    \end{subfigure}
    
    \begin{subfigure}{0.32\textwidth}
    \centering
    \includegraphics[width=\textwidth, trim={1cm 0cm 1.8cm 0cm}]{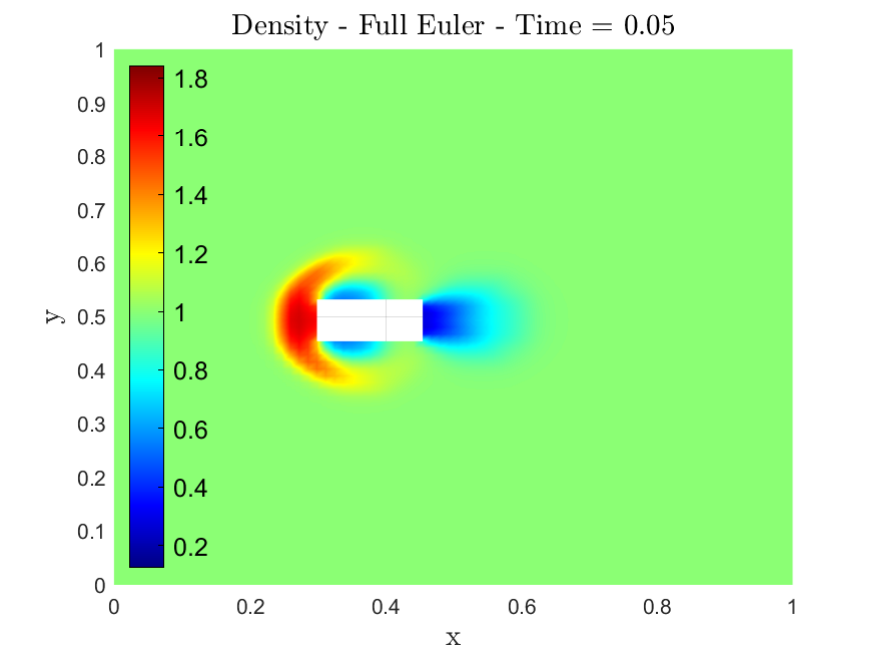}
    \end{subfigure}
    \begin{subfigure}{0.32\textwidth}
    \centering
    \includegraphics[width=\textwidth, trim={1cm 0cm 1.8cm 0cm}]{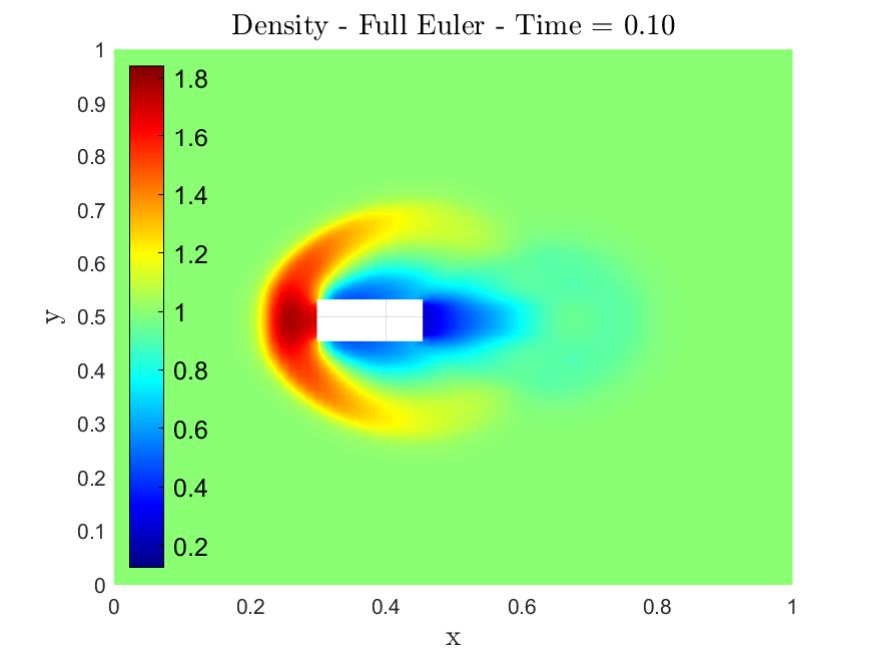}
    \end{subfigure}
    \begin{subfigure}{0.32\textwidth}
    \centering
    \includegraphics[width=\textwidth, trim={1cm 0cm 1.8cm 0cm}]{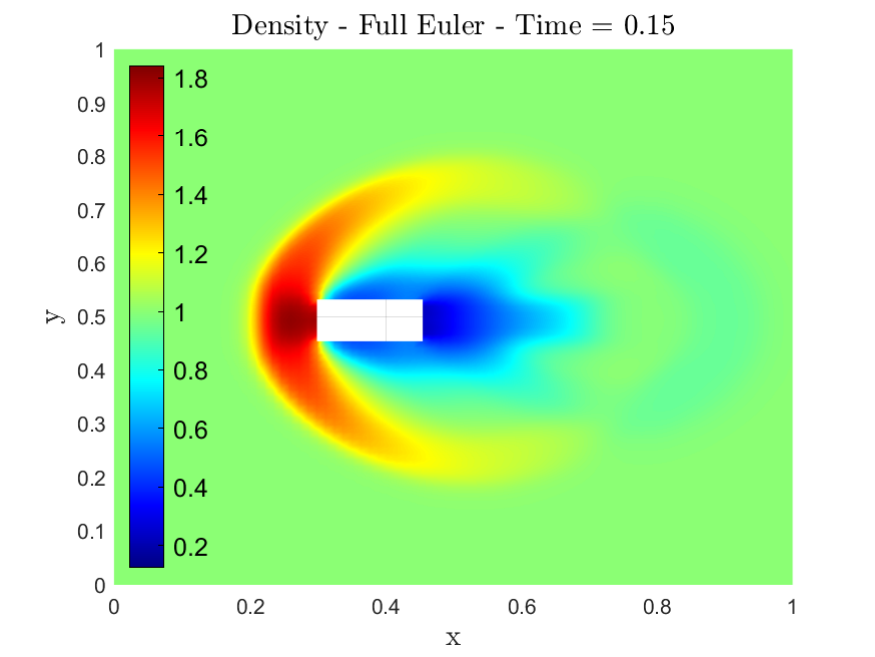}
    \end{subfigure}
    \caption{In the first row domain adaptation is shown, whereas the other rows present the densities obtained using the hybrid (second row), the full Kinetic (third row) and full Euler (fourth row) solvers, respectively, for the fluid flux towards a fixed rectangular obstacle with initial condition \eqref{eq::initial_Rectangular}. The solutions are displayed at time $t=0.05$, $t=0.10$ and $t=0.15$.}
    \label{figRec_density}
\end{figure}

\begin{figure}
    \centering
    
    \begin{subfigure}{0.32\textwidth}
    \centering
    \includegraphics[width=\textwidth, trim={1cm 0cm 1.8cm 0cm}]{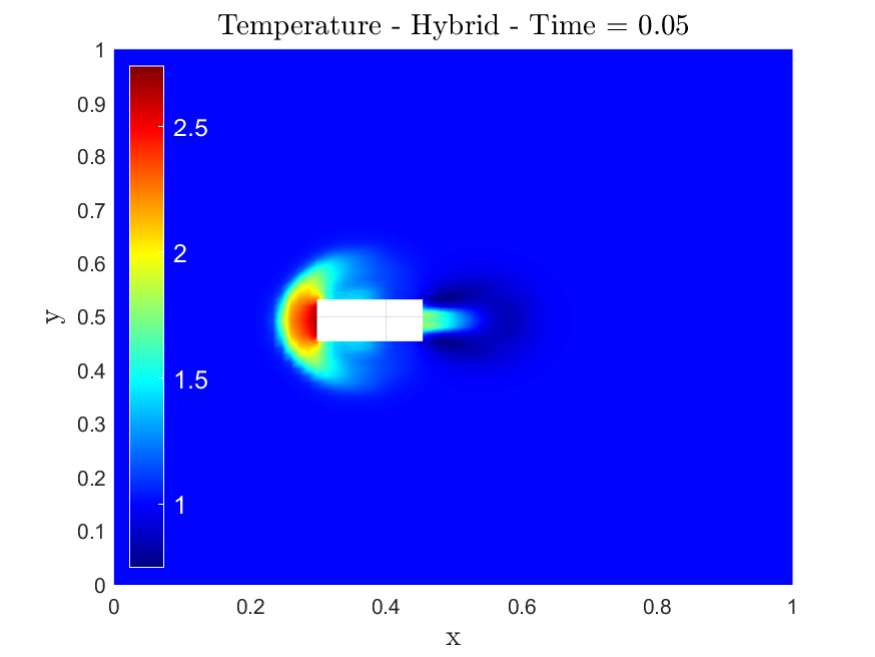}
    \end{subfigure}
    \begin{subfigure}{0.32\textwidth}
    \centering
    \includegraphics[width=\textwidth, trim={1cm 0cm 1.8cm 0cm}]{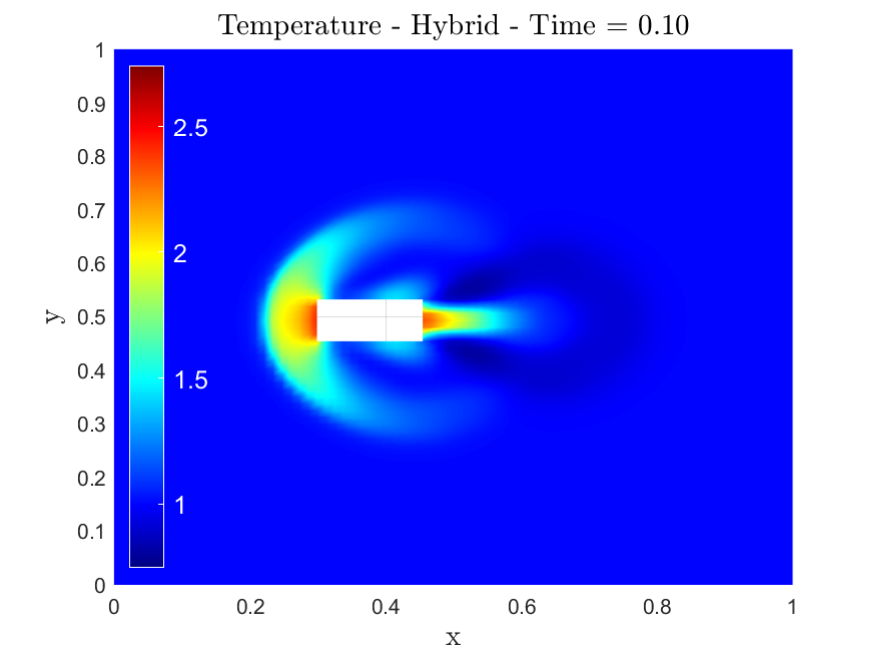}
    \end{subfigure}
    \begin{subfigure}{0.32\textwidth}
    \centering
    \includegraphics[width=\textwidth, trim={1cm 0cm 1.8cm 0cm}]{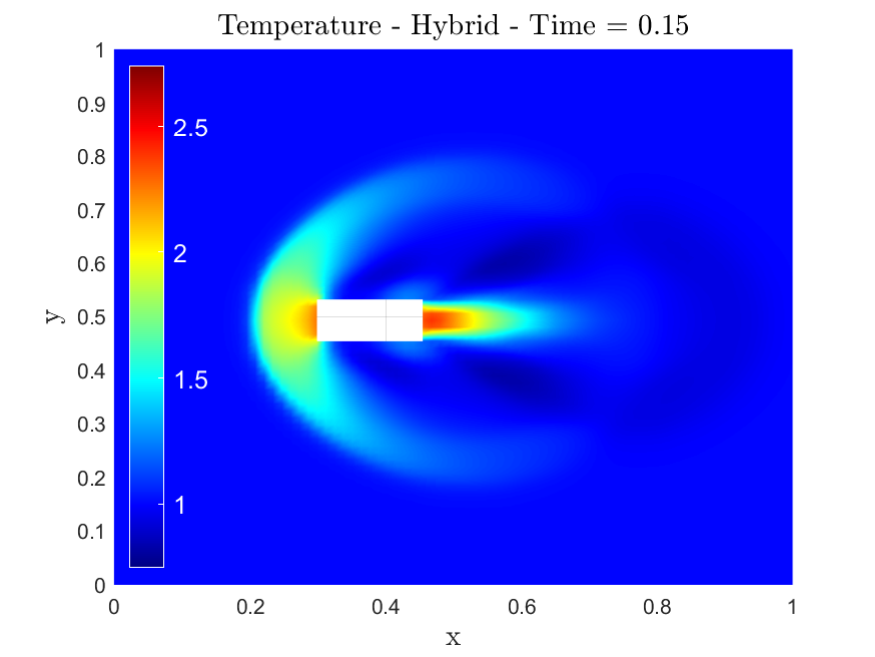}
    \end{subfigure}
    
    \begin{subfigure}{0.32\textwidth}
    \centering
    \includegraphics[width=\textwidth, trim={1cm 0cm 1.8cm 0cm}]{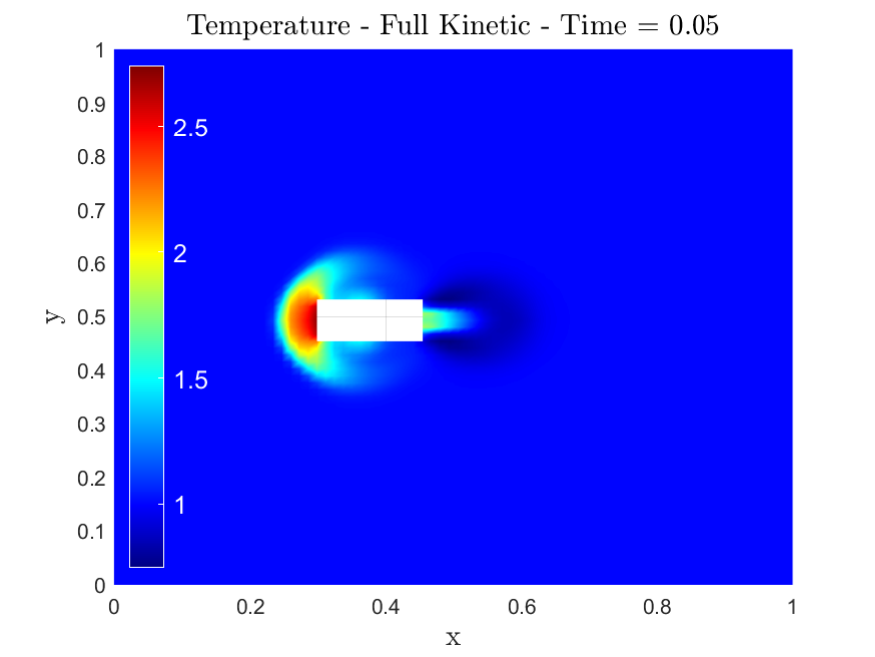}
    \end{subfigure}
    \begin{subfigure}{0.32\textwidth}
    \centering
    \includegraphics[width=\textwidth, trim={1cm 0cm 1.8cm 0cm}]{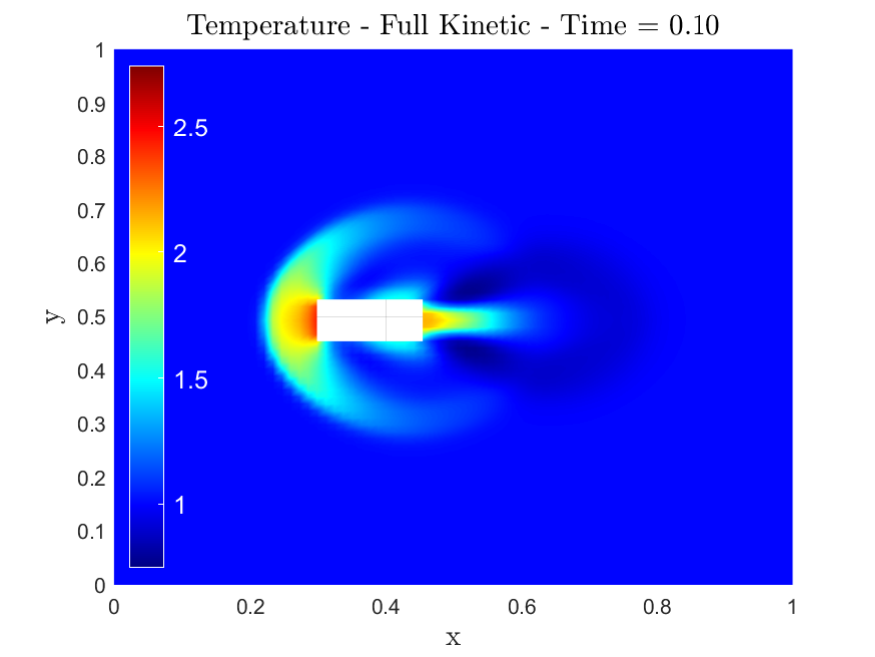}
    \end{subfigure}
    \begin{subfigure}{0.32\textwidth}
    \centering
    \includegraphics[width=\textwidth, trim={1cm 0cm 1.8cm 0cm}]{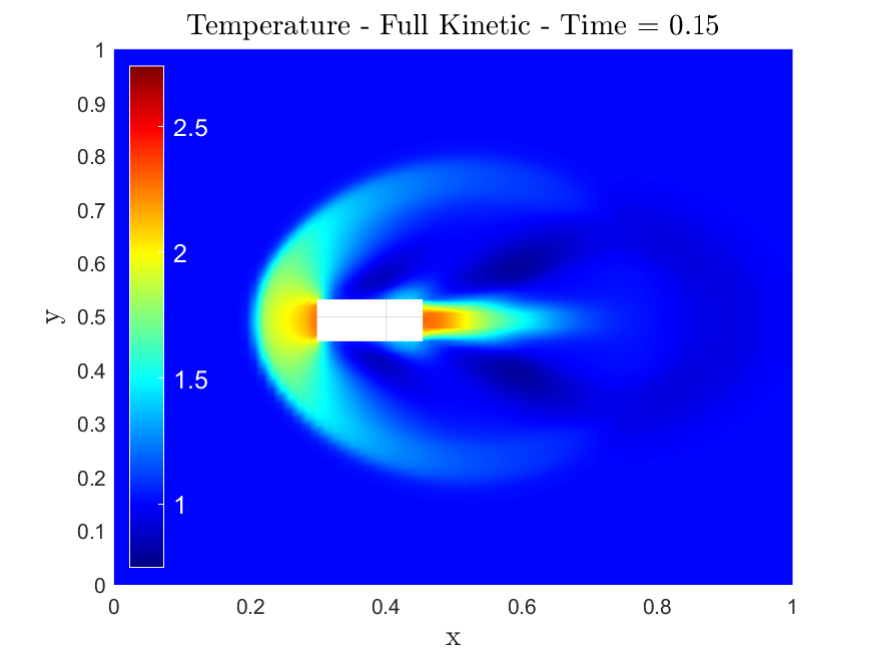}
    \end{subfigure}
    
    \begin{subfigure}{0.32\textwidth}
    \centering
    \includegraphics[width=\textwidth, trim={1cm 0cm 1.8cm 0cm}]{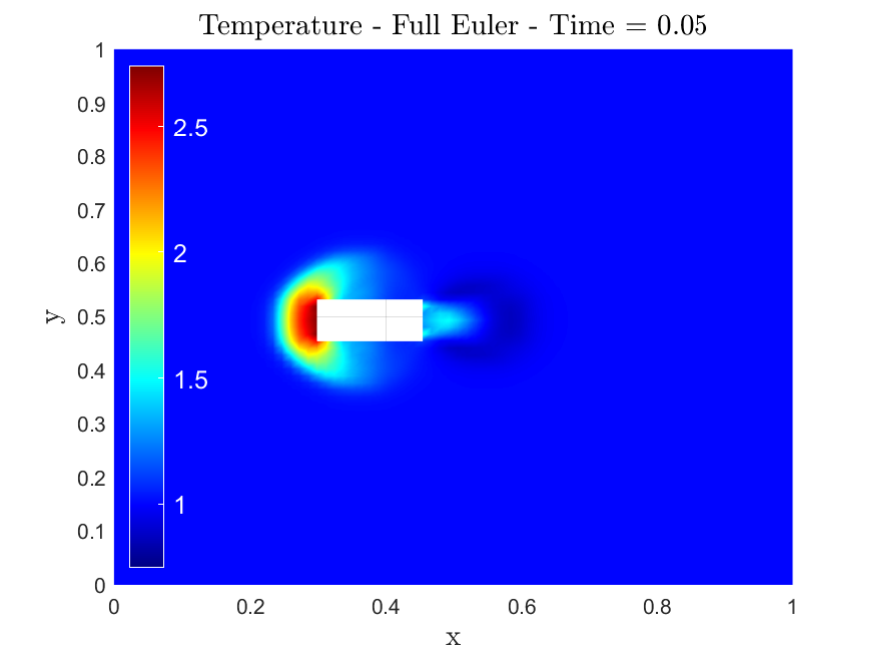}
    \end{subfigure}
    \begin{subfigure}{0.32\textwidth}
    \centering
    \includegraphics[width=\textwidth, trim={1cm 0cm 1.8cm 0cm}]{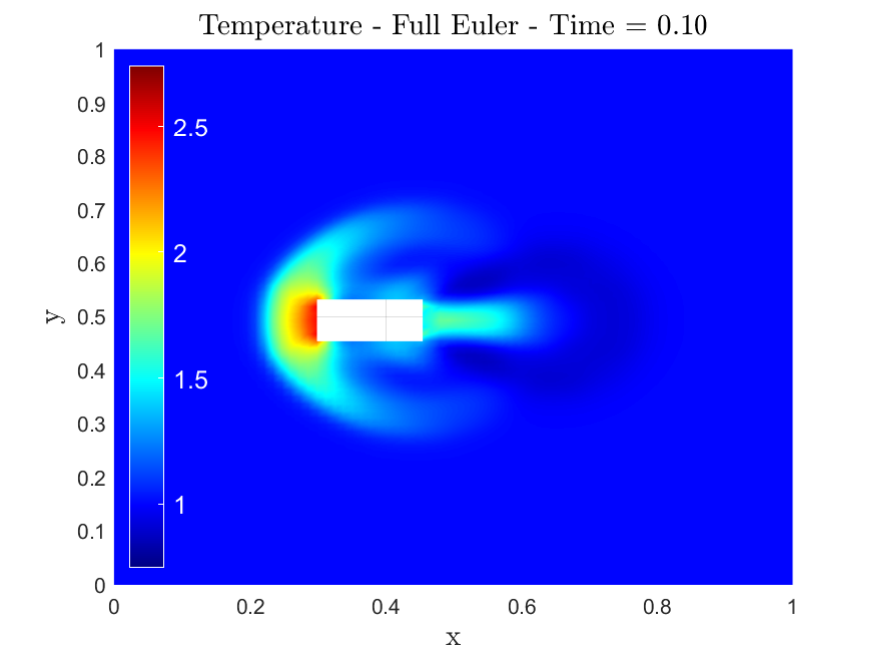}
    \end{subfigure}
    \begin{subfigure}{0.32\textwidth}
    \centering
    \includegraphics[width=\textwidth, trim={1cm 0cm 1.8cm 0cm}]{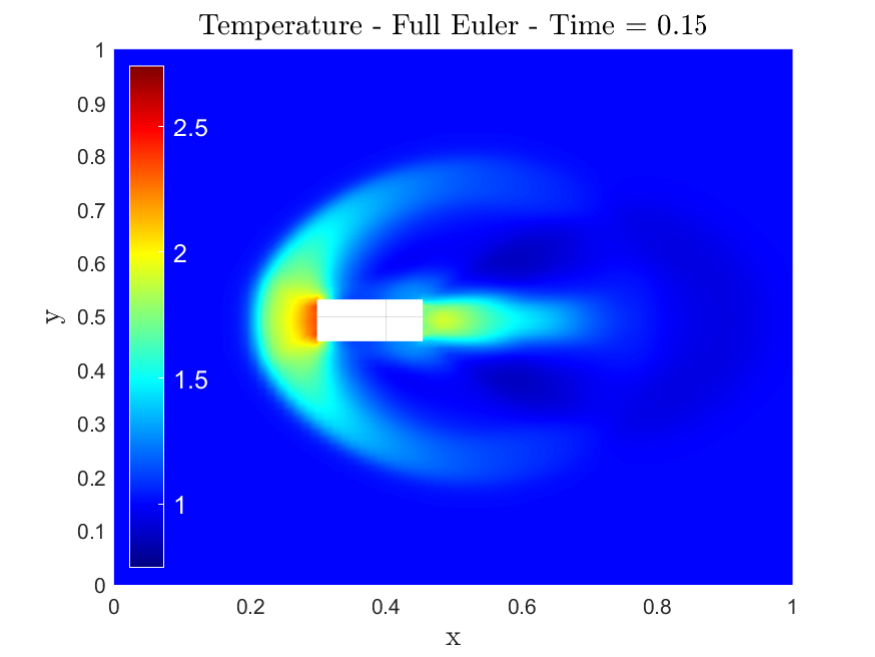}
    \end{subfigure}
    \caption{Temperatures obtained using the hybrid (first row), the full Kinetic (second row) and full Euler (third row) solvers, respectively, for the fluid flux towards a fixed rectangular obstacle with initial condition \eqref{eq::initial_Rectangular}. The solutions are displayed at time $t=0.05$, $t=0.10$ and $t=0.15$.}
    \label{figRec_temperature}
\end{figure}

In Figure \ref{figRec_density} we compare the density obtained using the hybrid scheme, the full Euler solver and the full Kinetic (ES-BGK) solver. In the same Figure on the first row it is also reported the domain adaptation. As already observed for the Section \ref{Test5}, also in this case the solutions are clearly consistent with each other and the kinetic regime is dynamically activated in the vicinity of the discontinuity, generated by the presence of the obstacle. This dynamical update ensures that the computational resources are concentrated only close to the fixed rectangular object, without unnecessary computational cost in the smoother regions. In Figure \ref{figRec_temperature} are reported the same quantities but for the temperature.

\subsection{Analysis of the speedup}
\begin{figure}
    \centering
    \includegraphics[scale=0.7]{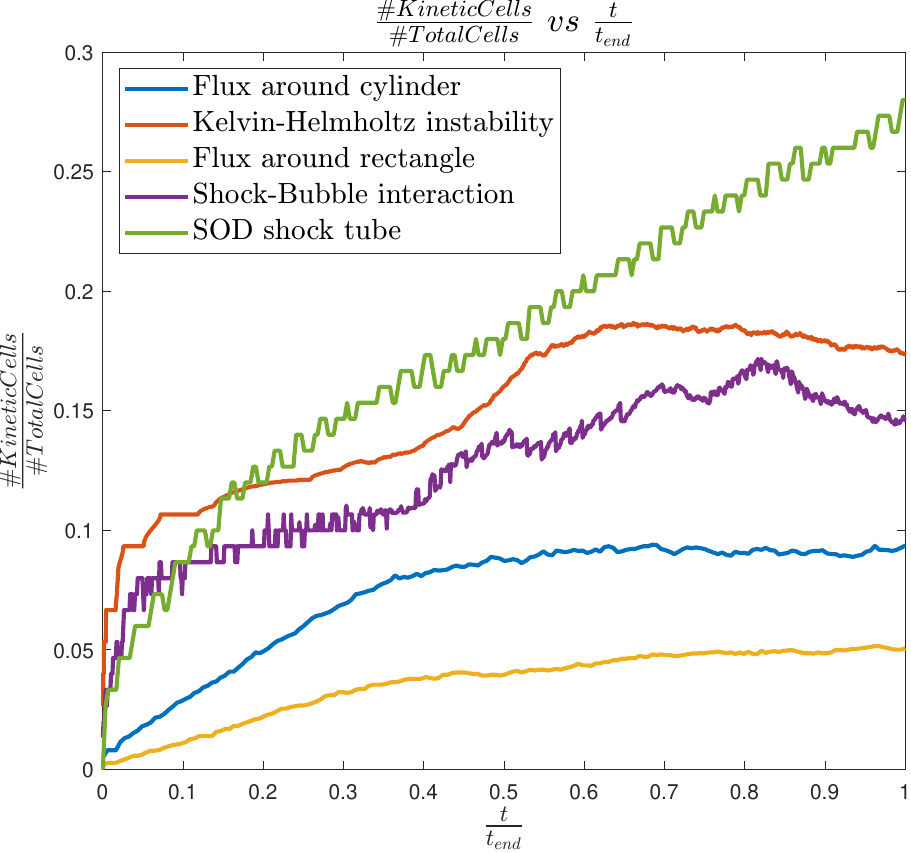}
    \caption{Fraction of kinetic cells over total number of cells for each simulation. It is worth noting that the vertical axis extends up to 0.3.}
    \label{figSpeed}
\end{figure}

In Figure \ref{figSpeed} the evolution in time of the ratio between the number of kinetic cells over the total number of cells is reported, for each simulation. It should be noted that the vertical axis extends up to 0.3, thereby indicating that, throughout the simulations, the number of kinetic cells remains consistently very small over time. Most of the computational domain continues to be treated using the macroscopic model and only a limited number of cells near discontinuities and other non-equilibrium regions require the kinetic solver. This highlights the efficiency of the adaptive approach, ensuring accuracy in capturing the overall dynamics.\\
This synergy between high-order coupled time integration and adaptive domain decomposition makes the method extremely effective for long-time simulations in rarefied gas dynamics. 

\section{Conclusion}
In this paper, we have proposed a high-order domain decomposition method for the ES-BGK model of the Boltzmann equation, which dynamically detects regions of equilibrium (solved using Euler equations) and non-equilibrium (where the ES-BGK model is used). The main novelty of this work is the development of a coupled strategy between the macroscopic and the kinetic solvers, which preserves the overall temporal order of accuracy of the scheme. Our approach is based on a coupled IMEX method across different subdomains and allows high accuracy and computational efficiency. It is important to remark that our strategy is not limited only to the coupled system Euler and ES-BGK, but to any system of this kind. 
Further perspectives include the application of this high-order coupling to different types of equations, and combining such domain decomposition methods with other time integration strategies. 

\section*{Acknowledgments}
The authors would like to thank Gabriella Puppo, Lorenzo Pareschi and Thomas Rey for fruitful discussions and for their comments on this manuscript.\\  
The authors received funding from the European Union's Horizon Europe research and innovation program under the Marie Skłodowska-Curie Doctoral Network DataHyking (Grant No. 101072546).
The authors acknowledges the support of the Mathematics and Computer Science Department of the University of Ferrara through access to their computational resources. TT is member of GNCS-INdAM research group.

\section*{Declarations}
\subsection*{Conflict of interest}
The authors declare that they have no financial or non-financial conflicts of interest that are relevant to the content of this article.

\bibliography{references}
\bibliographystyle{acm}

\end{document}